\documentclass[twoside,10pt,a4paper,german]{article}
\usepackage{german}

     \def \id{\mbox{\rm id }}

\def \non{\nonumber}

\def \qed{\mbox{ }\hfill{\bf q.e.d.}\par}
\def \beq{\begin{equation}}
\def \endeq{\end{equation}}
\def \beqa{\begin{eqnarray}}
\def \eeqa{\end{eqnarray}}
\def \beqanon{\begin{eqnarray*}}
\def \eeqanon{\end{eqnarray*}}
\def \bay{$$ \begin{array}}
\def \eay{\end{array} $$ }

\def \supp{\mbox{\rm supp }}

\def \dist{\mbox{\rm dist }}

\def \rot{\mbox{\rm rot }}
\def \div{\mbox{\rm div }}

\def \arccos{\mbox{\rm arc cos }}

\def \i{{\rm i \,}}

\def \ra{\rightarrow}
\def \Ra{\Rightarrow}

\def \sp#1#2{\langle \, #1 \, , \, #2 \, \rangle}

\def \ol#1{\overline{#1}}

\def \Abb#1#2#3#4#5
 {\begin{array}{ccccc}
 #1 & : & #2 & \longrightarrow & #3 \\
 {}   & {}  & #4 & \longmapsto & #5
    \end{array}}

\def \zweimat#1#2#3#4{\left( \begin{array}{cccc} #1&#2\\#3&#4 \end{array} \right)}
\def \zweivec#1#2{\left( \begin{array}{cc}#1\\#2 \end{array}\right)}
\def \d#1{\frac{d}{dx_#1}}

\def \alle#1{\bigwedge\limits_{#1}}
\def \gibt#1{\bigvee\limits_{#1}}

\def \abstand{\mbox{\  } \vspace{0.5cm} \mbox{\ } \new \noindent }
\def \abstandk{\mbox{\  } \vspace{0.3cm} \mbox{\ } \new \noindent }

\def \eps{\epsilon}
\def \new{\newline}
\def \vp{\varphi}
\def \Lra{\Leftrightarrow}
\def \mit{\ |\ }
\def \achtung#1{   }
\def \formel#1{\label{#1} \mbox{ }  }
\def \hsp{\hspace{-0.1cm}}
\def \zu{\hspace{-0.1cm}}
\def \hspz{\hspace{-0.2cm}}
\def \rrot{{\rm Rot}}
\def \ddiv{{\rm Div}}

\def \omeg#1{\omega_m^{ #1 }}
\def \ommeg#1{(\omega_m^{ #1 })^{-1}}
\def \hnull#1#2{\stackrel{\circ}{H} {}\zu_{ #1 }^{ #2 }}
\def \rrote{\rrot\ }
\def \ddivh{\ddiv\ }

\def \cht{\check{T}}
\def \chcht{\check{t}}
\def \chn{\check{N}}
\def \chgt{\check{\gamma}_T}
\def \chgn{\check{\gamma}_N}
\def \chtau{\check{\tau}}
\def \chrho{\check{\rho}}

\def \calf{{\cal F}}
\def \calL{{\cal L}}
\def \call{{\rm l}}
\def \calm{{\cal M}}
\def \cals{{\cal S}}
\def \calj{{\cal J}}

\def \tils{\tilde{S}}

\def \rrr{{\rm R}}
\def \rrrp#1{\stackrel{\bullet}{\rrr} {}\zu^{ #1 }}

\def \nnn{{\rm N}}
\def \zzz{{\rm Z}}
\def \ccc{{\rm C}}

\def \eh{\frac{1}{2}}

\def \malle{\mbox{ alle }}
\def \mfur{\mbox{ f\"ur }}
\def \mmit{\mbox{ mit }}
\def \beweis{{\bf Beweis: }}
\def \id{{\rm id}}
\def \d{{\rm d}}
\def \dx{{\rm dx}}
\def \dh{{\rm dh}}
\def \dg{{\rm dg}}
\def \dy{{\rm dy}}
\def \dvp{{\rm d\vp}}
\def \dr{{\rm dr}}
\def \gt{\gamma_T}
\def \gn{\gamma_N}

\def \normu#1#2{|| #1 ||_{ #2 }}
\def \spu#1#2#3{< #1 , #2 >_{ #3 }}

\def \ssqn#1{\sum_{I\in{\cal S}(q,N-1)} #1 }
\def \ssqmn#1{\sum_{I\in{\cal S}(q-1,N-1)} #1 }

\def \sqqn{\sum_{I\in{\cal S}(q,N)}}

\def \rq{R^{q,\Gamma_1}}

\def \rqm{R^{q-1,\Gamma_1}}

\def \rqn{\stackrel{o}{R}{}\zu ^q}
\def \rqnn{\stackrel{o}{R}{}\zu_0 ^{q}}

\def \rqmn{\stackrel{o}{R} {}\zu^{q-1}}

\def \rnull#1#2{\stackrel{o}{R}{}\zu_{ #1 }^{ #2 } }

\def \dq{D^{q,\Gamma_2}}
\def \dqp{D^{q+1,\Gamma_2}}

\def \dqn{\stackrel{o}{D} {}\zu^q}

\def \dnull#1#2{\stackrel{\circ}{D} {}\zu_{ #1 }^{ #2 } }

\def \hqeh{H_{1/2}^q}
\def \hqehm{H_{-1/2}^q}

\def \hqpehm{H_{-1/2}^{q+1}}

\def \hqpdh{H_{3/2}^{q+1}}

\def \hvnull#1#2{{H}_{ #1 }^{ #2} }
\def \rvnull#1#2{{R}_{ #1 }^{ #2} }
\def \dvnull#1#2{{D}_{ #1 }^{ #2} }
\def \cnull#1#2{\stackrel{\circ}{C}{}\zu_{ #1 }^{ #2 }}
\def \cqe{C_{\infty}^{q,\Gamma_1}}

\def \cqem{C_{\infty}^{q-1,\Gamma_1}}

\def \cqz{C_{\infty}^{q,\Gamma_2}}
\def \cqzp{C_{\infty}^{q+1,\Gamma_2}}

\def \cqn{\stackrel{\circ}{C}{}\zu_{\infty}^q}
\def \cq{C^q_{\infty}}
\def \cqp{C^{q+1}_{\infty}}

\def \cqn{\stackrel{\circ}{C}{}\zu^q_{\infty}}

\def \cq{C^q_{\infty}}
\def \cqp{C^{q+1}_{\infty}}

\newtheorem{satz}{Satz}[section]
\newtheorem{lemma}[satz]{Lemma}
\newtheorem{korollar}[satz]{Korollar}
\newtheorem{defi}[satz]{Definition}
\newtheorem{bemerk}[satz]{Bemerkung}

\title{\sc Die Maxwellgleichung\\ mit wechselnden Randbedingungen}
\author{\sf Peter Kuhn\footnote{\sf Dies ist ein preprint
meiner Dissertation (Dr.~rer.~nat.), welche dem Fachbereich 6, Mathematik und Informatik,
der Universit{\"a}t Essen, Deutschland, im August 1999 vorgelegt
und bei Shaker (Aachen, Deutschland) im Januar 2000 publiziert wurde.}}

\begin{document}

\date{August 1999}
\maketitle
\thispagestyle{empty}
\newpage
\thispagestyle{empty}
\tableofcontents
\newpage
\setcounter{page}{1}

\section{Einleitung}
\markboth{EINLEITUNG}{INHALT}
\subsection{Inhalt}
Sei $S $ ein
beschr\"anktes Gebiet
im $\rrr^3$, dessen Rand $\partial S$ in die
beiden Komponenten $\Gamma_1$ und $\Gamma_2$
zerlegt ist.
Seien ferner
$\eps$ und  $\mu$ (Dielektrizit\"at und Permeabilit\"at)
gleichm\"a{\ss}ig positiv
definite $3\times 3$ Matrizen und $\vec{N}$
die \"au{\ss}ere Normale.
Wir betrachten das Maxwellsche Randwertproblem
zu gegebenen Feldern
$\vec{J}$ und $\vec{K}$ L\"osungen $\vec{E}$ und $\vec{H}$
der Gleichungen
\beqa
 \rot \vec{E} +\i\omega \mu \vec{H} &=& \vec{J} \non\\
 \rot \vec{H} -\i\omega \eps \vec{E} &=& \vec{K} \non\\
  \vec{N}\wedge \vec{E}&=&0 \mbox{ in } \Gamma_1\non\\
  \vec{N}\cdot \vec{E} &=& 0 \mbox{ in } \Gamma_2\non
\eeqa
zu finden.
Um eine L\"osungstheorie f\"ur die Maxwellschen Gleichungen
aufzubauen, kann man diese in die Theorie  alternierender
Differentialformen beliebigen Ranges in beliebigen
Raumdimensionen einbetten: Identifizieren
wir die Felder
$\vec{E},\vec{K},\vec{N}$ mit den 1--Formen
$E:=\vec{E}\cdot {\rm \vec{ds}}$,
$K:=-\vec{K}\cdot {\rm \vec{ds}}$,
$N:=\vec{N}\cdot {\rm \vec{ds}}$
und $\vec{H},\vec{J}$
mit den $2$--Formen
$H:=\vec{H}\cdot {\rm \vec{dF}}$,
$J:=\vec{J}\cdot {\rm \vec{dF}}$, wobei
\beqa
 {\rm \vec{ds}} :=
\left[  \begin{array}{c}
\dx^1\\
\dx^2\\
\dx^3
\end{array} \right] \ ,
  \ {\rm \vec{dF}} :=
\left[  \begin{array}{c}
\dx^2\wedge\dx^3\\
\dx^3\wedge\dx^1\\
\dx^1\wedge\dx^2
\end{array} \right] \ ,\non
\eeqa
so erf\"ullen diese
\beqa \left.\begin{array}{rcl}
 \d E+\i\omega\mu H&=&J\\
 \delta H+\i\omega\eps E&=&K \\
  N\wedge E &=& 0\mbox{ in }\Gamma_1\\
  N\wedge * E &=& 0\mbox{ in }\Gamma_2
  \end{array}\right\} \formel{main}
\eeqa
(positiv definite Matrizen $\sigma$ werden durch die Vorschrift
$\sigma( \vec{E}\cdot\vec{ds})=(\sigma\vec{E})\cdot\vec{ds}$
zu positiv definiten Transformationen von
Differentialformen). Die \"au{\ss}ere Ableitung
$\d$ bezeichnen wir in Zukunft mit rot, die Koableitung
$\delta$ mit div.
Durch eine geeignete schwache Formulierung
lassen sich die Gleichungen in
$(\ref{main})$ mit Hilbert\-raum\-me\-tho\-den
behandeln: \abstandk
Wir schreiben $L_2^q(S)$ f\"ur die Menge der
$q$--Formen,
deren
Komponentenfunktionen quadratintegrabel sind,
$L_{2,\sigma}^q(S):=\sigma^{-1/2}L_{2}^q(S)$
und
$R^{q}(S)$ bzw. $D^q(S)$ f\"ur die Menge der Formen aus $L_2^q(S)$,
deren Rotation bzw. Divergenz Elemente aus $L_2^{q+1}(S)$ bzw.
$L_2^{q-1}(S)$ sind.
Die Verallgemeinerung der Randbedingung
$N\wedge E=0$ in $\Gamma_1$ bzw. $N\wedge *E=0$ in $\Gamma_2$
kennzeichnen wir mit einem oberen Index: $R^{q,\Gamma_1}(S)$
bzw.
$D^{q,\Gamma_2}(S)$.
Geringe
Voraussetzungen an die Trennmenge
$\ol{\Gamma}_1\cap\ol{\Gamma}_2$
(vgl. Satz \ref{dichtesatz})
liefern
einen selbstadjungierten Maxwelloperator
\beqa
  D(M)&:=&R^{q,\Gamma_1}(S)
  \times
  D^{q+1,\Gamma_2}(S)\subset L_{2,\eps}^q(S) \times L_{2,\mu}^q(S)
   \non\\
  M&:=&\zweimat{0}{\i\eps^{-1}\div}{\i\mu^{-1}\rot}{0} \ ,\non
\eeqa
und die Gleichungen aus $(\ref{main})$ gehen \"uber  in
\beqa
M\zweivec{E}{H}-\omega\zweivec{E}{H}=
\zweivec{\i\eps^{-1}K}{\i\mu^{-1}J} \ .\formel{main2}
\eeqa
Der Maxwelloperator wird
vom Raum
\beqa
  R^{q,\Gamma_1}(S) \cap
 \ol{\eps^{-1}\div D^{q+1,\Gamma_2}(S)} \times
  D^{q+1,\Gamma_2}(S) \cap
 \ol{\mu^{-1}\rot R^{q,\Gamma_1}(S)} \formel{reduziert}
\eeqa
und dessen orthogonalem Komplement
\beqa
  R_0^{q,\Gamma_1}(S)\times D_0^{q+1,\Gamma_2}(S)
   \formel{reduziert2}
\eeqa
reduziert
(der untere Index 0 steht f\"ur Rotations-- bzw.
Divergenzfreiheit).
Da die Behandlung der durch $(\ref{reduziert2})$ reduzierten
Gleichung evident ist, ist vorwiegend
der durch
$(\ref{reduziert})$
reduzierte Maxwelloperator Gegenstand
unserer Betrachtungen.
F\"ur einen glatten Rand $\partial S$ und eine
glatte Trennmenge werden wir zeigen, da{\ss}
die Einbettung
\beqa
 R^{q,\Gamma_1}(S)\cap D^{q,\Gamma_2}(S) \hookrightarrow L_2^q(S)
 \formel{einlkompakt}
\eeqa
kompakt ist. Dieses Ergebnis hat als Konsequenz, da{\ss}
die R\"aume $\rot R^{q,\Gamma_1}(S)$, $\div D^{q,\Gamma_2}(S) $
abgeschlossen und die Dirichlet--Neumann--Felder
\beqa
  R_0^{q,\Gamma_1}(S)\cap D_0^{q,\Gamma_2} (S)\formel{einldirineu}
\eeqa
endlich dimensional sind.
Das Spektrum des Maxwelloperators besteht dann
nur aus isolierten Punkten in $\rrr$.
Dar\"uber hinaus existiert
ein kompakter L\"osungsoperator,
so da{\ss} f\"ur die Gleichung
$(\ref{main2})$ die Fredholmsche
Alternative gilt (vgl. \cite{weck}). \new
\abstandk
Wir werden
die Tangentialspuren von Formen aus $R^q(S)$ untersuchen.
Schon bekannt ist, da{\ss} im Falle eines glatten
Gebietes $S$ ein linearer, stetiger und
surjektiver Spuroperator
von $R^q(S)$ nach $R_{-1/2}^q(\partial S)$, die Menge
der Funktionale auf
$H_{1/2}^{q} (\partial S)$,
deren Tangentialrotation
Funktionale auf $H_{1/2}^{q+1}(\partial S)$ sind, existiert.
Hier\-f\"ur liefern wir einen weiteren Beweis.
Wir werden f\"ur eine glatte
Trennmenge die Existenz eines linearen stetigen Spuroperators
$R^q(S)\ra R_{-1/2}^{q,\Gamma_2}(\partial S)$ zeigen.
Letzter Raum besteht aus den Einschr\"ankungen der Funktionale
aus $R_{-1/2}^{q}(\partial S)$ auf
$H_{1/2}^{q,\Gamma_2}(\partial S)$, Elemente aus
$H_{1/2}^{q}(\partial S)$, die fast \"uberall in $\Gamma_2$
verschwinden.
Mit einer \"ahnlichen Technik l\"osen wir dann das
statische Maxwellsche Problem
\beqa
  \rot E=F, \div E=G, N\wedge E =\lambda \ .\non
\eeqa
Dualit\"at liefert entsprechende Resultate f\"ur
den Raum $D^q(S)$. \new
Im Falle einer leeren Trennmenge lassen sich aus
den obigen Aussagen
ein linearer stetiger
Fortsetzungsoperator $R_{-1/2}^{q}(\Gamma_1)\ra R^q(S)$
sowie eine
L\"osungstheorie
f\"ur
das Problem
\beqa
 \rot E=F,\div E=G,N\wedge E=\lambda \mbox{ in } \Gamma_1,
 N\wedge *E =\theta \mbox{ in } \Gamma_2\formel{maxmisch}
\eeqa
ableiten.
\abstandk
Abschlie{\ss}end betrachten wir zwei weitere
Probleme im Zusammenhang mit wechselnden Randbedingungen.
Zum einen zeigen wir
f\"ur ein glattes Gebiet im $\rrr^N$, dessen
Randst\"uck $\Gamma_1$ in $K$
glatte Zusammenhangskomponenten zerf\"allt, da{\ss}
die Dimension des Raumes
$(\ref{einldirineu})$ f\"ur $q=1$ gerade
$K-1$ ist. Zum anderen wollen wir
die Eigenformen des Maxwelloperators
auf der halben Kreislinie und dem Halbkreis
bestimmen und Aussagen \"uber deren
Regularit\"atseigenschaften machen.
\subsection{Geschichte}
\markright{GESCHICHTE}
Die Verallgemeinerung der Maxwellgleichungen auf
Differentialformen geht auf Weyl \cite{weyl} zur\"uck,
der mit Integralgleichungsmethoden eine
L\"osungstheorie f\"ur den homogenen isotropen Fall
($\eps=\mu=1$) aufstellen konnte.
\abstandk
In glatten Gebieten folgt die kompakte Einbettung
\beqa
  R^{q,\partial S}(S) \cap  D^q(S)\hookrightarrow L_2^q(S)\formel{einlkomp}
\eeqa
mit dem Rellichschen Auswahlsatz
aus der stetigen Einbettung
\beqa
 R^{q,\partial S }(S) \cap  D^q(S)\hookrightarrow H_1^q(S) \ .\formel{einlreg}
\eeqa
Einen solchen Regularit\"atsbeweis lieferte
Leis in \cite{leis} f\"ur glatte Gebiete im $\rrr^3$.\new
F\"ur die auf Differentialformen verallgemeinerten
Maxwellgleichungen konnte
Weck  in \cite{weck2}, \cite{weck}
eine gro{\ss}e Klasse von nicht glatten Gebieten
(verallgemeinerte Kegelgebiete) angeben,
in denen die Einbettung $(\ref{einlkomp})$
kompakt ist.
Er behandelte den inhomogenen anisotropen Fall
($\eps$ und $\mu$ geeignete
Transformationen von Differentialformen) mittels einer
vollst\"andigen Induktion \"uber die Raumdimension.
Dabei zeigte er auch die Unabh\"angigkeit der
kompakten Einbettung von $\eps$ und $\mu$.
\new
Ein Beweis f\"ur Gebiete im $\rrr^3$ mit der
eingeschr\"ankten Kegeleigenschaft wurde von
Weber in \cite{weber2} gef\"uhrt.
Die Voraussetzungen an das Gebiet wurden
lediglich f\"ur die Existenz
eines Calderonschen Fortsetzungsoperators
$H_2(S)$ nach $H_2(\rrr^3)$
ben\"otigt, um dann (auf Felder aus $\nabla H_2$ ) den
Rellichschen Auswahlsatz anzuwenden.\new
Witsch ersetzte in \cite{witsch}
diese Kombination, Fortsetzungsoperator und kompakte Einbettung,
durch einen kompakten Fortsetzungsoperator
$H_2(S)$ nach $H_1(\rrr^N)$, f\"ur dessen Existenz
er die Voraussetzungen in \cite{weber2} weiter
zu $p$--$cusp$ Gebieten mit $p<2$ abschw\"achen konnte.\new
Einen elementaren Beweis brachte Picard in
\cite{picard3} im Fall der Weylschen Verallgemeinerung
f\"ur Lipschitz--Gebiete,
eine gr\"o{\ss}ere Menge, als die der
Gebiete mit der eingeschr\"ankten Kegeleigenschaft.
Nachdem er die Unabh\"angigkeit der kompakten Einbettung
von Lipschitz--Transformationen gezeigt hatte,
lokalisierte er das Problem
und konnte es dann auf die Einheitskugel \"ubertragen.
Dort f\"uhrten schon bekannte Resultate
zum Ziel.\new
Eine Vereinfachung des Beweises aus \cite{weck} wurde
in \cite{wwp} f\"ur $S\subset\rrr^3$ gef\"uhrt.
Dar\"uber hinaus wurde
durch einen anderen Induktionsanfang die
Klasse der Gebiete  mit
kompakter Einbettung $(\ref{einlkomp})$ nochmals erweitert.
Dies f\"uhrte unter anderem zu Teilmengen aus $\rrr^3$,
die lokal Lipschitz--hom\"oomorph zu
Gebieten sind, die aus endlich vielen Zusammenhangskomponenten
von $p$--$cusp$'s, $p<2$ oder Kegelspitzen bestehen.
\abstandk
Spuroperatoren $R^q(S)\ra R^{q}_{-1/2}(\partial S)$,
Fortsetzungss\"atze und L\"osungstheorien f\"ur
das statische Problem
\beqa
 \rot E=F\ ,\ \div E =G\mbox{ in } S \ ,\
 N\wedge E =f \mbox{ in }  \partial S\non
\eeqa
wurden von Georgescu in
\cite{georgescu} und von Paquet in \cite{paquet}
f\"ur Differentialformen auf kompakten
glattberandeten Mannigfaltigkeiten
untersucht.
\new
Alonso und Valli fanden in
\cite{alonso}  einen Weg, den Fortsetzungsoperator
f\"ur Gebiete im $\rrr^3$
durch L\"osen geeigneter Differentialgleichungen
herzuleiten.
Im Falle einer leeren Trennmenge
charakterisierten sie die Tangentialspuren
der Felder aus $R^q(S)$ auf einem  Randst\"uck $\Gamma_1$
und brachten eine L\"osungstheorie f\"ur das
Problem $(\ref{maxmisch})$.\abstandk
Weitere Untersuchungen des statischen Problems
findet man in
\cite{kresspotential} und
\cite{picard4}.
Die hierbei auftretenden harmonischen Felder
wurden in
\cite{duff1}, \cite{duff2} und \cite{martensen} mit
klassischen Methoden
behandelt.
F\"ur nicht glatte Gebiete $S$
hat Picard in \cite{picard2}, \cite{picardhodge} gezeigt, da{\ss} sich
die Dimensionen
der harmonischen Differentialformen oder Neumann--Felder
$(\ref{einldirineu})$ im Falle
$\Gamma_2=\partial S$, $\Gamma_1=\emptyset$ und
$q$ beliebig
durch die Betti--Zahlen des Gebietes ausdr\"ucken lassen;
genauer
\beqa
 \dim (R_0^q(S)\cap D_0^{q,\partial S}(S)) =\beta_q \ ,\non
\eeqa
wobei $\beta_q$ gerade die $q$--te Betti Zahl ist.
\new
Den Fall gemischter Randbedingungen
und leerer Trennmenge
betrachtete Kress in \cite{kress} f\"ur Gebiete im $\rrr^3$.
\abstandk
Saranen untersuchte in \cite{saranen} die G\"ute der
L\"osungen der Maxwellgleichungen in Kegelspitzen,
indem er nach den Eigenformen auf dem Kegeldeckel
entwickelte und die Koeffizienten untersuchte. Hierbei
benutzte er die Resultate aus \cite{weck}.
\subsection{Vorgehensweise}
\markright{VORGEHENSWEISE}
In Kapitel \ref{chvorbereitung} werden wir
zun\"achst
grundlegende Bezeichnungen einf\"uhren
und einige Werkzeuge f\"ur deren Anwendung bereitstellen.
Eine besondere Bedeutung kommt hier den
Approximationseigenschaften zu.
W\"ahrend im Falle homogener Randbedingungen die
Segmenteigenschaft gen\"ugt, um die R\"aume $R^q(S)$ und $D^q(S)$
durch glatte Formen anzun\"ahern, m\"ussen wir im Falle
wechselnder Randbedingungen zus\"atzlich \"ahnliche
Voraussetzungen
an die Randst\"ucke $\Gamma_1$ oder $\Gamma_2$
stellen (Satz \ref{dichtesatz}).
Dies ist notwendig, um sp\"ater
auf den Satz von Stokes in Lemma \ref{walze}
zugreifen zu k\"onnen.
Dieses Lemma liefert ein zweites wichtiges Werkzeug:
Mit Hilfe von  \cite{wwspherical}
stellen wir hier Formeln zur Verf\"ugung, die den
Zusammenhang zwischen der Rotation auf dem Rand und
der Rotation im Inneren spezieller Gebiete, den Kegelspitzen,
darlegen.
\abstandk
Nach diesen Vorbereitungen zeigen wir in Kapitel \ref{kompeinb},
Satz \ref{einistkompakt}
f\"ur einen glatten
Rand $\partial S$ und eine
glatte Trennmenge
die Kompaktheit der
Einbettung in $(\ref{einlkompakt})$.
Den Beweis, f\"uhren wir
wie in \cite{wwp} (vgl. auch \cite{weck} und \cite{weck2})
per Induktion \"uber die
Raumdimension:
Gilt die kompakte Einbettung, so k\"onnen wir
nach Eigenformen entwickeln. Aus
diesem Entwicklungsresultat in $(N-1)$--dimensionalen
Mannigfaltigkeiten folgt schlie{\ss}lich
die kompakte Einbettung in
$N$--dimensionalen Gebieten mit glattem Rand und glatter
Trennmenge.
Bei diesem Dimensionssprung kommen uns die
oben erw\"ahnten Hilfsmittel zugute.
\abstandk
In Kapitel \ref{chregularitaet} werden wir
einen Regularit\"atssatz f\"ur Formen herleiten.
Hier halten wir uns im wesentlichen an den
Beweis aus \cite{weber} und modifizieren diesen
an den Stellen, an denen von der speziellen
Situation im $\rrr^3$ Gebrauch gemacht wird.
Dazu benutzen wir eine Spiegelungstechnik
wie in \cite{wentzig}.
\abstandk
Mit Hilfe dieses Regularit\"atsresultates
werden wir in den folgenden beiden Kapiteln
S\"atze \"uber Spuren, Fortsetzungen
(Kapitel \ref{chspursatz}) und
L\"osungstheorie zum statischen Maxwellproblem
mit homogenen Randbedingungen
(Kapitel \ref{chrwp}) beweisen.
Die hier angewandte Technik
basiert auf der Formulierung geeigneter koerzitiver
Hil\-bert\-raumprobleme und geht
auf \cite{alonso} zur\"uck.
\abstandk
Um in Kapitel \ref{chfelder} f\"ur
$q=1$ die Dimension, der im Falle wechselnder
Randbedingungen
auftretenden Dirichlet--Neumann--Felder $(\ref{einldirineu})$
zu bestimmen,
verallgemeinern wir die Methode aus \cite{picard4}.
\abstandk
In Kapitel \ref{cheigenformen} berechnen wir zun\"achst
f\"ur die halbe Kreislinie die Eigenformen
des Maxwelloperators. Mit
Hilfe des Entwicklungsresultates
aus Kapitel \ref{kompeinb} k\"onnen
wir die Eigenformen f\"ur den Halbkreis
nach diesen entwickeln. Die Koeffizienten dieser Entwicklung
erf\"ullen dann die Besselsche Differentialgleichung,
\"uber deren L\"osungen, die Besselfunktionen, viele
Arbeiten verfa{\ss}t wurden.
Im Falle homogener Randbedingungen
wurde ein \"ahnliches  Verfahren
in \cite{saranen} angewandt.
\section{Vorbereitung}
\label{chvorbereitung}
\markboth{VORBEREITUNG}{BEZEICHNUNGEN}
\subsection{Bezeichnungen}
Mit $\rrr$\label{pxrrr} bzw. $\ccc$\label{pxccc} bezeichnen wir
die Menge
der reellen bzw. komplexen
Zahlen, mit $\nnn$\label{pxnnn} die Menge der
nat\"urlichen Zahlen ohne $0$. F\"ur komplexe Zahlen $z$
ist $\ol{z}$\label{pxolz} die Konjugation. Die imagin\"are Einheit
nennen wir $\i$\label{pxi}.
Falls $U$ eine Menge und $n$ eine nat\"urliche Zahl ist,
definieren
wir rekursiv
$U^1:=U$, $U^n:=U^{n-1}\times U$\label{pxuhochn}.
F\"ur die Normen in $\rrr^{N}$ und $\ccc^N$ schreiben wir
$|\cdot|$\label{pxnormr}.
\abstandk
Sind $U,V$ Teilmengen eines metrischen Raumes
$(X,d)$,
so ist $\ol{U}$ \label{pxolu} der Abschlu{\ss} und
$\partial U$\label{pxupartial} der Rand von $U$. Ist nicht klar,
bez\"uglich welcher Metrik der Abschlu{\ss} zu bilden
ist, versehen wir $\ol{U}$ mit einem
oberen Index $d$.
Wir sagen $U\subset\subset V$\label{pxkompakt}, wenn
$\ol{U}$ kompakt und $\ol{U}\subset V$ gilt.
Die Abstandsfunktion bezeichnen wir mit dist\label{pxdist} und setzen
$\dist(U,\emptyset):=\infty$.
\abstandk
F\"ur zwei Mengen $U,V$ ist $\calf(U,V)$\label{calf} die Menge aller
Abbildungen $f$,
deren Definitionsbereich $D(f):=U$\label{pxdefinition} ist,
und deren Wertebereich
$R(f)$\label{pxwerte} in $V$ liegt; $N(f)$ ist der Nullraum.
Gilt $U'\subset U$, so
bezeichnen wir mit  $f_{|_{U'}}$\label{pxeinsch} die Einschr\"ankung
von $f$ auf $U'$.
Wir schreiben $\supp f$\label{pxsupp} f\"ur den Tr\"ager einer
komplexwertigen Funktion $f$.
F\"ur $G\subset \rrr^N$, $G$ offen definieren wir weiter\new
\begin{tabular}{rcp{10.0cm}}
  $C_\infty(G)$\label{pxcunendlich}&& Raum der unendlich oft differenzierbaren
    komplexwertigen Funktionen\\
  $\cnull{\infty}{\ } (G)$\label{pxcnullunendlich}&:=&$\{\vp\in C_\infty(G) \mit
    \supp\vp\subset\subset G\}$\\
  $C_\infty(\ol{G})$\label{pxcunendol}&:=&$ \{\vp_{|_G}\mmit
  \vp\in\cnull{\infty}{\ }(\rrr^N) \}$ \\
  $L_p(G)$\label{pxlzwei}&& Raum der \"Aquivalenz\-klassen aller
   Lebes\-gue--me{\ss}\-baren Funk\-tionen $f$ mit
     $\normu{f}{L_p(G)} :=(\int_{G}|f(x)|^p\dx )^{1/p}<\infty $,
   $p=1,2$ \\
 $\spu{f}{g}{L_2(G)}$ &:=&$\int_Gf(x)\ol{g(x)}\dx$\\
 $H_m(G)$\label{pxhm}&& Sobolevr\"aume (siehe \cite[Definition 3.1]{wloka})
 mit Norm \new $\normu{\cdot}{H_m(G)} $.
\end{tabular}
\abstand
Wir schreiben $H_1\oplus H_2$\label{pxoplus} f\"ur die orthogonale
Summe zweier Unterr\"aume $H_1,H_2$
eines Hilbert\-raumes.
Der zu einem linearen Operator $A$ adjungierte
Operator ist $A^*$\label{pxadjungierter}.
Haben wir in einem
Raum $H$ ein Skalarprodukt $\spu{\cdot}{\cdot}{H} $
erkl\"art, so setzen wir
$\normu{x}{H} := (\spu{x}{x}{H} )^{1/2}$ f\"ur $x\in H$.
Die Dimension eines Vektorraumes $V$ ist $\dim V$\label{pxdimension}.
\abstandk
Mit $c$\label{pxc} bezeichnen wir Konstanten, die sich im Laufe eines
Beweises \"andern k\"onnen, deren
\"Anderungen aber unabh\"angig von
aus dem Kontext ersichtlichen Eigenschaften sind.
\abstandk
Elemente
$(i_1,\cdots,i_q)\in\{1,\cdots,N\}^q$ mit $i_k\not=i_l$
f\"ur $k\not=l$
nennen wir Multiindizes der L\"ange $|I|:=q$.\label{pxibetrag}
Gilt $i_1<\cdots<i_q$ f\"ur
einen Multiindex $I=(i_1,\cdots,i_q)$, so ist dieser geordnet.
Die
Menge der geordneten Multiindizes
bezeichnen wir mit $\cals(q,N)$\label{pxcals} und
sagen $j\in I$, falls $j\in\{i_1,\cdots,i_q\}=:{\cal I}$,
$j\not\in I$ entsprechend. Im ersten Fall
ist $I-j=:(\hat{i}_1,\cdots \hat{i}_{q-1})$
der geordnete Multiindex der
L\"ange $q-1$ mit $\{\hat{i}_1,\cdots \hat{i}_{q-1},j \}={\cal I}$,
im zweiten Fall
schreiben wir $I+j=:(\hat{i}_1,\cdots \hat{i}_{q+1})$
f\"ur den geordneten Multiindex der
L\"ange $q+1$ mit
$\{\hat{i}_1,\cdots \hat{i}_{q+1}\}\setminus \{j\}={\cal I}$.
F\"ur einen
ungeordneten Multiindex $I$
schreiben wir
$\sigma(I)$\label{pxsigma} f\"ur das Vorzeichen der Permutation,
welche die Ordnung wiederherstellt.
Somit
gilt f\"ur Multiindizes $I$ der L\"ange $p$ und $J$ der
L\"ange $q$
\beqa
 \sigma(I,J)=(-1)^{pq}\sigma(J,I)\ ,\non
\eeqa
wobei
\beqa
 I,J&:=& (I,J):=(i_1,\cdots,i_p,j_1,\cdots,j_q) \non\\
 \mfur && I=(i_1,\cdots,i_p)\ ,
  \ J=(j_1,\cdots,j_q) \non
\eeqa
f\"ur die Konkatenation von $I$ und $J$ steht.
Wir bezeichnen mit $\calj$\label{pxcalj} die Abbildung,
die einen ungeordneten Multiindex sortiert und
mit $I'$\label{pxistrich} den Multiindex, der die Gleichung
$\calj(I,I')=(1,\cdots,N)$ erf\"ullt.
Nat\"urliche Zahlen wollen wir mit Multiindizes
der L\"ange 1 identifizieren.
\abstand
F\"uhren wir die Vorzeichen der
''Divergenz'' (siehe $(\ref{dddiv})$ und Seite \pageref{pagedelta})
und von ''$**$''  in den Raumdimensionen
$N$ und $N-1$ ein durch\label{pxsigmaq}
\beqa  \begin{array}{rclcrclc}
 \sigma_{q}&:=&(-1)^{N(q-1)}&,&\sigma'_{q}&:=&(-1)^{(N-1)(q-1)}&,\\
 \kappa_{q}&:=&(-1)^{q(N-q)}&,&\kappa'_{q}&:=&(-1)^{q(N-1-q)} &,
 \end{array} \non
\eeqa
so gelten
\beqa
 \begin{array}{rclcrclcrclc}
 \kappa_{q+2}&=&\kappa_{q}&,&\sigma_{q+2}&=&\sigma_q&,&
 \kappa_{q}&=&\kappa_{N-q}& ,\\
 \sigma_{N-q}&=&\sigma_{q+1}&,&
 \kappa_{q}\sigma_{q+1}&=&(-1)^q&,&\sigma_q\sigma_{q+1}&=&(-1)^N&,\\
 \sigma_q\kappa_{q}&=&(-1)^{N+q}&,&&&&&&&&
  \end{array}           \non
\eeqa
$\kappa',\sigma'$ entsprechend, und
\beqa
  \kappa'_{q-1}\sigma_{q}=1\ ,
  \ \sigma'_{q}\kappa_q=(-1)^{N+1} \ .\non
\eeqa
Das Kronecker--Symbol bezeichnen wir mit
$\delta_{i,j}$\label{pxkroneck}.
\subsection{Differentialformen}
\markright{DIFFERENTIALFORMEN}
\label{chdifferentialformen}
Sei in dieser Arbeit stets $M$ eine vollst\"andige
$N$--dimensionale
reelle differenzierbare
Mannigfaltigkeit,
versehen mit einer Orientierung und Riemannscher Metrik,
kurz Mannigfaltigkeit,
und $S$ eine offene Teilmenge mit
kompaktem Abschlu{\ss} in $M$.
Die folgenden Aussagen entnehmen wir \cite{bishop}
oder \cite{jaenich}.
\abstandk
Aus den Voraussetzungen an
$M$ folgt die Existenz einer Metrik $d_M$ auf $M$\label{pxdm}.
Wir nennen Paare $(V,h)$\label{pxkarte} Karten um $x$ in $M$
oder Koordinatenumgebung um $x$,
wenn $V$ eine offene Umgebung von $x$ in $M$ ist,
und die Abbildung $h$ diese Umgebung diffeomorph
auf eine offene Teilmenge des $\rrr^N$ abbildet.
Wir treffen folgende Konvention:
\begin{enumerate}
\item[\ ]
Diffeomorphismen $V\ra U$ sind Einschr\"ankungen
invertierbarer Abbildungen $V_0\ra U_0$ mit $V_0,U_0$ offen,
$V\subset\subset V_0$, $U\subset\subset U_0$,
die in beiden Richtungen
unendlich oft differenzierbar sind.
\end{enumerate}
Die Tangenten in einem Punkt $x$ an $M$, die wir
als
Derivationen auf
$C_\infty(m)$ (das ist die Menge der reellwertigen Funktionen $f$,
die in einer
Umgebung $U:=U(f)\subset M$
von $m$ definiert und unendlich oft differenzierbar sind )
auffassen k\"onnen, spannen einen $N$--dimensionalen
linearen Raum ${\cal T}_x$\label{pxtan1} oder ${\cal T}M_x$\label{tanm1} auf.
F\"ur $x\in M$ bezeichnen wir
den komplexen Raum der alternierenden kovarianten
Tensoren vom Rang $q$ zum Tangentialraum von $x$
mit
$A^q(x)$\label{pxaq} und dessen
B\"undel mit $A^q(M)$.
Elemente aus $A^q(M)$ nennen wir $q$--Formen oder Formen.
Ist $q<0$ oder $q>N$, so
identifizieren wir
solche mit der Nullabbildung.
Auf dem Raum $A^q(M)$ ist ein
\"au{\ss}eres Produkt
$\wedge:A^q(M)\times A^p(M)\ra A^{q+p}(M)$\label{pxwedge}
(punktweise) erkl\"art mit der Eigenschaft
\beqa
 \alle{\Phi\in A^q(M)}\alle{\Psi\in A^p(M)}
  \Phi\wedge \Psi=(-1)^{qp} \Psi\wedge\Phi \ .\non
\eeqa
In einer Koordinatenumgebung $(V,h)$ um $x$ bilden
die Differentiale $\dh^{i}$\label{pxdxi} der Koordinatenfunktionen $h_i$
eine Basis von $A^1(x)$, damit auch von $A^1(S\cap V)$, und wir
k\"onnen in $S\cap V$ eine Form $\Phi$
eindeutig darstellen
durch
\beqa
 \Phi=\sum_{I\in{\cal S}(q,N)} \Phi_I\dh^I \formel{philokal}
\eeqa
mit $\Phi_I:V\ra \ccc$
und $\dh^{I}:=\dh^{i_1}\wedge\cdots\wedge \dh^{i_q}$
f\"ur $I:=(i_1,\cdots,i_q)$.
Wegen der Anforderungen an $M$ ist
$A^1(x)$ mit einer Orientierung und  Bilinearform
versehen. F\"ur eine positiv orientierte
Orthonormalbasis
$\{\dh^i\mit i=1,\cdots, N\}$ erkl\"aren wir
punktweise den Sternoperator mittels
\beqa
 *\dh^{I}=\sigma(I,I') \dh^{I'} \ .\formel{sterndef}
\eeqa
Dieser ist unabh\"angig von der Wahl der Karten und
liefert einen Isomorphismus
$*:A^{q}(S)\ra A^{N-q}(S)$ mit den Eigenschaften
\beqa
   **\Phi&=&\kappa_q\Phi \non \\
   \Phi\wedge \Psi&=&*\Phi\wedge*\Psi \non\\
  \ *(\vp \Phi)&=&\vp * \Phi \non
\eeqa
f\"ur $\Phi\in A^q(S)$, $\Psi\in A^{N-q}(S)$ und
$\vp\in A^0(S)$.\abstandk
Wir sagen $\vp\in C_m(S)$\label{pxcm},
$m\in\nnn\cup\{0\}\cup\{\infty \}$, wenn f\"ur
eine Karte, dann alle Karten
$(V,h)$ die Funktionen $\vp\circ h^{-1}$ in
$C_m(h(S\cap V))$ liegen.
Sind die Komponentenfunktionen einer Form $\Phi$
in der Darstellung $(\ref{philokal})$
aus $C_m(S)$, so schreiben wir
$\Phi\in C_m^q(S)$\label{pxcq} und definieren
weiter\label{pxcqnull}\label{pxcqol}
\beqa
   \cnull{m}{q} (S)&:=&\{\Phi\in C_m^q(S) \mit
     \supp\Phi\subset\subset S\} \non\\
    C_m^q(\ol{S}) &:=&\{\Phi_{|_S} \mit
        \Phi\in\cnull{m}{q} (M)\} \ .\non
\eeqa
Ist $q=0$,
verzichten wir manchmal auf
den oberen Index $q$.\abstandk
Die \"au{\ss}ere Ableitung $\d$ hat die Eigenschaften
 \beqa
  \d(\Phi\wedge\Psi)&=&\d\Phi\wedge\Psi+(-1)^q\Phi\wedge\d \Psi
 \formel{dwedge} \\
   \d\d\Phi&=&0 \formel{rotrotistnull}
\eeqa
f\"ur alle $\Phi\in C_\infty^q (S)$,
$\Psi\in C_\infty^p (S)$ und ist
lokal erkl\"art durch
\beqa
 \d\Phi &=&\sum_{I\in\cals(q,N)}
  \sum_{j=1}^N\partial_j\Phi_I\dh^j\wedge\dh^{I}\non\\
  &=&\sum_{I\in{\cal S}(q+1,N)}\sum_{j\in I}\sigma(j,I-j)
  \partial_j\Phi_{I-j} \dh^{I}
 \formel{rotlokal}
\eeqa
mit $\partial_j\Phi_I:=\frac{\partial}{\partial h_j}\Phi_I
 :=\d\Phi_I(\partial_j) $
f\"ur das Differential $\d$ und den Tangentialvektor
$\partial_j$ mit
$\dh_i(\partial_j)=\delta_{i,j}$, $\Phi$ wie in $(\ref{philokal})$.
Auf Formen $\Phi\in C_\infty (S)$
wirkt $\d$ wie das Differential.
Die Koableitung $\delta:=\sigma_q*\d*$
\label{pagedelta}
hat im Falle
einer positiv orientierten Orthonormalbasis
$\{\dh^i,i=1,\cdots,N\}$
lokal die Darstellung
\beqa
  \delta\Phi&=&
  \sum_{I\in{\cal S}(q-1,N)}\sum_{j\not\in I}\sigma(j,I)
  \partial_j\Phi_{I+j} \dh^{I} \ .
 \formel{divlokal}
\eeqa
Gilt f\"ur eine Karte $(V,h)$ und Zahlen
$a_i,b_i$
\beqa
 Q:=\{x\in M \mit a_i<h_{i}(x)<b_i\ ,\ i=1,\cdots,N\}\subset V
\ ,\non
\eeqa
so definieren wir f\"ur $\Phi= \Phi_I\dh^{I}\in \cnull{\infty}{N} (M) $
\beqa
 \int_Q\Phi:=\int_{a_1}^{b_1}\cdots\int_{a_N}^{b_N}
  \Phi_I(h_1,\cdots,h_N)\dh^N\cdots\dh^1   \ .\formel{intphi}
\eeqa
Ist $\xi_\alpha$ eine der \"Uberdeckung $Q_\alpha$
untergeordnete Zerlegung der 1, so ist der Ausdruck
\beqa
 \int_M E :=\sum_\alpha \int_{Q_\alpha} \xi_\alpha \Phi
 \non
\eeqa
unabh\"angig von der Wahl der Karten.
Integration
\"uber Teilmengen von $M$ realisieren wir wie
\"ublich mit der charakteristischen Funktion.
F\"ur
$\Phi\in \cnull{\infty}{N-1} (M)$ gilt
\beqa
 \int_M\d\Phi =0 \ .\formel{ddint}
\eeqa
Eine unendlich oft differenzierbare Abbildung
$\tau:S\subset M\ra \tilde{S}\subset \tilde{M}$
induziert eine Abbildung ($x\in S$)\label{pxpushdef}
\beqa
 \Abb{\tau_*}{{\cal T}M_x}{{\cal T}\tilde{M}_{\tau(x)}}{t}{\tau_*t}
 \non\\
 \mmit  (\tau_*t)(f):=t(f\circ\tau) \ .\non
\eeqa
Wir bezeichnen den Raum der $q$--Tupel von Tangentialvektoren
aus ${\cal T}_x$ bzw. ${\cal T}M_x$ mit ${\cal T}^q_x$ bzw. ${\cal T}^qM_x$ und
erkl\"aren
f\"ur $\Phi\in A^q(\tils)$ die Form
$\tau^*\Phi\in A^q(S)$
durch
\beqa
  (\tau^*\Phi)_x(v)
  =\Phi_{\tau(x)}(\tau_*v) \formel{pulldef}
\eeqa
f\"ur alle Tangentialvektoren $v\in {\cal T}^q_x$.
Hierbei verstehen wir den Ausdruck $\tau_*v$ komponentenweise.
Die Abbildung $\tau^*:A^q(\tilde{S})\ra A^q(S)$
hat die Eigenschaften
\beqa
 \alle{\vp\in A^0(\tilde{S})}\tau^*\vp&=&\vp\circ \tau \non\\
 \alle{\Phi\in A^q(\tilde{S})}\alle{\Psi\in A^p(\tilde{S})}
 \tau^* (\Phi\wedge \Psi)&=&\tau^* \Phi\wedge \tau^*\Psi
 \non\\
 \alle{\Phi\in C_{\infty}^{q} (\tilde{S})}
 \d\tau^*\Phi&=&\tau^*\d\Phi \formel{dtausch} \\
 \alle{\Phi\in C_{\infty}^{N} (\ol{\tilde{S}})}
 \int_S \tau^*\Phi&=&\int_{\tilde{S}}\Phi
 \ ,\non
\eeqa
wobei die letzte Aussage nur f\"ur
orientierungserhaltende Diffeomorphismen $\tau$
gilt.
F\"ur solche
erkl\"aren wir
die linearen Transformationen
$\eps,\mu:A^q({S})\ra A^q(S)$\label{pxtrans} durch
\beqa
 \eps&:=&\eps_\tau:=\eps^q_\tau:=\kappa_{q}*\tau^**(\tau^*)^{-1}
 \formel{transeins}\\
 \mu&:=&\mu_\tau:=\mu^q_\tau:=\kappa_{q}\tau^**(\tau^{-1})^**
 \ .\formel{transzwei}
\eeqa
Der Kettenregel entnehmen wir die Eigenschaften
\beqa
 *\eps\tau^*&=&\tau^**\formel{taustern} \\
 \eps\mu&=&\id \ .\non
\eeqa
Lokal wirkt $\tau^*$ in folgender Weise:
Bildet $\tau$ die Koordinatenumgebung $V\subset M$
diffeomorph auf
$W\subset\tilde{M}$ ab, und sind $h_i:V\ra V_1\subset\rrr^N$
und $g_i:W\ra W_1\subset \rrr^N$ die Koordinatenabbildungen,
$F:V_1\ra W_1,x\mapsto g\circ\tau\circ h^{-1} (x)$
und
\beqa
 \Phi(w):=\sum_{I\in\cals(q,N)} \Phi_I(w)\dg^{I}\ ,\non
\eeqa
$w\in W$,
so gilt f\"ur $v\in V$ nach \cite{weck2}
\beqa
 \tau^*\Phi(v)&=&\sum_{I\in\cals(q,N)} \sum_{|J|=q}
  \sigma(J)E_{\calj(J)}(\tau(v))\partial_IF_J(h(v))\dh^{I}
  \formel{tauweck}\\
  \partial_IF_J(x)&:=&
  \partial_{i_1}F_{j_1}(x)\cdots
  \partial_{i_q}F_{j_q}(x)\ . \non
\eeqa
F\"ur Koordinaten $x_i=\tau_i(y)$ im $\rrr^N$ folgt daraus
\beqa
 \tau^*dx^{i}=d\tau_i(y)=\sum_{j=1}^{N}\partial_j\tau_i(y)dy^j
 \ .\formel{taudx}
\eeqa
\abstand
Wir betrachten noch den Spezialfall
der Inklusion:
Ist $\partial S$ eine differenzierbare
Untermannigfaltigkeit von $M$, so
gilt f\"ur die Einbettung
\beqa
 \Abb{\iota}{\partial S}{M}{x}{x} \non
\eeqa
der Satz von Stokes
\beqa
 \int_S\d\Phi =\int_{\partial S} \iota^*\Phi
 \ .\formel{dint}
\eeqa
Wir benutzen ohne weiteren Kommentar
die folgenden Konventionen:
\begin{enumerate}
 \item[i)] F\"ur eine offene Teilmenge $T_0$ von
 $T$
 sei $\iota:T_0\ra T$ die Inklusion
   $x\mapsto x$.
 Die Einschr\"ankung $\iota^*$ einer Form
 $\Phi\in A^q(T)$ auf $T_0$
 bezeichnen wir wieder mit $\Phi$
 und bemerken, da{\ss} die Einschr\"ankung
 nicht nur mit \"au{\ss}erem Produkt und
 \"au{\ss}erer Ableitung, sondern auch mit Sternoperator
 und Koableitung tauscht.
\item[ii)] Ist $\Phi\in A^q(T_0)$, $T_0\subset T$,
so bezeichnen wir
die Fortsetzung von $\Phi$ auf $T$ durch 0 auch wieder
mit $\Phi$.
\end{enumerate}
Gilt in ii) $\dist(\supp\Phi, T\setminus T_0 )>0$
und geh\"ort $\Phi$ zu irgendeinem der im
folgenden eingef\"uhrten R\"aume von $q$--Formen
auf $T_0$ (z.B. $\Phi\in \cnull{\infty}{q} (T_0)$), so
gilt dies auch f\"ur die Nullfortsetzung
($\Phi\in\cnull{\infty}{q} (T)$).
\subsection{Geometrische Voraussetzungen}
\markright{GEOMETRISCHE VORAUSSETZUNGEN}
\label{chgebiet}
Bezeichne
$U_N(R)$ die Kugel um den Nullpunkt des $\rrr^N$ mit Radius $R$
und \label{pxun}
\beqa
 U_N^{+}(R)&:=&\{x\in U_N(R) \mit x_N>0\} \non\\
 U_N^{-}(R)&:=&\{x\in U_N(R) \mit x_N<0\} \non\\
 U_N^{0}(R)&:=&\{x\in U_N(R) \mit x_N=0\} \non\\
 U_N^{0,+}(R)&:=&\{x\in U_N^0(R) \mit x_1>0\} \non\\
 U_N^{0,-}(R)&:=&\{x\in U_N^0(R) \mit x_1<0\} \non\\
 U_N^{0,0}(R)&:=&\{x\in U_N^0(R) \mit x_1=0\} \ .\non
\eeqa
Die
Sph\"are im $\rrr^N$ mit Radius $R$ nennen wir $S_N(R)$.
Im Falle $R=1$ verzichten wir auf die Angabe des Radius.
\abstandk
\begin{defi} \formel{segment} \new
\begin{enumerate}
\item[i)]
Ist $\partial S$ eine
$(N-1)$--dimensionale glatte Untermannigfaltigkeit von $M$,
und existiert um jedes $x\in\partial S$ eine
''Randkarte'', das ist eine Karte
$(V,h)$ f\"ur $M$
mit
\beqa \left. \begin{array}{lcl}
 h(x)&=&0 \\
 h(\ol{V})&=&\ol{U_N}  \\
 h(\partial S\cap V)&=& U_N^0\\
 h(S\cap V)&=& U_N^-\ ,
  \end {array}    \right\} \formel{karte}
\eeqa
so nennen wir $S$ glatt.
\item[ii)]
 Wir sagen $S$ besitzt die Segmenteigenschaft, falls
 um jedes $x\in\partial S$ eine Karte $(V,h)$ f\"ur $M$,
 ein $\rho\in(0,1)$ und
 ein Vektor $y\in\rrr^N$ existieren mit
 \beqa \left. \begin{array}{rcl}
   h(\ol{V})&=&\ol{U_N} \\
  (U_N(\rho) \cap \ol{h(S\cap V)})+\tau y&\subset&
   h(S\cap V) \mfur\malle \tau\in(0,1)\  .
  \end{array} \right\}  \formel{karte5}
 \eeqa
\end{enumerate}
\end{defi}
\begin{defi}   \formel{randeig}  \new
Sei $(S,\Gamma_1,\Gamma_2)\subset M\times M\times M$. Wir sagen
$(S,\Gamma_1,\Gamma_2)\in\calm(M)$\label{pxcalm}, falls
$ \Gamma_1$ und $\Gamma_2 $ relativ offene Teilmengen  in
$\partial S$ sind und die Eigenschaften
\beqa
  \Gamma_1\cap\Gamma_2&=&\emptyset\non\\
  \partial\Gamma_1=\partial\Gamma_2&=:&\gamma\non\\
  \Gamma_1\cup\Gamma_2\cup\gamma &=&\partial S\mbox{ disjunkt }
  \non
\eeqa
besitzen.
Die Menge $\gamma$ nennen
 wir Trennmenge.
\end{defi}
\begin{defi} \formel{randeig2} \new
Sei $(S,\Gamma_1,\Gamma_2) \in \calm(M)$ und bezeichne $\gamma$
die Trennmenge.
\begin{enumerate}
\item[i)] Ist $S$ glatt,
so hei{\ss}t $(S,\Gamma_1,\Gamma_2)$ Gebiet mit
 \"Ubergangsrand.
 \item[ii)]
Ist zus\"atzlich
$\gamma$
eine
$(N-2)$--dimensionale
glatte Untermannigfaltigkeit von $\partial S$,
und existiert um jedes $x\in\gamma$ eine
''glatte \"Ubergangsrandkarte''
$(V,h)$ f\"ur $M$ mit $(\ref{karte})$
und
\beqa    \left.
 \begin{array}{lcl}
   h(\Gamma_1\cap V)&=& U_N^{0,-} \\
   h(\Gamma_2\cap V) &=& U_N^{0,+} \\
   h(\gamma\cap V) &=& U_N^{0,0} \ ,
    \end{array}  \right\} \formel{karte2}
\eeqa
so hei{\ss}t $(S,\Gamma_1,\Gamma_2)$ glatt.
\item[iii)] Ein Gebiet $(S,\Gamma_1,\Gamma_2)$
mit \"Ubergangsrand hei{\ss}t
S--Gebiet\label{pxsgebiet} (Segment--Gebiet), wenn f\"ur ein $j\in\{1,2\}$
und um jedes $x$ aus der Trennmenge
eine \"Ubergangsrandkarte
$(V,h)$ mit $(\ref{karte})$,
ein Vektor $y\in\rrr^N$
und ein $\rho\in(0,1)$ existieren
mit
\beqa
 (U_N(\rho)\cap \ol{h(\Gamma_j\cap V)})+\tau y \subset
 h(\Gamma_j\cap V) \mfur\malle \tau\in(0,1)\ .\formel{segeig}
\eeqa
\end{enumerate}
\end{defi}
\begin{bemerk} \formel{bemerkung5} \new
Erf\"ullt $(S,\Gamma_1,\Gamma_2)$ die Bedingungen in
Definition \ref{randeig2}, iii) mit $j=1$ und
dem Vektor $y$, so auch f\"ur $j=2$ mit dem Vektor $-y$.
Daher ist die Bezeichnung S--Gebiet
unabh\"angig von $j$ in $(\ref{segeig})$.
\end{bemerk}
Wir fassen die oben gemachten Definitionen zusammen und
sammeln Eigenschaften der in dieser Arbeit
betrachteten Mengen $S$:\new
Um jeden Punkt $x\in S$ bzw. $y\in\Gamma_i$ $(i=1,2)$ bzw.
$z\in \gamma$ existiert eine Karte $(V,h)$ mit $x$ bzw. $y$ bzw.
$z\in V$. Wir k\"onnen o.B.d.A. annehmen, da{\ss}
$\partial S\cap V=\emptyset$ f\"ur Karten $(V,h)$ um
$x\in S$ und $\ol{S}\cap V\subset S\cup\Gamma_i$ f\"ur Karten um
$y\in\Gamma_i$ gilt. Im ersten Fall nennen wir $(V,h)$ eine
interne Karte, im zweiten Fall eine interne Randkarte
und im dritten Fall $(z\in\gamma)$ eine \"Ubergangsrandkarte.
Letztere beiden nennen wir gemeinsam Randkarten.
Da $S$ kompakt ist, gen\"ugt eine endliche Kollektion
$\{(V_k,h_k) \mit k=1,\cdots,K\}$ von Karten, um
$\ol{S}$ mit $\{V_k \mit k=1,\cdots K\}$ zu \"uberdecken.
Hierzu sei
$\{\xi_k\mit k=1,\cdots,K\}$\label{pxxik} eine untergeordnete
Zerlegung der 1.
Wir k\"onnen o.B.d.A. weiter annehmen,
da{\ss} stets
\beqa
  h_k(V_k)=U_N(0,1)\ \mbox{ und }
 \ \supp  \xi_k\circ h_k^{-1}\subset
  U_N(0,\frac{1}{3})\cap\ol{h_k(S\cap V_k)}
 \formel{suppined}
\eeqa
f\"ur alle $k$ erf\"ullt sind.
Je nach Regularit\"at der Geometrie unterscheiden wir drei
Typen von Daten:
\begin{enumerate}
\item[i)] Glatte Gebiete mit glattem \"Ubergangsrand:
Hier haben die internen Randkarten die Eigenschaft
$(\ref{karte})$ und die \"Ubergangsrandkarten
die Eigenschaft $(\ref{karte2})$.
\item[ii)] S--Gebiete: Interne Randkarten haben die Eigenschaft
$(\ref{karte})$, w\"ahrend \"Ubergangsrandkarten die
Eigenschaft $(\ref{segeig})$ haben.
\item[iii)] Z--Gebiete (werden in Abschnitt \ref{kompeinb}
definiert)
\item[iv)] Gebiete, an deren Rand wir keine Voraussetzungen
stellen m\"ochten: Von den Karten $(V,h)$
fordern wir lediglich $h:\ol{V}\ra \ol{U_N}$.
\end{enumerate}
%
\markright{DIE SOBOLEVR\"AUME $H_m^q$}
\subsection{Die Sobolevr\"aume $H_m^q$}   \label{chhmq}
Mit den oben erw\"ahnten Karten $(V_k,h_k)$, $k=1,\cdots, K$
definieren wir
f\"ur $m\in [0,\infty)$
die Sobolevr\"aume
$H_m^q(S)$\label{pxhmq} als die Menge der Formen $E\in A^q(S)$ mit
\beqa
  \normu{E}{H_m^q(S)}:=
 (\sum_{k=1}^K\sum_{I\in\cals(q,N)}
  \normu{E_I^k}{H_m(h_k(S\cap V_k))}^2 )^{1/2} <\infty \ ,
  \formel{hnorm}
\eeqa
wobei $E_I^k$ die Komponenten von
$(h_k^{-1})^*E$
bzgl. kartesischer Koordinaten sind (nach unserer
Konvention identifizieren wir hier die Form $E$ mit
ihrer
Ein\-schr\"an\-kung auf $S\cap V_k$).
Aus den Transformationss\"atzen, $(\ref{tauweck})$ f\"ur
Formen und \cite[Satz 4.1]{wloka} im skalaren Fall,
folgt, da{\ss} die Definition unabh\"angig von
der gew\"ahlten \"Uberdeckung ist und
verschiedene \"Uberdeckungen
\"aquivalente Normen liefern. Ebenso entnehmen wir
$(\ref{tauweck})$, da{\ss} f\"ur einen Diffeomorphismus
$\tau:T\ra S$  und $E\in H_m^q(S)$
\beqa
 c'\normu{E}{H_m^q(S)} \le \normu{\tau^*E}{H_m^q(T)}
 \le c\normu{E}{H_m^q(S)}
 \formel{pullnorm}
\eeqa
mit von $E$ unabh\"angigen Konstanten $c,c'>0$
erf\"ullt ist.
Die Vollst\"andigkeit wird ebenso auf $H_m^q(S)$
\"ubertragen wie folgende Aussagen:
\beqa
 C_\infty^{q}(S)\cap H_m^q(S) &&\mbox{ dicht in } H_m^q(S)
 \formel{hdicht3} \\
 \cnull{\infty}{q} (S)&&\mbox{ dicht in } H_0^q(S)
  \formel{ledicht}\\
 \alle{\Phi\in C_\infty^{p}(\ol{S})}\ \gibt{ c>0 }\ \alle{E\in H_m^q(S)}
 &&\normu{\Phi\wedge E}{H_m^{q+p}(S)} \le c\normu{E}{H_m^q(S)}
 \formel{multstetig}  \\
 \gibt{ c>0 }\ \alle{E\in H_m^q(S)}&& \normu{*E}{H_m^{N-q}(S)}
  \le c\normu{E}{H_m^q(S)} \formel{sternlein}
\eeqa
Wir zeigen nur $(\ref{hdicht3})$.
Zerlegen wir die Form $E\in H_m^q(S)$ in
$\sum_{k=1}^K\xi_kE$, so gen\"ugt es, $\xi_kE$ durch
Elemente aus $C_\infty^{q}(S)\cap H_m^q(S)$
zu approximieren, wobei $\{\xi_k \subset C_\infty(\ol{S})\}$
die der
\"Uberdeckung $\{V_k\}$ untergeordnete Zerlegung der 1 ist.
Wir setzen $U:=h_k(S\cap V_k)$.
Aus $(\ref{suppined})$  und \cite[Satz 3.5]{wloka} folgt,
da{\ss} wir die Komponenten $E_I^k\in H_m(U)$ von
$(h_k^{-1})^*\xi_kE$ durch Funktionen
$\Phi_I^l\in C_\infty(U)\cap H_m(U)$ approximieren
k\"onnen, deren Tr\"ager o.B.d.A. kompakt enthalten ist in
$U_N(1/3)\cap U$.
Nach $(\ref{pullnorm})$ konvergiert
die Folge
$h_k^*\Phi^l$ mit $\Phi^l:=\sum_{I\in\cals(q,N)}\Phi_I^l\dx^{I}$
in $H_m^q(S\cap V_k)$ gegen $\xi_kE$.
 Wegen der Bemerkung
nach unserer
Konvention liegen die Nullfortsetzungen der $\Phi^l$
im Raum $C_\infty^{q}(S)\cap H_m^q(S)$ und approximieren $\xi_kE$ in
$H_m^q(S)$.
Mit der gleichen Technik folgt die zweite
Behauptung aus den entsprechenden Aussagen \"uber $L_2$. Mit
\beqa
 \normu{\sum_{J\in{\cal S}(p,N)}\Phi_J\dx^J\wedge\sqqn E_I \dx^{I} }{
 H_m^{q+p}(U)} ^2&\le&  \sum_{J\in{\cal S}(p,N)} \sqqn
 \normu{\Phi_JE_I}{H_m(U)} ^2 \non\\
 &\le& c\normu{\sqqn E_I \dx^{I} }{H_m(U) } ^2 \non
\eeqa
f\"ur Teilmengen $U\subset\rrr^N$
und \cite[Lemma 3.2]{wloka}
erhalten wir $(\ref{multstetig})$.
Analog: $(\ref{sternlein})$
\abstand
Wir definieren (vgl. $(\ref{ledicht})$)\label{pxlzq}
\beqa
 L_2^q(S)&:=&H_0^q(S)\non\\
 \spu{E}{H}{q,S} &:=&\int_S E\wedge*\ol{H}
 \mfur E,H\in L_2^q(S) \ .\non
\eeqa
Die Normen in $L_2^q(S)$ und $H_0^q(S)$ sind \"aquivalent.
Nach \cite{weck} ist $L_2^q(S)$ ein Hilbert\-raum.
In diesem sind die Transformationen
$\eps_\tau,\mu_\tau$ aus $(\ref{transeins})$,
$(\ref{transzwei})$
zu einem Diffeomorphismus $\tau:S\ra T$
zul\"assig:
\begin{defi} \formel{dummy} \new
Lineare symmetrische gleichm\"a{\ss}ig positiv
definite und beschr\"ankte Transformationen
auf $A^q(x)$ nennen wir zul\"assig,
wenn f\"ur alle Karten $(V,h)$ um $x$
die Abbildungen $\eps_{I,J}$ mit
\beqa
 \eps(x)\sum_{I\in\cals(q,N)} \Phi_I(x)\dh^{I}
 =\sum_{I,J\in\cals(q,N)} \eps_{I,J}(x)\Phi_J(x)\dh^{I}
 \formel{matrixdar}
\eeqa
me{\ss}bar sind.
\end{defi}
\abstand
Weitere Eigenschaften \"ubertragen sich, wenn wir
an den Rand st\"arkere Voraussetzungen stellen.
Besitzt $S$ Segmenteigenschaft, liefert
die gleiche Argumentation wie beim Beweis von
$(\ref{hdicht3})$ mit
\cite[Satz 3.6]{wloka}
\beqa
   C_\infty^q(\ol{S})\mbox{ dicht in }  H_m^q(S) \formel{hdicht}
\eeqa
und mit \cite[Satz 3.7]{wloka}:
\begin{lemma}  \formel{eininhnull}        \new
 Besitze $S$ die Segmenteigenschaft, und sei $T$ offen
 mit $S\subset\subset T\subset\subset M$.
 Ferner sei $\Phi\in H_m^q(T)$
 mit $\Phi=0$ in $T\setminus S$. Dann gilt\label{pxhnullq}
 \beqa
  \Phi\in\hnull m q (S) :=\ol{ \cnull{\infty}{q} (S) } ^{H_m^q(S)}
  \ .\non
 \eeqa
\end{lemma}
\abstand
Mit Hilfe der eingeschr\"ankten ''Randkarten''  k\"onnen wir
den Raum $H_m^q(\partial S)$  ein\-f\"uh\-ren.
Um Spurs\"atze auf Differentialformen zu \"ubertragen,
bringen wir die bekannten Spurs\"atze auf
eine geeignetere Form.
Mit \cite[Satz 8.7]{wloka}
und einem Approximations\-argument
zeigt man, da{\ss} f\"ur $m\in\nnn$ ein linearer stetiger
Spuroperator\label{pxtnull}
\beqa
 t_0&:&\{u\in H_m(U_N^-)\mit \supp u\subset\subset
 U_N(2/3)\} \non\\
 &&\ra \{v\in H_{m-1/2}(U_N^{0})\mit \supp v\subset\subset
 U_N(2/3)\}
  \non
\eeqa
existiert mit $ t_0\Phi(x')=\Phi(x',0)$ f\"ur alle
$\Phi\in C_m(\ol{U_N^-})$ mit
$\supp\Phi\subset\subset U_N(2/3)$,
$x'\in\rrr^{N-1}$.
Nach Multiplikation mit
$\psi\in\cnull{\infty}{\ } (U_N(2/3))$, $\psi=1$ in $U_N(1/3)$
erhalten wir nach \cite[Satz 8.8]{wloka}
einen linearen stetigen Fortsetzungsoperator\label{pxchtnull}
\beqa
 \chcht_0&:&\{u\in H_{m-1/2}(U_N^0)\mit
 \supp u\subset\subset U_N(1/3)\} \non\\
 &&\ra \{v\in H_m(U_N^{-})\mit \supp v\subset\subset
 U_N(2/3)\} \ .\non
\eeqa
Es gelten
$t_0\chcht_0=\id$ und
$\chcht_0{\Phi}\in \cnull{m}{\ }(\ol{U_N^-(2/3)})$
f\"ur $\Phi\in \cnull{m}{\ }(\ol{U_N^0(1/3)})$.
Die letzte Eigenschaft entnimmt man dem Beweis zu
\cite[Satz 8.8]{wloka}.
Die vorletzte Eigenschaft folgt aus
\beqa
 \tilde{T}(\psi\Phi)&=&\tilde{T}(\psi)\tilde{T}(\Phi)
 \non\\
 \tilde{T}(\psi\tilde{Z}\Phi)
  &=&\psi_{|_{U_N^0(2/3)}}\tilde{T}\tilde{Z}\Phi=\Phi\non
\eeqa
f\"ur
$\Phi\in H_{m-1/2}(U_N^0)$ mit
$\supp\Phi\subset\subset U_N^0(1/3)$,
den Spuroperator $\tilde{T}$ und den
Fortsetzungsoperator $\tilde{Z}$ aus \cite{wloka}.
Wir k\"onnen dann zeigen:
\begin{lemma}  \formel{hspur} \new
 Seien $S$ glatt, $\iota$ die Inklusion
 $\partial S\ra M$ und $m\in \nnn$.
 Dann existiert ein linearer und stetiger
 Spuroperator $T:H_m^q(S)\ra H_{m-1/2}^{q}(\partial S)$
 mit
 \begin{enumerate}
  \item[i)] $T \Phi=\iota^*\Phi$\label{pxspur}
      f\"ur alle $\Phi\in C^q_m(\ol{S})$
  \item[ii)] $  \d T\Phi=T \d \Phi$
  f\"ur alle $\Phi\in C^q_\infty(\ol{S})$.
 \end{enumerate}
 \end{lemma}
\beweis ii) folgt aus i), und wegen $(\ref{hdicht})$ gen\"ugt
es, Linearit\"at und Stetigkeit f\"ur Formen
$\Phi$ aus $C_m^q(\ol{S})$ zu
zeigen.
Ist $(V_k,h_k)$ eine ''Randkarte '' (siehe $(\ref{suppined})$)
f\"ur $S$, so ist
$(\partial S\cap V_k,\tilde{h}_k)$ eine Karte f\"ur $\partial S$,
wobei
\beqa
 \tilde{h}_k:=\hat{\iota}^{-1}\circ h_k\circ \iota \ \Ra
 \ \iota=h_k^{-1}\hat{\iota} \tilde{h_k}\formel{ttauscht3}
\eeqa
mit $\hat{\iota}:\rrr^{N-1}\ra\rrr^{N-1}\times\{0\}$,
$x'\mapsto (x',0)$.
Wir setzen
\beqa
 T\Phi&:=&\iota^*\Phi \ .\formel{tkdef}
\eeqa
Da $T$ linear ist, gen\"ugt es, die Stetigkeit der Abbildung
\beqa
 \Phi\mapsto \iota^*\xi_k\Phi \formel{abboben}
\eeqa
zu untersuchen ($\xi_k$ wie in $(\ref{suppined})$).
Geh\"ort $\xi_k$ zu einer internen Karte,
brauchen wir nichts zu zeigen.
F\"ur Randkarten folgt aus
\beqa
  \hat{\iota}^*\dx^{I}&=&\left\{\begin{array}{rcl}
      \dx^{I}&\mbox{ falls }&N\not\in I\\
      0&\mbox{ falls }& N\in I
  \end{array}\right.   \non\\
  \hat{\iota^*}f&=&f\circ\hat{\iota}=t_0f \mfur f\in
  C_m(\ol{U_N^-})\ ,\ \supp f\subset\subset U_N(2/3)
  \non
\eeqa
und dem skalaren Spursatz die Stetigkeit
von $\hat{\iota}$, aus
$(\ref{ttauscht3})$ und
$(\ref{pullnorm})$ die Stetigkeit
der Abbildung $(\ref{abboben})$.
\qed
\begin{lemma}   \formel{invhspur} \new
 F\"ur $m\in\nnn$ existiert ein linearer
 stetiger Fortsetzungsoperator\label{pxinvspur}
 \beqa
  \cht:H_{m-1/2}^{q}(\partial S)\ra H_m^q(S) \non
 \eeqa
 mit
 $T\cht=\id$.
\end{lemma}
\beweis
Mit $\tilde{\ }$ bezeichnen wir
wieder
Einschr\"ankungen auf $\partial S$ bzw.
$U_N^0$.
Wie oben folgt aus dem skalaren
Fortsetzungssatz, da{\ss} die Abbildung
\beqa
 \cht&:=&\sum_{k=1}^K\cht_k \non\\
 \cht_k&:=& h_k^*\chcht_0(\tilde{h}_k^{-1})^*
 \tilde{\xi_k}\non\\
 \mmit\ \ \chcht_0\sum_{I\in\cals(q,N-1)}\Phi_I\dx^{I}
    &:=&\sum_{I\in\cals(q,N-1)}(\chcht_0\Phi_I)\dx^{I} \non
\eeqa
linear und stetig ist.
F\"ur
$\Phi\in C_m^q (\partial S)$
erhalten wir mit
$(\ref{ttauscht3})$ und $(t_0\chcht_0)=\id$
\beqa
 T\cht_k\Phi&=&
   \iota^*h_k^*\chcht_0(\tilde{h}_k^{-1})^*
 \tilde{\xi_k}\Phi \non\\
   &=&\tilde{h}_k^*\hat{\iota}^*\chcht_0
    (\tilde{h}_k^{-1})^*
 \tilde{\xi_k}\Phi
  = \tilde{\xi_k}\Phi  \ ,\non
\eeqa
also auch $T\cht\Phi=\Phi$. Aus $(\ref{hdicht3})$ folgt
schlie{\ss}lich die Behauptung.
\qed
\abstand
Wir k\"onnen nun den stetigen und
linearen Normalenspuroperator
\beqa
  \Abb{N:=N^q}{H_m^q(S)}{H_{m-1/2}^{q-1}(\partial S)}{
  E}{\sigma_q*T* E} \formel{definspur}
\eeqa
mit der stetigen linearen Rechtsinversen\label{pxnspur}
\beqa
  \Abb{\chn:=\chn^q}{H_{m-1/2}^{q-1}(\partial S)}{H_m^q(S)}{
  E}{\kappa_q*\cht * E} \non
\eeqa
definieren.
Die damit gemachte Behauptung folgt aus
\beqa
 *T**\cht *\vp=\kappa_q*T\cht*\vp=\kappa_q\kappa_{q-1}'\vp
 \mfur \vp\in H_{m-1/2}^{q-1}(\partial S) \ .\non
\eeqa
F\"ur $E\in C_\infty^q(\ol{S})$ erf\"ullt der Normalenspuroperator
$\div N E=-N\div E$:
\beqa
 \div N E&=&\sigma'_{q-1}\sigma_q*\d **T* E\non\\
 &=&  \sigma'_{q-1}*T\d*E\non\\
 &=& -(-1)^N*T**\d*E\non\\
 &=&-N\div E \formel{tauschdiv}
\eeqa
Testen mit $v\in T^{N-q}\partial S$ liefert
$T*\cht=0$.
Es folgt
\beqa
  T\chn =0\ .\formel{tnisnull}
\eeqa
\subsection{Die R\"aume $R^{q,\Gamma}$ und $D^{q,\Gamma}$}
\markright{DIE R\"AUME $R^{q,\Gamma}$ UND $D^{q,\Gamma}$ }
\label{chrqdq}
Wir betrachten die Abbildungen
\beqa
 &&\Abb{\rot}{C_\infty^{q}(S)}{
 C_\infty^{q+1}(S)}{\Phi}{\d\Phi} \formel{rrrot}\\
&&\Abb{\div}{C_\infty^{q}(S)}{
 C_\infty^{q-1}(S)}{\Phi}{\delta\Phi\ .}
 \formel{dddiv}
\eeqa
Diese erf\"ullen
f\"ur  $\Phi\in C_\infty^{q}(\ol{S})$,
$\Psi\in C_\infty^{q+1}(\ol{S})$ nach
$(\ref{dwedge})$
\beqa
\spu{\Phi}{\div\Psi}{q,S} &+&\spu{\rot\Phi}{\Psi}{q+1,S} \non\\
   &=&\int_S\Phi\wedge *\sigma_{q+1}*\d *\ol{\Psi}
   +\int_S\d\Phi\wedge *\ol{\Psi} \non\\
   &=&(-1)^q\int_S \Phi\wedge \d *\ol{\Psi}
   + \int_S\d\Phi\wedge *\ol{\Psi}\non\\
   &=&\int_S\d(\Phi\wedge *\ol{\Psi}) \ .\formel{phipsi3}
\eeqa
F\"ur eine Teilmenge $\Gamma\subset \partial S$ sei
$C_\infty^{q,\Gamma}(\ol{S})$\label{pxcqgamma} die Menge der Einschr\"ankungen
von Formen aus $C_\infty^{q}(M)$ auf $S$, deren Tr\"ager
einen positiven Abstand zu $\Gamma $ besitzt.
Im Falle $(S,\Gamma_1,\Gamma_2)\in\calm(M)$
gilt f\"ur  $\Phi\in\cqe(\ol{S})$,
$\Psi\in\cqzp(\ol{S})$ nach $(\ref{phipsi3})$ und $(\ref{ddint})$
\beqa
  \spu{\Phi}{\div\Psi}{q,S} +\spu{\rot\Phi}{\Psi}{q+1,S}
  = 0 \ .\formel{phipsi}
\eeqa
Wir definieren\label{pxrqgamma}
\beqa
 \rq(S)&:=&\{E\in L_2^q(S) \mit
  \gibt{F\in L_2^{q+1}(S)}\alle{\Phi\in\cqzp(\ol{S})}
  \spu{E}{\div\Phi}{q,S}
  =\spu{F}{\Phi}{q+1,S} \} \non\\
 \dq(S)&:=&\{E\in L_2^q(S) \mit
  \gibt{G\in L_2^{q-1}(S)}\alle{\Phi\in\cqem(\ol{S})}
  \spu{E}{\rot\Phi}{q,S}
  =\spu{G}{\Phi}{q-1,S} \} \non
\eeqa
und setzen $\rot E:=-F$ bzw. $\div E:=-G$.
Als Adjungierte
der dicht definierten Operatoren
$-\div_{|_{C_\infty^{q+1,\Gamma_2}(\ol{S})}}$ bzw.
$-\rot_{|_{C_\infty^{q-1,\Gamma_1}(\ol{S})}}$ sind
$\rot_{|_{\rq(S)}}$ bzw.
$\div_{|_{\dq(S)}}$ wohldefiniert, und aus $(\ref{phipsi})$
folgt, da{\ss} diese zum einen selbst wieder dicht definiert sind
und zum anderen kein Konflikt mit den  Definitionen
$(\ref{rrrot})$ bzw. $(\ref{dddiv})$ besteht.
Das gleiche gilt dann auch f\"ur die Adjungierten,
deren Definitionsbereiche wir mit\label{pxhrqgamma}
\beqa
 \begin{array}{lcrclcr}
  \hat{R}^{q,\Gamma_1}(S)&:=&D((\div_{|_{\dqp(S)}})^*)&,
  &\rot_{|_{\hat{R}^{q,\Gamma_1}(S)}}&:=&
   (-\div_{|_{\dqp(S)}})^* \non\\
  \hat{D}^{q,\Gamma_2}(S)&:=&D((\rot_{|_{\rqm(S)}})^*)&,
  &\div_{|_{\hat{D}^{q,\Gamma_2}(S)}}&:=&
   (-\rot_{|_{\rqm(S)}})^* \non
 \end{array}
\eeqa
bezeichnen wollen.
Die Spezialf\"alle
\beqa
        (S,\Gamma_1,\Gamma_2)=(S,\emptyset,\partial S)
        \mbox{ bzw. }
        (S,\Gamma_1,\Gamma_2)=(S,\partial S,\emptyset)\non
\eeqa
liefern
die \"ublichen R\"aume\label{pxrq}
\beqa
 \begin{array}{rlcrclcrc}
  &R^q(S)&:=&R^{q,\emptyset}(S)&,&
  \dqn(S)&:=&\hat{D}^{q,\partial S}(S)&, \\
 \mbox{ bzw. } &D^q(S)&:=&D^{q,\emptyset}(S)&,&
  \rqn(S)&:=&\hat{R}^{q,\partial S}(S)&,
 \end{array}                        \non
\eeqa
und auch hier gilt, da{\ss} $\rot_{|_{R^q(S)}}$ eine Fortsetzung
von $\rot_{|_{\rq(S)}}$ ist ($\div $ analog).
Nach den obigen Betrachtungen werden die R\"aume
$R^q(S)$, $\rqn(S)$, $\rq(S) $ und $\hat{R}^{q,\Gamma_1}(S)$
bzw.
$D^q(S)$, $\dqn(S)$, $\dq(S) $ und $\hat{D}^{q,\Gamma_2}(S)$,
versehen mit dem Skalarprodukt\label{pxspurq}
\beqa
 \spu{E}{H}{R^q(S)} &:=&
 \spu{E}{H}{q,S} +\spu{\rot E}{\rot H}{q+1,S} \non\\
 \mbox{ bzw. }\spu{E}{H}{D^q(S)} &:=&
 \spu{E}{H}{q,S} +\spu{\div E}{\div H}{q-1,S} \ ,\non
\eeqa
zu Hilbertr\"aumen, und es gelten wegen
$\rot_{|_{\hat{R}^{q,\Gamma_1}(S)}}
=(\rot_{|_{C_\infty^{q,\Gamma_1}(\ol{S})}})^{**}$ ($\div$ analog)
\beqa
 \hat{R}^{q,\Gamma_1}(S)&=&\ol{\cqe(\ol{S})} ^{R^q(S)}
  \formel{rqdicht}\\
 \hat{D}^{q,\Gamma_2}(S)&=&\ol{\cqz(\ol{S})} ^{D^q(S)} \ .
 \formel{dqdicht}
\eeqa
Aus $(\ref{rotrotistnull})$, der Definition von $\delta$,
$(\ref{rrrot})$, $(\ref{dddiv})$, $(\ref{rqdicht})$ und
$(\ref{dqdicht})$ folgt
\begin{lemma} \formel{rotinrnull}\new
Es gelten die Inklusionen
\beqa
\begin{array}{rclcrcl}
 \rot R^{q,\Gamma_1}(S)&\subset& R_0^{q+1,\Gamma_1}(S)&,&
 \rot \hat{R}^{q,\Gamma_1}(S)&\subset&\hat{R}_0^{q+1,\Gamma_1}(S)\\
 \div D^{q,\Gamma_2}(S)&\subset& D_0^{q-1,\Gamma_2}(S)& ,&
 \div \hat{D}^{q,\Gamma_2}(S)&\subset& \hat{D}_0^{q-1,\Gamma_2}(S).
 \end{array}
 \non
\eeqa
\end{lemma}
Das Skalarprodukt im Raum $R^q(S)\cap D^q(S)$
erkl\"aren wir
durch\label{pxspurqdq}
\beqa
 \spu{E}{H}{R^q(S)\cap D^q(S)} &:=& \non\\
 \spu{E}{H}{q,S} \hsp&+&\hsp \spu{\rot E}{\rot H}{q+1,S}
 +\spu{\div E}{\div H}{q-1,S} \ .\non
\eeqa
Wir definieren noch\label{rqrest}
\beqa \begin{array}{lclclcl}
  R^q_0(S)&:=&\{E\in R^q(S) \mit \rot E=0\}&,&
  \ \rqn_0(S)&:=&\rqn(S)\cap R^q_0(S)\\
  \rq_0(S)&:=& \rq(S) \cap R^q_0(S)&,&
  \ \hat{R}^{q,\Gamma_1}_0(S)&:=&
   \hat{R}^{q,\Gamma_1}(S)\cap R^q_0(S) \\
  D^q_0(S)&:=&\{E\in D^q(S) \mit \div E=0\}&,&
  \ \dqn_0(S)&:=&\dqn(S)\cap D^q_0(S)\\
  \dq_0(S)&:=& \dq(S) \cap D^q_0(S)&,&
  \ \hat{D}^{q,\Gamma_2}_0(S)&:= &
   \hat{D}^{q,\Gamma_1}(S)\cap D^q_0(S)
 \end{array}                              \non
\eeqa
und sammeln weitere Eigenschaften:
F\"ur $(S,\Gamma_1,\Gamma_2)\in\calm(M)$,
$E\in D^{q,\Gamma_1}(S)$ und
$\Phi\in C_{\infty}^{\Gamma_2,N-q+1}(\ol{S})$ gelten
\beqa
 \spu{*E}{\div\Phi}{N-q,S} &=&
  \int_SE\wedge\div\ol{\Phi}
  =\sigma_{N-q+1}\int_SE\wedge*\rot(*\ol{\Phi}) \non\\
  &=&-\sigma_{N-(q-1)}\int_S\div E\wedge **\ol{\Phi}
  =-\int_S\rot *E\wedge *\ol{\Phi} \ .\non
\eeqa
Es folgt
\beqa  \begin{array}{rcl}
 E\in D^{q,\Gamma_1}(S)&\Lra& *E\in R^{N-q,\Gamma_1}(S) \\
 E\in R^{q,\Gamma_1}(S)&\Lra& *E\in D^{N-q,\Gamma_1}(S)\ ,
 \end{array}
 \formel{sternr}
\eeqa
wobei sich die anderen Behauptungen analog oder durch Anwenden des
Sternoperators ergeben.
\abstand
F\"ur $\Phi\in C_\infty^q(S)$ und
$\vp\in C_\infty(S)$ k\"onnen wir den
Ausdruck $\rot(\vp\Phi)$ bilden
und erhalten mit
$(\ref{dwedge})$
\beqa
 \rot(\vp\Phi)&=&\rot\vp\wedge\Phi+\vp\rot\Phi\ .\formel{rotformel}
\eeqa
Erf\"ullt $\vp(*\Phi)$
mit gen\"ugend glattem $\Phi\in A^q(S)$
und $\vp\in C_\infty(S)$ die Formel $(\ref{rotformel})$,
so gilt
\beqa
  \div(\vp\Phi)
    &=&\sigma_q*(\rot\vp\wedge *\Phi)+\vp\div\Phi\ .
    \formel{divformel}
\eeqa
Dies impliziert:
\begin{lemma}
   \formel{randver} \new
Seien $(S,\Gamma_1,\Gamma_2)\in\calm(M)$
und $\vp\in C_\infty^{0,\Gamma_2}(\ol{S})$.
Dann gelten
\beqa \begin{array}{crcl}
 i)&E\in\rq(S)&\Ra&\vp E\in R^{q,\partial S}(S)
   \\
 &\rot(\vp E)&=&\rot\vp\wedge E+\vp\rot E\\
   &\gibt{c>0}\normu{\vp E}{R^q(S)} &\le& c\normu{E}{R^q(S)} \\
 ii)&E\in D^{q,\Gamma_1}(S)&\Ra&\vp E\in D^{q,\partial S}(S)
 \\
 &\div(\vp E)&=&\sigma_q*(\rot\vp\wedge *E)+\vp\div E\\
 &\gibt{c>0}\normu{\vp E}{D^q(S)} &\le&c \normu{E}{D^q(S)} \ ,
 \end{array} \non
\eeqa
wobei die Konstante $c$ nur von den Schranken f\"ur
$\vp$ und den ersten Ableitungen von $\vp$ abh\"angt.
\end{lemma}
\beweis
Um i) zu zeigen, w\"ahlen wir $q\in \{1,\cdots,N\}$,
$\Phi\in C_\infty^{q+1}(\ol{S})$
und erhalten mit
$\vp\Phi\in C_\infty^{q,\Gamma_2}(\ol{S})  $
nach $(\ref{divformel})$
\beqa
 \spu{\vp E}{\div\Phi}{q,S} &=&\spu{E}{\vp \div\Phi}{q,S} \non\\
  &=&\spu{E}{\div(\vp\Phi)}{q,S}
   -\sigma_{q+1}\spu{E}{*(\rot\vp\wedge *\Phi)}{q,S} \non\\
  &=&-\spu{\rot E}{\vp\Phi}{q+1,S}
   -\sigma_{q+1}\int_S E\wedge **(\rot\vp\wedge *\ol{\Phi})\ . \non
\eeqa
Das Integral formen wir weiter um:
\beqa
  \int_S E\wedge **(\rot\vp\wedge *\ol{\Phi})
  &=&\kappa_{q}\int_S E\wedge \rot\vp\wedge *\ol{\Phi}\non\\
  &=&\sigma_{q+1}\spu{\rot\vp\wedge E}{\Phi}{q+1,S} \non
\eeqa
Da die Absch\"atzungen aus $(\ref{multstetig})$ und
$(\ref{sternlein})$ folgen,
erhalten wir somit  i) f\"ur alle $q$.
Im zweiten Fall liegt $*E$ nach $(\ref{sternr})$ in
$R^{N-q,\Gamma_1}(S)$. Nach i) erf\"ullt
$*\vp E=\vp*E$ die Formel $(\ref{rotformel})$,
und die Behauptung ii) folgt aus
$(\ref{divformel})$
und $(\ref{sternr})$.
\qed
\abstandk
Besitzt der Tr\"ager von $E\in\rq(S)$ bzw.
$E\in D^{q,\Gamma_1}(S)$ zus\"atzlich
positiven Abstand zu $\Gamma_2$,
so k\"onnen wir ein $\chi\in C_\infty^{0,\Gamma_2}(\ol{S})$ finden, mit
$\chi(x)=1$ f\"ur $x\in\supp E$. Wegen $E=\chi E$ in $S$ folgt
aus Lemma \ref{randver}:
\begin{lemma} \formel{supprandver}  \new
Sei $(S,\Gamma_1,\Gamma_2)\in\calm(M)$.
Ferner besitze der Tr\"ager der Form
$E\in\rq(S)$ bzw. $E\in D^{q,\Gamma_1}(S)$
einen positiven Abstand zum Randst\"uck $\Gamma_2$.
Dann gilt
$E\in R^{q,\partial S}(S)$
bzw. $E\in D^{q,\partial S}(S)$.
\end{lemma}
\abstand
Auf $H_m^q(S)$ l\"a{\ss}t sich wegen $(\ref{hdicht3})$ und
\beqa
 \alle{\Phi\in C_\infty^q(S)\cap H_m^q(S)}
 \normu{\rot\Phi}{H_{m-1}^{q+1}(S)}
 \le c \normu{\Phi}{H_m^q(S)} \non
\eeqa
(nach $(\ref{rotlokal})$, $(\ref{dtausch})$ und
$(\ref{pullnorm})$) die Rotation\label{pxrotzw} als stetige
Fortsetzung von $\rot_{|_{C_\infty^q(S)\cap H_m^q(S)}}$
erkl\"aren, damit auch die Divergenz\label{pxdivzw}
\beqa
 \div E:=\sigma_q*\rot* E\ .\non
\eeqa
Mit $(\ref{sternlein})$ erhalten wir
\begin{lemma}    \formel{dummy2}     \new
Die Abbildungen
\beqa
  \rot&:&H_m^q(S)\ra H_{m-1}^{q+1}(S) \non\\
  \div&:&H_m^q(S)\ra H_{m-1}^{q-1}(S) \non
\eeqa
sind linear und stetig.
\end{lemma}
\begin{lemma} \formel{taudivrot} \new
 Sei $(S,\Gamma_1,\Gamma_2)\in\calm(M)$.
 F\"ur einen orientierungserhaltenden
 Diffeomorphismus $\tau:S\ra T$,
 $\eps:=\eps_\tau$ (vgl. $(\ref{transeins})$) und
 $\hat{\Gamma}_i:=\tau(\Gamma_i)$
 gelten f\"ur $q=1,\cdots,N$
 \beqa \begin{array}{rrcl}
  i)&E\in R^{q,\hat{\Gamma}_1}(T)&\Ra&\tau^*E\in R^{q,\Gamma_1}(S)\\
  &\rot\tau^*E&=&\tau^*\rot E \\
  &\gibt{c>0}\normu{\tau^* E}{R^q(S)} &\le&c\normu{E}{R^q(T)} \\
 ii)& E\in D^{q,\hat{\Gamma}_1}(T)&\Ra&\tau^*E\in
  \eps^{-1}D^{q,\Gamma_1} (S)
  \\
  &\div\eps\tau^*E&=&\eps\tau^*\div E\\
  &\gibt{c>0}\normu{\tau^*E}{q,S} ^2 +\normu{\div\eps\tau^*E}{q-1,S} ^2&\le&
  c\normu{E}{D^q(T)} ^2 \ .
 \end{array} \non
 \eeqa
\end{lemma}
\beweis
Gilt f\"ur ausreichend glattes $\Phi\in A^q(T)$
\beqa
 \rot\tau^**\Phi=\tau^*\rot*\Phi \ ,\formel{bedingung}
\eeqa
so auch
\beqa
 \div\eps\tau^*\Phi&=&\sigma_q\kappa_{q}
   *\rot **\tau^**\Phi=\sigma_q*\tau^*\rot*\Phi\non\\
   &=&\sigma_q\kappa_{q-1}*\tau^***\rot *\Phi\non\\
   &=&\eps\tau^*\div \Phi\ .\formel{epsdiv}
\eeqa
Aus $(\ref{dtausch})$
und $(\ref{taustern})$
folgen
f\"ur $\Phi\in C_\infty^{q+1,\Gamma_2}(\ol{S})$,
$E\in R^{q\hat{\Gamma}_1}(T)$ und $q=1,\cdots, N$
\beqa
 \eps_{\tau^{-1}}(\tau^{-1})^*\Phi&\in&
 C_\infty^{q+1,\hat{\Gamma}_2}
   (\ol{T})
   \non\\
   \spu{\tau^*E}{\div\Phi}{q,S} &=&
     \int_S\tau^*E\wedge *\div\ol{\Phi}
     = \int_{T}E\wedge(\tau^{-1} )^**\div\ol{\Phi} \non\\
    &=&\int_{T}E\wedge *\eps_{\tau^{-1}}(\tau^{-1})^*\div\ol{\Phi}
     =\int_{T}E\wedge *\div
     \eps_{\tau^{-1}}(\tau^{-1})^*\ol{\Phi} \non\\
    &=&-\int_{T}\rot E\wedge *
     \eps_{\tau^{-1}}(\tau^{-1})^*\ol{\Phi}
     =-\int_{T}\rot E \wedge(\tau^{-1})^**\ol{\Phi} \non\\
     &=&-\int_S\tau^*\rot E\wedge *\ol{\Phi}
      =-\spu{\tau^*\rot E}{\Phi}{q+1,S} \ .\non
\eeqa
Da die Absch\"atzung aus $(\ref{pullnorm})$
folgt, ist i) f\"ur alle $q$
bewiesen.
Hiermit und mit $(\ref{sternr})$
erhalten wir
\beqa
 E\in D^{q,\hat{\Gamma}_1}(T)&\Lra&*E\in R^{N-q,\hat{\Gamma}_1}(T)\non\\
   &\Lra&\tau^**E\in R^{N-q,\Gamma_1}(S)
  \mbox{ und } \rot\tau^**E=\tau^*\rot *E\non\\
  &\Lra&\eps\tau^*E=\kappa_{q}*\tau^**E\in D^{q,\Gamma_1}(S)
  \mbox{ und } \div\eps\tau^*E=\eps\tau^*\div E\ , \non
\eeqa
wobei die letzte Behauptung aus
$(\ref{epsdiv})$ folgt und die Absch\"atzung in ii) impliziert.
\qed
\abstandk
Eine besondere Bedeutung haben die R\"aume $R^0$:
F\"ur Teilmengen $U\subset\rrr^N$ seien $E=e(x)\in R^0(U)$ und
$\rot E=F=\sum_{i=1}^N F_i\dx^{i}$ sowie
$\Phi=\sum_{i=1}^N\Phi_i\dx^{i}\in\cnull{\infty}{1} (U)$.
Nach $(\ref{divlokal})$ gilt
\beqa
 0&=&\spu{E}{\div\Phi}{0,U} +\spu{F}{\Phi}{1,U} \non\\
  &=&\int_Ue(x)\div\ol{\vec{\Phi} }(x) \dx +
   \int_U\vec{F} (x)\cdot\ol{\vec{{\Phi}} }(x) \dx \ ,\non
\eeqa
wobei $\vec{F} $ f\"ur das aus den Komponenten
von $F$ erstellte Feld und $\cdot$ f\"ur das innere
Produkt in $\ccc^N$ steht, $\vec{\Phi} $ analog.
Der Operator $\div$ ist hier im \"ublichen Sinne zu
verstehen.
Wir erhalten
\beqa
 *D^N(U)=R^{0}(U)=H_1(U) \ .\formel{risth}
\eeqa
\subsection{Die R\"aume $\rvnull{-1/2}{q} $ und $\dvnull{-1/2}{q} $}
\markright{DIE R\"AUME $\rvnull{-1/2}{q} $ UND $\dvnull{-1/2}{q} $}
F\"ur $m\in (0,\infty)$ bezeichne
$\hvnull{-m}{q} (S)$\label{pxhminus} den
Dualraum von
$\hnull{m}{q} (S)$ und
$\spu{\Lambda}{\Phi}{\hvnull{-m}{q} (S)} $\label{pxspuhminus},
$ \Lambda\in \hvnull{-m}{q} (S)$,
$\Phi\in \hnull{m}{q} (S)$
die Dualit\"at.
F\"ur diese und alle weiteren auftretenden Dualit\"aten
fordern wir stets
Antilinearit\"at in der zweiten
Komponente.
Erkl\"aren wir Rotation, Divergenz
und Sternoperator\label{pxdefhminus}
durch
\beqa
  \spu{\rot\Lambda}{\Phi}{\hvnull{-(m+1)}{q+1} (S)}
  &:=&-\spu{\Lambda}{\div\Phi}{\hvnull{-m}{q} (S)}
  \mfur \Phi\in \hnull{m+1}{q+1} (S)\formel{trot}\\
   \spu{\div\Lambda}{\Phi}{\hvnull{-(m+1)}{q-1} (S)}
  &:=&-\spu{\Lambda}{\rot\Phi}{\hvnull{-m}{q} (S)}
  \mfur \Phi\in \hnull{m+1}{q-1} (S)\non\\
    \spu{*\Lambda}{\Phi}{\hvnull{-m}{N-q} (S)}
  &:=&\kappa_{q}\spu{\Lambda}{*\Phi}{\hvnull{-m}{q} (S)}
  \mfur \Phi\in \hnull{m}{N-q} (S) \ ,\formel{sternrot}
\eeqa
so gelten
\beqa
 \spu{*\Lambda}{*\Phi}{\hvnull{-m}{N-q} (S)}
  &=&\kappa_{q}\spu{\Lambda}{**\Phi}{\hvnull{-m}{q} (S)}
  =\spu{\Lambda}{\Phi}{\hvnull{-m}{q} (S)} \non\\
 \spu{**\Lambda}{\Phi}{\hvnull{-m}{q} (S)}
  &=&\kappa_{q}\spu{*\Lambda}{*\Phi}{\hvnull{-m}{N-q} (S)}
  =\kappa_{q}\spu{\Lambda}{\Phi}{\hvnull{-m}{q} (S)} \non\\
 \div\Lambda&=&\sigma_q*\rot*\Lambda \ .\formel{divlambda}
\eeqa
Wir definieren\label{pxrvnull}
\beqa
 \rvnull{-1/2}{q} (S)&:=&\{\Lambda\in \hvnull{-1/2}{q} (S)\mit
  \rot\Lambda\in \hvnull{-1/2}{q+1} (S) \}  \non\\
 \dvnull{-1/2}{q} (S)&:=&\{\Lambda\in \hvnull{-1/2}{q} (S)\mit
  \div\Lambda\in \hvnull{-1/2}{q-1} (S) \}  \ ,\non
\eeqa
wobei $\rot\Lambda\in \hvnull{-1/2}{q+1} (S)$ bedeutet,
da{\ss} wir $\rot\Lambda$ stetig auf $\hvnull{-1/2}{q+1} (S)$
fortsetzen k\"onnen.
Aus $(\ref{divlambda})$ folgt
\beqa
 \rvnull{-1/2}{q} (S) =*\dvnull{-1/2}{N-q} (S) \ .\formel{sterndv}
\eeqa
F\"uhren wir in den R\"aumen $\rvnull{-1/2}{q} $
und $\dvnull{-1/2}{q} $
die Normen\label{pxrminusnorm}
\beqa
  \normu{\Lambda}{\rvnull{-1/2}{q} (S)} &=&
      \normu{\Lambda}{\hvnull{-1/2}{q} (S)}
      +\normu{\rot\Lambda}{\hvnull{-1/2}{q+1} (S)} \non\\
   \normu{\Lambda}{\dvnull{-1/2}{q} (S)} &=&
      \normu{\Lambda}{\hvnull{-1/2}{q} (S)}
      +\normu{\div\Lambda}{\hvnull{-1/2}{q-1} (S)} \non
\eeqa
ein, so ist die durch $(\ref{sterndv})$ induzierte Abbildung
isometrisch.\new
\subsection{Approximationseigenschaften}
\markright{APPROXIMATIONSEIGENSCHAFTEN}
Dem Abschnitt \ref{chrqdq} k\"onnen wir entnehmen,
da{\ss} stets $\hat{R}^{q,\Gamma_1}(S)\subset R^{q,\Gamma_1}(S)$
erf\"ullt ist. In diesem Abschnitt wollen wir
Kriterien f\"ur die Gleichheit finden.
Dies ist nach $(\ref{rqdicht})$ gleichbedeutend mit
der Frage, unter
welchen Voraussetzungen wir Formen aus $R^{q,\Gamma_1}(S)$
durch Folgen von Formen aus $C_\infty^{q,\Gamma_1}(\ol{S})$
approximieren k\"onnen.
\abstandk
F\"ur Formen $E\in R^q(U)$
mit $\supp E\subset\subset U$
k\"onnen wir durch Anwendung von Mollifiern
(vgl. \cite[Theorem 2.1]{agmon})
die Existenz von Folgen
$(\Phi_n)\subset \cqn(U)$ zeigen mit
\beqa
 \Phi_n\ra E \mbox{ in } R^q(U) \ .\formel{suppdicht}
\eeqa
Daraus folgen Approximationseigenschaften
weiterer R\"aume, wenn wir an $U$ strengere
Voraussetzungen stellen:
\begin{lemma} \formel{lemm0}  \new
 Besitzt $S$ die Segmenteigenschaft, so gilt
 \beqa
  R^q(S)=\ol{C_\infty^q(\ol{S})}^{R^q(S)} \ .\non
 \eeqa
\end{lemma}
\beweis Diese Aussage k\"onnen wir mit der gleichen Technik
beweisen wie die entsprechende Aussage in den
skalaren Sobolevr\"aumen (siehe \cite[Theorem 2.1]{agmon}).
\qed
\abstandk
Hat $S$ Segmenteigenschaft, so folgt daraus
\beqa
 R^{q,\partial S}(S)=\hat{R}^{q,\partial S}(S)
 =\rqn(S)\ .\formel{lemm2}
\eeqa
Wir erhalten weiter:
\begin{lemma} \formel{lemm3} \new
Besitze $S$ die Segmenteigenschaft. Seien ferner
$\Gamma_1\subset\partial S$ offen und
$E\in R^q(S)$ mit $\dist(\supp E,\Gamma_1)>0$.
 Dann gilt
 \beqa
  E\in \hat{R}^{q,\Gamma_1}(S)
    \ .\non
  \eeqa
\end{lemma}
\beweis Wir w\"ahlen $\chi\in C_\infty^{0,\Gamma_1}(\ol{S})$
mit $\chi=1$ in $\supp E$ und
gem\"a{\ss} Lemma \ref{lemm0}
eine Folge $(\Phi_n)\subset C_\infty^q(\ol{S})$
mit $\Phi_n\ra E$ in $R^q(S)$. Aus $(\ref{rotformel})$
folgt
\beqa
 \chi\Phi_n\ra \chi E =E \mbox{ in } R^q(S) \non
\eeqa
und damit die Behauptung. \qed
\begin{lemma}  \formel{dichtelemma}          \new
 Sei $(S,\Gamma_1,\Gamma_2)\in\calm(M)$.
 Ferner existiere ein
 $i\in\{1,2\}$, so da{\ss}
 \beqa
  \alle{q\in\zzz} R^{q,\Gamma_i}(S)=\hat{R}^{q,\Gamma_i}(S)
 \mbox{ oder } \alle{q\in\zzz} D^{q,\Gamma_i}(S)=\hat{D}^{q,\Gamma_i}(S)\non
 \eeqa
 erf\"ullt ist.
 Dann gelten f\"ur alle $i$ und alle $q$
 \begin{enumerate}
  \item[i)] $R^{q,\Gamma_i}(S)=\hat{R}^{q,\Gamma_i}(S)$
   \item[ii)] $D^{q,\Gamma_i}(S)=\hat{D}^{q,\Gamma_i}(S)$ .
  \end{enumerate}
\end{lemma}
\beweis Gelte $R^{q,\Gamma_1}(S)=\hat{R}^{q,\Gamma_1}(S)$
f\"ur alle $q$.
Aus $E\in D^{q,\Gamma_1}(S)$ und
$(\ref{sternr})$ folgt $*E\in R^{N-q,\Gamma_1}(S)$, und
wir finden eine Folge
$(\Phi_n)\subset C_\infty^{N-q,\Gamma_1}(\ol{S})$,
die $*E$ in $R^{N-q}(S)$ approximiert.
Wir erhalten $\kappa_{q}*\Phi_n\in\cqe(\ol{S})$ und
mit $(\ref{sternlein})$
\beqa
 \normu{\div E-\div(\kappa_{q}*\Phi_n)}{q-1,S}
  &=&\normu{\rot*E-\rot\Phi_n}{N-q+1,S}
   \longrightarrow 0\non\\
  \normu{E-\kappa_{q}*\Phi_n}{q,S}
  &=&\normu{*E-\Phi_n}{N-q,S}
   \longrightarrow 0\non
\eeqa
f\"ur $n\ra\infty$, also
$D^{q,\Gamma_1}(S)=\hat{D}^{q,\Gamma_1}(S)$ f\"ur alle $q$.
\abstand
Somit gibt es zu jedem $H\in D^{q+1,\Gamma_1}(S)$
eine Folge $(H_k)\subset C_\infty^{q+1,\Gamma_1}(\ol{S})$,
die in $D^{q+1}(S)$ gegen $H$ konvergiert. F\"ur
$E\in R^{q,\Gamma_2}(S)$ folgt
\beqa
 \spu{E}{\div H}{q,S} &\leftarrow& \spu{E}{\div H_k}{q,S} \non\\
   &=&-\spu{\rot E}{H_k}{q+1,S} \non\\
   &\ra&-\spu{\rot E }{H}{q+1,S} \ ,\non
\eeqa
also auch
\beqa
 \hat{R}^{q,\Gamma_2} (S) &=&R^{q,\Gamma_2}(S) \non\\
 \hat{D}^{q,\Gamma_2} (S) &=&D^{q,\Gamma_2}(S) \ .\non
\eeqa
Analog kann man die Behauptungen unter den
anderen Voraussetzungen zeigen.
\qed
\begin{defi} \formel{dichtdefi} \new
Gebiete $(S,\Gamma_1,\Gamma_2)$, welche die Voraussetzung
von Lemma \ref{dichtelemma} erf\"ullen, nennen wir
D--Gebiete und schreiben $(S,\Gamma_1,\Gamma_2)\in\calm_D(M)$.
\end{defi}
\begin{bemerk} \formel{bemerk2} \new
Aus Lemma \ref{taudivrot}
folgt,
da{\ss} D--Gebiete durch Diffeomorphismen wieder
auf D--Gebiete abgebildet werden.
\end{bemerk}
\begin{lemma} \formel{lemm4}     \new
Seien $U\subset\rrr^N$, $y\in\rrr^N$ und
\beqa
 \Abb{\eta_s:=\eta_{s,y}}{\rrr^N}{\rrr^N}{x}{x+sy} \ .\non
\eeqa
Ferner sei $s_0>0$ und $E\in R^q(U)$
mit $\eta_s^*E\in R^q(U)$ f\"ur alle
$s\in[0,s_0)$.
Dann ist die Abbildung
 \beqa
  {[0,s_0)}\ra {R^q(U)}\ ,\ {s}\mapsto {\eta_s^*E} \non
 \eeqa
stetig in $0$.
\end{lemma}
\beweis Wegen $\rot\eta_s^*E=\eta_s^*\rot E$ gen\"ugt es,
$ \normu{E-\eta_s^*E}{q,U} <\eps$ f\"ur $s$  klein
zu zeigen. Da aber f\"ur Koordinaten $z_i:=(x+sy)_i$
nach dem Transformationssatz $\eta_s^*\d z^i=\d x^i$ gilt,
folgt die Behauptung aus $(\ref{hnorm})$ und
der Stetigkeit der Verschiebung in $L_2(\rrr^N)$ .
\qed
\begin{satz}  \formel{dichtesatz}          \new
 Sei $(S,\Gamma_1,\Gamma_2)\in\calm(M)$ ein
 S--Gebiet.
 Dann gelten f\"ur $i=1$,
 $i=2$ und alle $q$
   \begin{enumerate}
  \item[i)] $R^{q,\Gamma_i}(S)=\hat{R}^{q,\Gamma_i}(S)$
   \item[ii)] $D^{q,\Gamma_i}(S)=\hat{D}^{q,\Gamma_i}(S)$ .
  \end{enumerate}
\end{satz}
\beweis  Nach Lemma \ref{dichtelemma} brauchen wir
nur i) mit $i=1$ zu zeigen.
Wir benutzen $(\ref{segeig})$ f\"ur $j=2$
(vgl. Bemerkung \ref{bemerkung5}).
Es gen\"ugt eine Form
$E\in R^{q,h(\Gamma_1\cap V)}(h(S\cap V))$ f\"ur
Karten $(V,h)$ mit
$\supp E\subset\subset U_N(1/3)$ durch
Formen $\Phi$ aus $C_\infty^{q,h(\Gamma_1\cap V)}(\ol{h(S\cap V)})$
mit $\supp\Phi\subset\subset U_N$
zu approximieren.
Nach Multiplikation mit
$\zeta\in\cnull{\infty}{\ } (U_N)$,
$\zeta=1$ in $U_N(1/3)$ folgt dies
f\"ur innere Karten
aus $(\ref{suppdicht})$ und f\"ur interne
Randkarten mit $\partial S\cap V\subset\Gamma_2$
aus Lemma \ref{lemm3}, mit $\partial S\cap V\subset\Gamma_1$
aus Lemma \ref{supprandver} und $(\ref{lemm2})$.
Seien nun $\rho$, $y$ wie in $(\ref{segeig})$,
$\gamma_i:=h(\Gamma_i\cap V)$
und
$\eta_s:=\eta_{s,-y}$.
Nach den Lemmata \ref{lemm4} und \ref{taudivrot}
gilt f\"ur vorgegebenes $\eps$ und
$s$ klein
\beqa
 \normu{\eta_s^*E-E}{R^q(U_N^-)}
 &\le& \frac{\eps}{2} \non\\
  \supp \eta_s^*E &\subset\subset& U_N(2/3)\cap \ol{h(S\cap V)}
 \formel{suppe} \\
  \eta_s^*E&\in& R^{q,\gamma_{1,s}}(U_N^-) \ ,\non
\eeqa
wobei
$\gamma_{i,s}:=(\gamma_i+sy)\cap U_N^0$.
Wegen $\dist(\gamma_1,\gamma_{2,s})>0$ existiert
$\chi\in C_\infty^{\gamma_{2,s}}(\ol{U_N^-})$ mit
$1-\chi\in C_\infty^{\gamma_1}(\ol{U_N^-})$.
Wir erhalten mit Lemma \ref{lemm3}
\beqa
 (1-\chi)\eta_s^*E\in\hat{R}^{q,\gamma_1}(U_N^-) \non
\eeqa
und mit $(\ref{suppe})$, Lemma \ref{supprandver},
und $(\ref{lemm2})$
\beqa
 \chi\eta_s^*E\in \rqn(U_N^-) \ .\non
\eeqa
Wir k\"onnen also bis auf eine Genauigkeit von
$\eps/4$ die Summanden $(1-\chi)\eta_s^* E$
durch Formen aus $C_\infty^{q,\gamma_1}(\ol{U_N^-})$
und $\chi\eta_s^* E$
durch Formen aus $\cqn(\ol{U_N^-})$ approximieren.
Die Multiplikation mit $\zeta$ (siehe oben) liefert
die Behauptung.
\qed
\subsection{Kegelspitzen} \label{chkegelspitzen}
\markright{KEGELSPITZEN}
F\"ur Elemente
$(B,\gamma_1,\gamma_2)\in\calm(S_N)$ sei\label{pxkegelspitze}
\beqa
  C_R(B)&:=&\{rm \mit r\in(0,R)\ ,\ m\in B\} \non\\
  C_R(B,\gamma_1,\gamma_2)&:=&
  (C_R(B),C_R(\gamma_1),C_R(\gamma_2)) \ .\non
\eeqa
Im Falle $R=1$ verzichten wir auf den Index $R$.
\abstandk
Wir wollen die Wirkung von Rotation und Divergenz auf den
Tangential-- bzw. Normalenanteil von Formen auf
Mengen $C_R(B)$
untersuchen. Dazu zitieren wir zun\"achst
einige
Resultate aus
\cite{wwspherical}:\new
Aus einer Koordinatenabbildung
$\vp:V\subset S_N\ra U\subset \rrr^{N-1}$
erhalten wir durch die Vorschrift
\beqa
\Abb{\tilde{\vp}}{\tilde{V}}{\rrr^+\times U\subset\rrr^N}{
 rt}{(r,u(t))} \non
\eeqa
eine Koordinatenabbildung f\"ur
$\tilde{V}:=\{rt\,|\,r\in\rrr^+,t\in V\}\subset\rrrp{N} :=\rrr^N\setminus\{0\}$\label{pxrpunkt}.
Ist $\{\Psi_1(t)\d u^1,\cdots, \Psi_{N-1}(t)\d u^{N-1}\}$ eine Orthonormalbasis
f\"ur $A^1(t)$, so ist
\beqa
  \{\d r,r\Psi_1(t)\d u^1,\cdots,r\Psi_{N-1}(t)\d u^{N-1}\} \non
\eeqa
eine
Orthonormalbasis f\"ur $A^1(rt)$, $rt\in\tilde{V}$.
Mit den  kartesischen Koordinaten  $x_i$ definieren
wir\label{pxwwsph}
\beqa
       X\hspace{-0.4cm}&:=&\hspace{-0.4cm}\sum_{n=1}^N x_n \dx^n \non\\
 \hat{R}:A^q(\rrr^N)&\longrightarrow &A^{q+1}(\rrr^N)\non\\
    E&\longmapsto & X\wedge E \non\\
 \hat{T}:A^q(\rrr^N)&\longrightarrow & A^{q-1}(\rrr^N)\non\\
 E&\longmapsto &\sigma_q*(X\wedge * E) \non\\
 m:A^q(\rrr^N)&\longrightarrow & A^q(\rrr^N)\non\\
 E&\longmapsto &|x|\cdot E\ .\non
 \formel{wwoperator}
\eeqa
F\"ur Formen $E\in A^q(\rrr^N)$ und $H\in A^{q-1}(\rrr^N)$
gelten dann
\beqa
 E\wedge *\hat{R}H&=&\kappa_q(*E)\wedge (X\wedge H) \non\\
   &=&\sigma_q X\wedge (*E)\wedge H\non\\
   &=&\hat{T}E\wedge *H \formel{wwrt} \\
  (\hat{R}\hat{T}+\hat{T}\hat{R})E&=&m^2E\ .\formel{wwmz}\\
  \dr&=&m^{-1}X\mfur r(x):=|x|\non
\eeqa
Eine Form $E\in A^q(\rrrp{N} )$ zerlegen wir eindeutig gem\"a{\ss}
\beqa
 E=\dr \wedge E^{\rho}+E^\tau\ ,E^\rho:=m^{-1}\hat{T}E\ ,
 E^\tau:=m^{-2}\hat{T}\hat{R}E\formel{wwmtr}
\eeqa
in ihren
Tangential-- und Normalenanteil.
Dies induziert die surjektiven Abbildungen\label{pxrhotau}
\beqa
 \rho:=\rho^q:A^q(\rrrp{N} )&\ra& \calf(\rrr^+,A^{q-1}(S_N))\non\\
 \tau:=\tau^q:A^q(\rrrp{N} )&\ra& \calf(\rrr^+,A^{q}(S_N)) \ ,\non
\eeqa
die
lokal durch (vgl. $(\ref{wwmtr})$)
\beqa
 \left[ \rho E(r) \right] (t)&:=&r^{-(q-1)}\ssqmn c_I^\rho(r,t)\Psi_I(t)\d u^I\non\\
  \mfur E^{\rho}&=&\ssqmn c_I^\rho(r,t)\Psi_I(t)\d u^I\non\\
 \left[ \tau E(r) \right] (t)&:=&r^{-q}\ssqn c_I^\tau(r,t)\Psi_I(t)\d u^I\non\\
  \mfur E^\tau&=&\ssqn c_I^\tau(r,t)\Psi_I(t)\d u^I
  \non
\eeqa
definiert sind und Rechtsinverse\label{pxchrhotau}
$\chrho:=\chrho^q$,
$\chtau:=\chtau^q$ besitzen mit
\beqa
 \rho\chrho&=&\id\mbox{ in }\calf(\rrr^+,A^{q-1}(S_{N}))\non\\
  \tau\chtau&=&\id\mbox{ in }\calf(\rrr^+,A^{q}(S_{N}))\non\\
  \chrho\rho +\chtau\tau&=&\id\mbox{ in } A^q(\rrrp{N} )\non\\
  \rho\chtau\hsp&=&\hsp 0\ ,
  \ \tau\chrho=0\ ,\ \chrho_{q+1}=\d r\wedge \chtau_q\ .\non
\eeqa
Weitere Resultate aus \cite{wwspherical} k\"onnen wir,
teils aufgrund ihrer
lokalen Eigenschaften, teils durch Multiplikation mit
charakteristischen Funktionen,
auf unsere Situation \"ubertragen:
F\"ur offene Teilmengen $B\subset S_N$ und $I:=(0,R)$
sei $\calL^q:=\calL^q_R$\label{pxcalfq} die Menge
der Bochner--me{\ss}baren Funktionen $f\in\calf(I,L_2^q(B))$
mit $\spu{f}{f}{\calL^q} <\infty$, wobei
\beqa
  \spu{u}{v}{\calL^q}& :=&\int_0^R r^{N-1}\spu{u(r)}{v(r)}{q,B} \dr
 \ .\non
\eeqa
Die Abbildungen
\beqa \begin{array}{rcccc}
 \rho&:&L_2^q(C_R(B))&\longrightarrow& \calL^{q-1} \\
 \tau&:&L_2^q(C_R(B))&\longrightarrow& \calL^{q} \\
 \chrho&:&\calL^{q-1}&\longrightarrow& L_2^q(C_R(B)) \\
 \chtau&:&\calL^{q}&\longrightarrow& L_2^q(C_R(B))
 \end{array}     \formel{wwsphiso}
\eeqa
sind stetig, die letzten beiden sogar isometrisch.
Die Zerlegung
\beqa
  L_2^q(C_R(B))=\chrho\rho L_2^q(C_R(B))  \oplus
   \chtau\tau L_2^q(C_R(B)) \formel{wwsph}
\eeqa
ist orthogonal.
Wir definieren \label{pagelzz}
\beqa
   L_{2,N}(I):=\{u\mbox{ me{\ss}bar }
   \mit \normu{u}{L_{2,N}(I)} :=\int_I r^{N-1}|u(r)|^2 \dr<\infty\}
    \non
\eeqa
und folgern aus dem Satz von Fubini:
\begin{lemma} \formel{approx} \new
Konvergiert im Raum $L_2^q(B)$ eine Folge $E_k$ gegen $E$,
so konvergiert f\"ur alle
$\vp$ aus $L_{2,N}(I)$ die Folge $\vp E_k$ in
$\calL^q$ gegen $\vp E$.
\end{lemma}
Seien $\rrot$\label{pxgrot} bzw. $\ddiv$\label{pxgdiv} die Rotation bzw. Divergenz
auf der Einheitssph\"are und\label{pxgm}\label{pxgd}
\beqa
        \hat{M}\Psi(r):=r\Psi(r)\ ,
        \ D\Psi(r):=\partial_r\Psi(r) \non
\eeqa
f\"ur gen\"ugend glattes
$\Psi\in \calL^q$.
Dann gelten f\"ur $\Phi\in C_\infty^q (Z)$
\beqa  \left.
\begin{array}{ccc}
 \rho\div \Phi&= &-\hat{M}^{-1}\ddiv\rho \Phi \\
 \tau\div \Phi&= &\hat{M}^{-(N-q)}D\hat{M}^{N-q}\rho \Phi+\hat{M}^{-1}\ddiv \tau \Phi \\
  \rho\rot \Phi&= &-\hat{M}^{-1}\rrot \rho \Phi+\hat{M}^{-q}D\hat{M}^q\tau \Phi\\
  \tau \rot \Phi&=&\hat{M}^{-1}\rrot\tau \Phi\ .
 \end{array}  \right\} \formel{rrotddiv}
\eeqa
Wir erhalten:
\begin{lemma} \formel{walze} \new
Seien $R\in(0,\infty)$, $I:=(0,R)$,
$(B,\gamma_1,\gamma_2)\in\calm(S_N)$,
$(Z,\Gamma_1,\Gamma_2):=C(B,\gamma_1,\gamma_2)$ und
$\vp\in\cnull{\infty}{\ } (I)$.
Ferner besitze $B$ die Segmenteigenschaft.
Dann gelten
\begin{enumerate}
 \item[i)] f\"ur $H\in \hat{D}^{q,\Gamma_2}(Z)$ und
   $e\in R^{q-1,\gamma_1}(B)$
\beqa
 \spu{\tau H}{\hat{M}^{-1}\vp\rrote e}{\calL^q}
 &+&\spu{\rho H}{\hat{M}^{-(q-1)}D\hat{M}^{q-1}\vp e}{\calL^{q-1}} \non\\
 &+&\spu{\tau\div H}{\vp e}{\calL^{q-1}} =0 \ ,\formel{walzee}
\eeqa
\item[ii)] f\"ur $H\in \hat{D}^{q,\Gamma_2}(Z)$ und
   $e\in R^{q-2,\gamma_1}(B)$
\beqa
  -\spu{\rho H}{\hat{M}^{-1}\vp\rrote e}{\calL^{q-1}}
 +\spu{\rho\div H}{\vp e}{\calL^{q-2}}=0 \ ,\formel{walzez}
\eeqa
\item[iii)] f\"ur $E\in \hat{R}^{q,\Gamma_1}(Z)$ und
   $h\in D^{q+1,\gamma_2}(B)$
\beqa
  \spu{\tau E}{\hat{M}^{-1}\vp\ddivh h}{\calL^q}
 +\spu{\tau\rot E}{\vp h}{\calL^{q+1}} =0\ , \formel{walzec}
\eeqa
\item[iv)] f\"ur $E\in \hat{R}^{q,\Gamma_1}(Z)$ und
   $h\in D^{q,\gamma_2}(B)$
\beqa
 \spu{\tau E}{\hat{M}^{-(N-q-1)}D\hat{M}^{N-q-1}\vp h}{\calL^q}
 &-&\spu{\rho E}{\hat{M}^{-1}\vp\ddivh h}{\calL^{q-1}} \non\\
 &+&\spu{\rho\rot E}{\vp h}{\calL^q}=0 \ . \formel{walzed}
\eeqa
 \end{enumerate}
\end{lemma}
\beweis
Wegen der Stetigkeit der Abbildungen $\tau$ und $\rho$
gen\"ugt es, die Aussagen f\"ur $H\in C_\infty^{q,\Gamma_2}(\ol{Z})$
und $E\in C_\infty^{q,\Gamma_1}(\ol{Z})$ zu zeigen.
Seien  $e_1\in C_\infty^{q-1,\gamma_1}(\ol{B})$,
$e_2\in C_\infty^{q-2,\gamma_1}(\ol{B})$ und
$\tilde{E}=\chrho\vp e_2+\chtau\vp e_1$.
Wegen
\beqa
  \dist(\supp(*H\wedge \tilde{E}),\partial Z) >0 \non
\eeqa
k\"onnen wir unter Erhaltung der Differenzierbarkeit
$*H\wedge \tilde{E}$ in $\ol{U_N(R)}$ zu Null fortsetzen,
und es gilt dann
$\iota^*(*H\wedge \tilde{E})=0$ f\"ur die Inklusion
$\iota:S_N(R)\rightarrow \ol{U_N(R)} $. Aus
$(\ref{dint})$ und $(\ref{dwedge})$ folgt dann
\beqa
 0&=& \int_{U_N(R)}\sigma_q\d(* H\wedge \tilde{E})\non\\
  &=& \spu{H}{\rot \tilde{E}}{q,Z} +\spu{\div H}{\tilde{E}}{q-1,Z} \non\\
 &=&
 \spu{\tau H}{\tau\rot \tilde{E}}{\calL^q}
  +\spu{\rho H}{\rho\rot \tilde{E}}{\calL^{q-1}} \non\\
  &&+\spu{\tau\div H}{\tau \tilde{E}}{\calL^{q-1}}
   +\spu{\rho\div H}{\rho \tilde{E}}{\calL^{q-2}} \ .\non
\eeqa
Aus $(\ref{rrotddiv})$ folgen mit $e_2=0$ bzw. $e_1=0$ die
Behauptungen i) bzw. ii) im Falle glatter Formen $e_1$ und
$e_2$. \new
Andernfalls
betrachten wir zun\"achst den ersten Term in $(\ref{walzee})$.
Der lokalen Darstellung der Abbildung $\tau$ entnehmen wir
\beqa
 \gibt{c>0} \alle{r\in (0,R) } \dist(\supp \tau H(r),\gamma_2)
  \ge c>0 \ .\non
\eeqa
Somit existiert eine Abbildung $\chi\in C_\infty(\ol{B})$ mit
\beqa
 \alle{r\in (0,R)} \chi=1 \mbox{ in } \supp \tau H(r)\ ,
 \ \dist(\supp\chi,\gamma_2)>0 \ .\non
\eeqa
Wir erhalten
\beqa
\spu{\tau H}{\hat{M}^{-1}\vp\rrote e}{\calL^q}
=\spu{\tau H}{\hat{M}^{-1}\vp\rrote \chi e}{\calL^q} \ .\non
\eeqa
Aus Lemma \ref{supprandver} und $(\ref{lemm2})$ folgt,
da{\ss} wir $\chi e$ in $R^{q-1}(B)$ durch eine Folge aus
$\cnull{\infty}{q-1}(B)\subset C_\infty^{q-1,\gamma_1}(\ol{B})$
approximieren k\"onnen. Lemma \ref{approx} liefert
schlie{\ss}lich die Konvergenz im Raum $\calL^q$.
Die anderen Terme in $(\ref{walzee})$ sowie die Terme in
$(\ref{walzez})$
k\"onnen wir genauso behandeln.
\abstandk
Im Fall iii) und iv) gehen wir analog vor:
F\"ur
$h_1\in C_\infty^{q+1,\gamma_2}(\ol{B})$,
$h_2\in C_\infty^{q,\gamma_2}(\ol{B})$ und
$\tilde{H}:=\chtau\vp h_1+\chrho\vp h_2$ erhalten wir
\beqa
  0&=& \spu{\tau E}{\tau\div\tilde{H} }{\calL^q}
  +\spu{\rho E}{\rho\div \tilde{H}}{\calL^{q-1}} \non\\
  &&+\spu{\tau\rot E}{\tau \tilde{H}}{\calL^{q+1}}
  +\spu{\rho\rot E}{\rho \tilde{H}}{\calL^{q}}  \ .\non
\eeqa
Setzen wir $h_2=0$ bzw. $h_1=0$, so folgen iii) bzw. iv) im
Fall glatter Formen $h_1$ und $h_2$. Das gleiche Argument wie
oben liefert dann die Behauptung.\new
\qed
\section{Die kompakte Einbettung} \label{kompeinb}
\markboth{DIE KOMPAKTE EINBETTUNG}{DIE KOMPAKTE EINBETTUNG}
Wir definieren rekursiv die Gebiete, f\"ur welche wir
die kompakte Einbettung zeigen wollen:
\begin{defi}    \formel{zgebiet} \new
 Ein S--Gebiet $(S,\Gamma_1,\Gamma_2)$
 hei{\ss}t
 Z--Gebiet\label{pxzgebiet} in $M$ (zul\"assig), falls um
 jedes $z$ aus der Trennmenge $\gamma$ eine
 \"Ubergangsrandkarte mit $(\ref{karte})$ existiert,
 so da{\ss} das Gebiet
 $(\tilde{S},\gamma_1,\gamma_2)$
 mit
 \beqa
  \tilde{S}&:=&\{x\in S_{N} \mit x_N<0\} \non\\
  \gamma_i&=&S_{N}\cap h(\Gamma_i\cap \ol{V}) \non\\
  h(\Gamma_i\cap {V})&=&\{ts \mit t\in(0,1), s\in\gamma_i\}\non
 \eeqa
 ein $Z$--Gebiet in $S_N$ ist. \new
 Ein Z--Gebiet
 $(S,\Gamma_1,\Gamma_2)\in\calm(S_2)$ ist
 der untere Halbkreis
 $S:=\{x\in S_2 \mit x_2<0\} $ und
 $\Gamma_1\cup\Gamma_2=\{(1,0),(-1,0)\}$. Hier ist die Trennmenge
 $\gamma$ leer.
\end{defi}
Aus Satz \ref{dichtesatz} folgt, da{\ss} Z--Gebiete
$(S,\Gamma_1,\Gamma_2)$ stets die Approximationseigenschaft
$R^{q,\Gamma_1}(S)=\hat{R}^{q,\Gamma_1}(S)$ erf\"ullen.
\begin{satz}
 \formel{einistkompakt} \new
 In Z--Gebieten ist f\"ur jede zul\"assige
 Transformation $\eps$ die Einbettung
 \beqa
  R^{q,\Gamma_1}(S)\cap \eps^{-1}D^{q,\Gamma_2}(S)\hookrightarrow
  L_2^q(S) \non
 \eeqa
 kompakt.
\end{satz}
Um diese Aussage zu beweisen, gehen wir
wie in \cite{wwp} vor und f\"uhren eine vollst\"andige
Induktion \"uber die Raumdimension durch.
Wir werden sehen (Lemma \ref{wwp15}),
da{\ss} wir nach den Eigenformen des Maxwelloperators
entwickeln k\"onnen, sofern die Einbettung in
Satz \ref{einistkompakt} kompakt ist.
Gilt dieses Entwicklungsresultat
in $(N-1)$--dimensionalen Z--Gebieten
$(S,\gamma_1,\gamma_2)$ mit $S:=\{x\in S_N \mit x_N<0\}$,
so zeigt Lemma \ref{lemma9}, da{\ss} beschr\"ankte Familien
aus $\hat{R}^{q,\Gamma_1}(Z) \cap
\hat{D}^{q,\Gamma_2}(Z)$ f\"ur alle $R<1$
in $L_2^q(Z_R)$ relativ kompakt sind,
wobei
$(Z,\Gamma_1,\Gamma_2)=C(S,\gamma_1,\gamma_2)$,
$Z_R:=C_R(S)$.
Nachdem wir in Lemma \ref{indanfang} die kompakte
Einbettung in eindimensionalen Z--Gebieten gezeigt haben,
erhalten wir schlie{\ss}lich durch Lokalisieren die Aussage
f\"ur $N$--dimensionale Z--Gebiete
(Beweis von Satz \ref{einistkompakt}). Zun\"achst
zitieren wir
einige Hilfsmittel aus \cite{wwp}:
\begin{defi} \mbox{ } \new
Seien $H_\pm$\label{pxhpm} zwei Hilbertr\"aume mit
Normen $|\cdot|_\pm$\label{pxhpmnorm} und Skalarprodukten
$\spu{\cdot}{\cdot}{_\pm} $\label{pxhpmspu}.
\begin{enumerate}
\item[i)]
 F\"ur
 zwei dicht definierte Operatoren
 \beqa
  {A_\pm}:{D(A_\pm)\subset H_\pm}\longrightarrow {H_\mp} \non
 \eeqa
 nennen wir $(A_+,A_-)$\label{pxapm} ein duales Paar in $(H_+,H_-)$, falls
 $A_\pm^*=A_\mp$.
 \item[ii)] Ein duales Paar $(A_+,A_-)$ hat die
 Kompaktheitseigenschaft, wenn
 die Einbettungen $D(A_\pm)\cap \ol{R(A_\mp)}$
 mit den Graphennormen  nach
 $\ol{R(A_\mp)}$ mit $|\cdot|_\pm$ kompakt sind.
 Ist nur eine Einbettung kompakt, so hat
 $(A_+,A_-)$ die partielle Kompaktheitseigenschaft.
\end{enumerate}
\end{defi}
\begin{lemma}    \formel{wwp8} \new
   F\"ur ein duales Paar $(A_+,A_-)$ in $(H_+,H_-)$ und
   topologische Isomorphismen
   \beqa
        T_\pm:\tilde{H}_\pm\ra H_\pm \non
    \eeqa
   ist $(T_-^{-1}A_+T_+,T_+^*A_-(T_-^{-1})^*)$ ein duales
   Paar in $(\tilde{H}_+,\tilde{H}_-)$.
\end{lemma}
\begin{lemma} \formel{wwp9} \new
 Hat das duale Paar $(A_+,A_-)$ in $(H_+,H_-)$ die partielle
 Kompaktheitseigenschaft, so hat es auch die
 Kompaktheitseigenschaft, und es gelten:
 \begin{enumerate}
  \item[i)] $R(A_\pm) =\ol{R(A_\pm)}$ und
   $H_\pm=N(A_\pm)\oplus R(A_\mp)$.
   \item[ii)] Es existieren Zahlen $c_\pm>0$ mit
   \beqa
   \alle{x\in D(A_\pm)\cap R(A_\mp)}
     |x|_\pm\le c_\pm |A_\pm x_\pm|_\mp \ .\non
   \eeqa
   \item[iii)] Es existieren Folgen (evtl. auch endliche oder
   leere Folgen)
   $(\lambda_k)$ in $(0,\infty)$ und
   $(\vp_k^{\pm})\subset D(A_\pm)$ mit
   \begin{enumerate}
    \item[a)] $\lambda_k\ra\infty$
    (im Falle einer unendlichen Folge)
    \item[b)] $\{ \vp_k^{\pm}\mit k\in\nnn\}$ ist
     ein vollst\"andiges Orthonormalsystem in $R(A_\mp)$
    \item[c)] $\spu{A_\pm\vp_k^{\pm}}{A_\pm u}{\mp}
    =\lambda_k\spu{\vp_k^{\pm}}{u}{\pm}$ f\"ur
    $k\in\nnn$ und $u\in D(A_\pm)$
    \item[d)] $\vp_k^{\mp}=\lambda_k^{-1/2}A_\pm\vp_k^{\pm}$.
   \end{enumerate}
    \end{enumerate}
\end{lemma}
Der einfacheren Lesbarkeit wegen
f\"uhren wir den
auf funktionalanalytischen Grundlagen
beruhenden Beweis.
\abstandk
\beweis
Wir setzen
\beqa
  D_\pm:=D(A_\pm)\ ,\ R_\pm:=R(A_\pm)\ ,
  \ N_\pm:=N(A_\pm)\ ,
 \   X_\pm:=D_\pm\cap \ol{R_\mp} \non\\
 \spu{x}{y}{X_\pm} :=
  \spu{A_\pm x}{A_\pm y}{\mp} +
  \spu{x}{y}{\pm} \ .\non
\eeqa
Aus $A_{\pm}^*=A_{\mp}$ folgen die Zerlegungen
\beqa
 H_\pm=N_\pm \oplus \ol{R_\mp} \ .\formel{wwpzerl}
\eeqa
Sei die Einbettung
\beqa
 X_+ \hookrightarrow
 \ol{R_-} \formel{wwpkomp}
\eeqa
kompakt. Finden wir keine Konstante $c_+$, so da{\ss}
die Absch\"atzung
\beqa
 \alle{x_+\in X_+}|x_+|\le c_+|A_+x_+|_- \label{peterab}
\eeqa
 erf\"ullt ist, gibt es eine Folge
$(x^+_k)\subset X_+$ mit
$|x^+_k|=1$ und $A_+x^+_k\ra 0$. Diese enth\"alt
wegen der kompakten Einbettung $(\ref{wwpkomp})$
eine konvergente Teilfolge mit Grenzwert $x_+\in \ol{R_-}$.
Aus der Abgeschlossenheit von $A_+$ folgt
$x_+\in N_+$. Wegen $(\ref{wwpzerl})$ verschwindet
der Grenzwert, im Widerspruch zur Stetigkeit der Norm.
\abstandk
Dies liefert aber auch die Abgeschlossenheit von $R_+$:
Ist $(y^-_k)\subset R_+$ eine Folge mit Grenzwert $y_-$ in $H_-$,
so existiert o.B.d.A. eine Folge
$(x^+_k)\subset X_+ $ mit $A_+x^+_k=y^-_k$.
Wegen der Absch\"atzung $(\ref{peterab})$ und der kompakten
Einbettung $(\ref{wwpkomp})$ besitzt $(x^+_k)$
eine konvergente Teilfolge mit Grenzwert $x_+$ in $\ol{R_-}$.
Die Abgeschlossenheit von $A_+$ impliziert
$x_+\in X_+$ und $A_+x_+=y_-$.
Damit haben wir bis auf $R_-=\ol{R_-}$
die Aussagen i) und ii) im Fall ''$+$'' gezeigt.
\abstandk
Nun zu iii):
Nach ii) liefert der Satz von Lax--Milgram
zu allen $y_+\in \ol{R_-}$ eine eindeutige L\"osung $x_+\in X_+$
von
\beqa
 \alle{ z_+\in X_+}
 \spu{A_+x_+}{A_+z_+}{-} & =&\spu{y_+}{z_+}{+}
  \formel{wwpgl3} \\
  |x_+|_{X_+}&\le& c|y_+|_+ \ .\formel{wwpab3}
\eeqa
Wegen
$D_+=N_+\oplus X_+$ gilt $(\ref{wwpgl3})$
f\"ur alle ${z_+\in D_+}$, und wir erhalten
$A_+x_+\in X_-$ und $A_-A_+x_+=y_+$.
Der L\"osungsoperator
\beqa
 \Abb{L}{\ol{R_-}\subset H_+}{\ol{R_-}\subset H_+}{y_+}{x_+} \ ,\non
\eeqa
der $y_+$ die L\"osung $x_+$ von $(\ref{wwpgl3})$
zuordnet, ist nach $(\ref{wwpab3})$ und Voraussetzung
kompakt und erf\"ullt $N(L)=0$.
Wegen
\beqa
 \spu{y_+}{Ly_+}{+} =\spu{A_-A_+x_+}{x_+}{+}
 =|A_+x_+|_+^2 \in [0,\infty) \non
\eeqa
f\"ur $y_+\in \ol{R_-}$ und $x_+=Ly_+$
ist $L$ selbstadjungiert und positiv.
Der Spektralsatz liefert
eine monoton wachsende (evtl. auch endliche oder leere) Folge
$(\lambda_k)\subset (0,\infty)$
und ein vollst\"andiges Orthonormalsystem $\{\vp_k^+\}$
in $\ol{R_-}$ mit
$\lambda_k\ra\infty$ f\"ur $k\ra\infty$
(im Falle einer unendlichen Folge) und
$L\vp_k^+=\lambda_k^{-1}\vp_k^+$.
F\"ur alle $u_+\in D_+$ gilt
\beqa
 \spu{A_+\vp^+_k}{A_+u_+}{-}
 =\spu{\lambda_k A_+L\vp^+_k}{A_+u_+}{-}
 =\lambda_k\spu{\vp^+_k}{u_+}{+} \ .\non
\eeqa
Wir definieren $\vp^-_k:=\lambda_k^{-1/2}A_+\vp^+_k\in X_-$,
und erhalten mit $u_-\in D_-$
\beqa
 A_-\vp^-_k&=&\lambda_k^{1/2}A_-A_+L\vp^+_k=
 \lambda_k^{1/2}\vp^+_k
 \non\\
 \spu{A_-\vp^-_k}{A_-u_-}{+} &=&
  \lambda_k^{1/2}\spu{\vp^+_k}{A_-u_-}{+}
  =\lambda_k\spu{\vp^-_k}{u_-}{-} \ .\non
\eeqa
F\"ur $x_-\in R_+$, $x_-=A_+x_+$ mit o.B.d.A. $x_+\in X_+$
gilt
\beqa
  \spu{x_-}{\vp^-_k}{-} =
  \spu{A_+x_+}{\vp^-_k}{-} =
   \lambda_k^{1/2}\spu{x_+}{\vp^+_k}{+} \ .\formel{peter5}
\eeqa
Setzen wir $x_-:=\vp^-_l,x_+:=\lambda_l^{-1/2}\vp^+_l$,
so folgt hieraus,
da{\ss} $\{\vp^-_k\}$ ein Orthonormalsystem ist.
Andererseits impliziert $(\ref{peter5})$
aber auch die Vollst\"andigkeit
(wegen der Vollst\"andigkeit von $\{\vp^+_k\}$ in $\ol{R_-}$).
Damit ist iii) gezeigt.
\abstandk
Ist $(y^-_k)$ eine beschr\"ankte Folge in $X_-$,
so besitzt die Folge $(x^+_k)$ mit
$A_+x^+_k=y^-_k$ wegen ii) und der Voraussetzung
eine konvergente Teilfolge, die wir
wieder mit $(x^+_k)$ bezeichnen.
Wir setzen $x=x^+_{k}-x^+_l$, $y$ analog
und erhalten mit iii) und der Schwarzschen Ungleichung
\beqa
 \spu{y}{y}{-}&=&\spu{y}{A_+x}{-} \non\\
    &=&\sum_{n=1}^\infty\spu{y}{\vp^-_n}{-}
     \cdot\spu{\vp^-_n}{A_+x}{-} \non\\
    &=&\sum_{n=1}^\infty\spu{y}{A_+\vp^+_n}{-}
     \cdot\spu{\vp^+_n}{x}{+} \non\\
    &\le& |x|_+\cdot|A_-y|_+\non\\
   &\le& |x^+_k-x^+_l|_+\cdot(|A_-y^-_k|_++|A_-y^-_l|_+)\ra 0 \ .\non
\eeqa
Dies liefert die Kompaktheit der Einbettung
$X_-\hookrightarrow R_+$,  und analog zum Fall ''$+$'' folgen
i) und ii) f\"ur den Fall ''$-$''.\qed
\abstandk
Weitere Resultate aus \cite{wwp} k\"onnen wir fast direkt
\"ubertragen:
In einem Gebiet $(S,\Gamma_1,\Gamma_2)\in\calm(M)$ bilden
die Operatoren
\beqa
 \Abb{\rot^q}{R^{q,\Gamma_1}(S)}{L_2^{q+1}(S)}{E}{\rot E} \non\\
  \Abb{\div^{q+1}}{\hat{D}^{q+1,\Gamma_2}(S)}{L_2^{q}(S)}{E}{\div E} \non
\eeqa
ein duales Paar in $(L_2^{q}(S),L_2^{q+1}(S))$.
Seien $\sigma,\mu,\eps$ zul\"assige Transformationen.
Bezeichnen wir mit $L_{2,\mu}^q(S)$\label{pxlzmu} den Raum $L_2^q(S)$,
versehen mit\label{pxlzmuspu}
$\spu{\cdot}{\cdot}{\mu,q,S} :=\spu{\mu\cdot}{\cdot}{q,S} $,
so sind die Abbildungen\label{pxjsigmu}
\beqa
 \Abb{\eps}{L_2^q(S)}{L_2^q(S)}{E}{\eps E} \non\\
 \Abb{j_{\sigma,\mu}}{L_{2,\sigma}^q(S)}{L_{2,\mu}^q(S)}{E}{E}
 \non
\eeqa
topologische Isomorphismen und besitzen die Adjungierten
\beqa
 \eps^*=\eps\ ,\ (j_{\sigma,\mu})^*=
  j_{\id,\sigma}\sigma^{-1}\mu j_{\mu,\id}  \ .\non
\eeqa
Nach Lemma \ref{wwp8} bilden
\beqa \left.\begin{array}{rcl}
  (\i\div^{q}\eps j_{\eps,\id}&,&\i j_{\id,\eps} \rot^{q-1}) \\
  (\i j_{\id,\mu}\mu^{-1} \rot^q j_{\eps,\id}&,&\i j_{\id,\eps}
   \eps^{-1} \div^{q+1}j_{\mu,\id})   \\
  (\i\rot^{q+1}\mu j_{\mu,\id}&,&\i j_{\id,\mu}\div^{q+2})
  \end{array} \right\} \formel{dualpaar}
\eeqa
duale Paare in $(L_{2,\eps}^{q}(S),L_2^{q-1}(S))$,
$(L_{2,\eps}^{q}(S),L_{2,\mu}^{q+1}(S))$ bzw.
$(L_{2,\mu}^{q+1}(S),L_{2}^{q+2}(S))$.
Aus $H=N(A)\oplus\ol{R(A^*)}$ (orthogonal) f\"ur einen
dicht definierten Operator $A:H\ra \tilde{H}$
und aus
$\ol{\eps^{-1}\div\hat{D}^{q+1,\Gamma_2}}(S)\subset
\eps^{-1}\hat{D}_0^{q,\Gamma_2}(S)$ (nach $(\ref{rotrotistnull})$)
folgt
\begin{lemma} \formel{orthzerl} \new
Es gelten die orthogonalen Zerlegungen
\beqa
 L_{2,\eps}^q(S) &=& \eps^{-1}\hat{D}_0^{q,\Gamma_2}(S)
 \oplus \ol{\rot R^{q-1,\Gamma_1} (S)}
  =R_0^{q,\Gamma_1}(S) \oplus
  \ol{\eps^{-1}\div \hat{D}^{q+1,\Gamma_2}(S)} \non\\
 L_{2,\eps}^{q}(S) &=& \ol{\rot R^{q-1,\Gamma_1}(S)}
  \oplus (\eps^{-1}\hat{D}_0^{q,\Gamma_2}(S)
  \cap R_0^{q,\Gamma_1}(S)) \oplus
  \ol{\eps^{-1}\div\hat{D}^{q+1,\Gamma_2}(S)} \non\\
 L_{2,\mu}^{q+1}(S)&=&
   \hat{D}_0^{q+1,\Gamma_2}(S)\oplus
   \ol{\mu^{-1}\rot R^{q,\Gamma_1}(S)}
   =\ol{\div\hat{D}^{q+2,\Gamma_2}(S)} \oplus
   \mu^{-1}R_0^{q+1,\Gamma_1}(S) \non\\
 L_{2,\mu}^{q+1}(S)&=&
  \ol{\mu^{-1}\rot R^{q,\Gamma_1}(S)}\oplus
  (\hat{D}_0^{q+1,\Gamma_2}(S) \cap
   \mu^{-1}R_0^{q+1,\Gamma_1}(S)) \oplus
   \ol{\div\hat{D}^{q+2,\Gamma_2}(S)}  \ .\non
\eeqa
\end{lemma}
\begin{lemma} \formel{wwpc1} \new
 Die Kompaktheitseigenschaft der dualen Paare in
 $(\ref{dualpaar})$ h\"angt nicht von $\eps$ und
 $\mu$ ab.
\end{lemma}
\beweis Sei $q$ beliebig.
Wir betrachten die Abbildungen
\beqa
 \Abb{A_+}{j_{\id,\eps}\eps^{-1}\hat{D}^{q,\Gamma_2}(S)
  \subset L_{2,\eps}^q(S)}{L_2^{q-1}(S)}{E}{\i\div\eps j_{\eps,\id} E}
  \non\\
  \Abb{A_-}{R^{q-1,\Gamma_1}(S)\subset L_2^{q-1}(S)}{L_{2,\eps}^q(S)}{
   E}{\i j_{\id,\eps}\rot E} \ ,\non
\eeqa
wobei wir sowohl hier als auch im folgenden den Index
$q$ bzw. $q-1$ bei $\div$ bzw. $\rot$
wieder fortlassen wollen.
Eine in der Graphennorm des Operators $A_-$ beschr\"ankte Folge
$(E_k)\subset D(A_-)\cap \ol{R(A_+)}$ ist eine Folge
in $R^{q-1,\Gamma_1}(S)\cap \ol{\div\hat{D}^{q,\Gamma_2} }(S)$,
f\"ur die der Ausdruck
\beqa
 \normu{E_k}{q-1,S} +(\spu{\eps j_{\id,\eps}\rot E_k}{
 j_{\id,\eps} \rot E_k}{\eps,q,S} )^{1/2} \formel{qwer}
\eeqa
beschr\"ankt ist.
Dies ist genau dann der Fall, wenn f\"ur den Ausdruck
$(\ref{qwer})$ f\"ur $\eps = \id$ eine Schranke existiert.
Damit ist die Frage, ob $(E_k)$ eine in
$\ol{R(A_+)}=\ol{\div\hat{D}^{q,\Gamma_2} }(S)$ konvergente
Teilfolge besitzt, unabh\"angig von $\eps$. Lemma \ref{wwp9}
liefert dann die Kompaktheitseigenschaft f\"ur das
duale Paar.
Analog k\"onnen wir die anderen dualen Paare behandeln.
\qed
\begin{lemma} \formel{wwp11} \new
  Seien $\eps,\mu$ zul\"assig und
  \beqa
  X&:=&R^{q,\Gamma_1}(S)\cap \eps^{-1}\hat{D}^{q,\Gamma_2}(S)\non\\
  \normu{E}{X,q} &:=&\normu{E}{q,S} +\normu{\div\eps E}{q-1,S}
  + \normu{\rot E}{q+1,S} \ .\non
 \eeqa
 Dann sind \"aquivalent:
  \begin{enumerate}
  \item[i)] Die Einbettung $X\hookrightarrow L_2^q(S)$ ist
   kompakt.
  \item[ii)] Die Einbettungen
  $\ol{\rot R^{q-1,\Gamma_1}(S)} \cap \eps^{-1}\hat{D}^{q,\Gamma_2}
  (S)$ und $\ol{\eps^{-1}\div \hat{D}^{q+1,\Gamma_2}(S)}
  \cap R^{q,\Gamma_1}(S)$ mit $\normu{\cdot}{X,q} $ nach
  $L_2^q(S)$ sind kompakt und die
  Dirichlet--Neumann--Felder
  ${\cal H}^q:=R^{q,\Gamma_1}_0(S)\cap \hat{D}_0^{q,\Gamma_2}(S)$
  endlich dimensional.
\end{enumerate}
\end{lemma}
\beweis Analog zu \cite{wwp}. Ebenso:
\begin{lemma} \formel{wwp15} \new
Gelten die Voraussetzungen von Lemma \ref{wwp11}
  und sei die Einbettung
  $X\hookrightarrow L_2^q(S)$ kompakt.
  Dann gibt es ein $N^q\in\nnn\cup \{0\}$\label{pxnq}, ein
 bzgl. $\spu{\eps\cdot}{\cdot}{q,S}$ bzw.
 $\spu{\mu\cdot}{\cdot}{q+1,S}$ vollst\"andiges  Orthonormalsystem
 $\{E_n^q, n\in\nnn\}$\label{pxenq} von
 $\eps^{-1} \hat{D}_0^{q,\Gamma_2}(S)$ bzw.
 $\{H_n^q, n>N^q\}$\label{pxhnq}
 von $\ol{\mu^{-1}\rot R^{q,\Gamma_1}(S) }$
 sowie eine
 Folge $\{\omega_n^q,n>N^q\}\subset (0,\infty)$,
 $\omega_n^q\stackrel{n\ra\infty}{\ra} \infty $ mit
 \beqa
  \mfur n\le N^q && E_n^q\in\rq_0(S) \cap
  \eps^{-1}\hat{D}_0^{q,\Gamma_2}(S)\non\\
  \mfur n>N^q &&
  (E_n^q,H_n^q)\in\rq(S)\times\hat{D}^{q+1,\Gamma_2}(S)\non\\
   &&\rot E_n^q+\i\omega_n^q\mu H_n^q =0\non\\
    &&\div H_n^q+\i\omega_n^q\eps E_n^q=0\
    \ .\formel{eigenloesung}
 \eeqa
 Im Falle endlicher oder leerer Orthonormalsysteme
 sind die Bezeichnungen
 entsprechend zu \"andern.
\end{lemma}
Letzteres wollen wir so auch ohne weiteren Kommentar
in Zukunft handhaben. Treten z.B. Reihen
\"uber leere bzw. endliche Orthonormalsysteme auf,
so sind diese durch 0 bzw. Summen zu ersetzen.
\abstandk
Wir bereiten das folgende Lemma vor und
erinnern an die in Abschnitt \ref{chkegelspitzen}
eingef\"uhrten Bezeichnungen.
Seien $\cals:=(\tilde{S},\gamma_1,\gamma_2)\in\calm_D(S_N)$
und $Z:=C_R(\cals)$ f\"ur ein $R\in(0,\infty)$.
Gelte ferner f\"ur $\cals$,
$\eps=\id$, $\mu=\id$ und alle $q$ die
Aussage von Lemma \ref{wwp15}. Wir zerlegen
$E\in\hat{R}^{q,\Gamma_1}(Z)
 \cap \hat{D}^{q,\Gamma_2}(Z)$ nach $(\ref{wwsph})$
orthogonal in $E=\chrho\rho E+\chtau\tau E$.
Nach Lemma \ref{orthzerl} und Voraussetzung
(beachte $\cals\in\calm_D(S_N)$)
k\"onnen wir $\rho E(r)$ und $\tau E(r)$ in Fourierreihen
entwickeln, so da{\ss}
\beqa
  E&=&\chrho \sum_{n\ge 1} a_n(r)E_n^{q-1}
  +\chrho \sum_{n>N^{q-2}} b_n(r)H_n^{q-2}\non\\
  &&+\chtau \sum_{n\ge 1} c_n(r)E_n^{q}
 +\chtau \sum_{n>N^{q-1}} d_n(r)H_n^{q-1} \formel{ezerlegung}
\eeqa
mit
\beqa \begin{array}{cclccl}
 a_n(r)&:=&\spu{\rho E(r)}{E_n^{q-1}}{q-1,\tilde{S}} &
 b_n(r)&:=&\spu{\rho E(r)}{H_n^{q-2}}{q-1,\tilde{S}} \\
 c_n(r)&:=&\spu{\tau E(r)}{E_n^{q}}{q,\tilde{S}} &
 d_n(r)&:=&\spu{\tau E(r)}{H_n^{q-1}}{q,\tilde{S}}
 \,.\end{array} \non\\
 \formel{koeffizienten0}
\eeqa
(Diese Definitionen gelten nur f\"ur solche $n$, f\"ur die
die zweite Komponente im Skalarprodukt erkl\"art ist;
vgl. Lemma \ref{wwp15}).
Analog:
 \beqa \begin{array}{cclccl}
 a_n^D(r)&:=&\spu{\rho \div E(r)}{E_n^{q-2}}{q-2,\tilde{S}} &
 b_n^D(r)&:=&\spu{\rho \div E(r)}{H_n^{q-3}}{q-2,\tilde{S}} \\
 c_n^D(r)&:=&\spu{\tau \div E(r)}{E_n^{q-1}}{q-1,\tilde{S}} &
 d_n^D(r)&:=&\spu{\tau \div E(r)}{H_n^{q-2}}{q-1,\tilde{S}} \\
 a_n^R(r)&:=&\spu{\rho \rot E(r)}{E_n^{q}}{q,\tilde{S}} &
 b_n^R(r)&:=&\spu{\rho \rot E(r)}{H_n^{q-1}}{q,\tilde{S}} \\
 c_n^R(r)&:=&\spu{\tau \rot E(r)}{E_n^{q+1}}{q+1,\tilde{S}} &
 d_n^R(r)&:=&\spu{\tau \rot E(r)}{H_n^{q}}{q+1,\tilde{S}}
 \end{array} \non\\
 \formel{koeffizienten}
\eeqa
Wir definieren $\call_2(\rho,R)$\label{pxcall}
als den Raum der Folgen
$(u_n)$ me{\ss}barer Funktionen mit
\beqa
 r^\rho |u_n|^2&\in&L_1(0,R) \non\\
 \normu{(u_n)}{{\call}_2(\rho,R)} &:=&
  (\sum_{n}\int_0^R r^\rho |u_n(r)|^2 \dr )^{1/2}<\infty\ .\non
\eeqa
Versehen mit\label{pxspucall}
\beqa
  \spu{(u_n)}{(v_n)}{{\call}_2(\rho,R)} &:=&
  \sum_{n}\int_0^R r^\rho u_n(r)\ol{v_n(r)} \dr \ ,\non
\eeqa
wird dieser zu einem Hilbert\-raum.
Treffen wir zus\"atzlich die Konvention, da{\ss}
die Folgen aus $\call_2$ je nach Kontext bei 1 oder
$N^q+1$ usw. starten (vgl. wieder Lemma \ref{wwp15}),
so sind die Abbildungen
\beqa  \left.
 \begin{array}{ccccccccc}
 A&:&{\call}_2({N-1},R) &\ra &\calL^{q-1}&,&(a_n)&
 \ra&\sum_{n\ge 1} a_n(r)E_n^{q-1}   \\
B&:&{\call}_2({N-1},R) &\ra &\calL^{q-1}&,&(b_n)&
 \ra&\sum_{n>N^{q-2}} b_n(r)H_n^{q-2}   \\
 C&:&{\call}_2({N-1},R) &\ra &\calL^{q}&,&(c_n)&
 \ra&\sum_{n\ge 1} c_n(r)E_n^{q}   \\
 D&:&{\call}_2({N-1},R) &\ra &\calL^{q}&,&(d_n)&
 \ra&\sum_{n>N^{q-1}} d_n(r)H_n^{q-1}
\end{array}  \right\}              \formel{abcd}
\eeqa
wohldefiniert und isometrisch.
Letzteres folgt z.B. f\"ur die Abbildung $A$ aus
\beqa
 \normu{A(a_n)}{\calL^{q-1}} ^2
  &=&\int_0^R r^{N-1} \sum_{n\ge 1} |a_n|^2 dr\non\\
  &=&\sum_{n\ge 1}\int_0^R r^{N-1}|a_n|^2 dr \non\\
  &=&\normu{(a_n)}{\call_2({N-1},R)} ^2 \ .\formel{labelneu}
\eeqa
\begin{lemma} \formel{dgl} \new
 F\"ur das Gebiet
   $\cals:=(\tilde{S},\gamma_1,\gamma_2)\in\calm_D(S_N)$,
    $\eps=\id$, $\mu=\id$ und alle $q$ gelte
    die Aussage von Lemma  \ref{wwp15}.
 Sei ferner $(Z,\Gamma_1,\Gamma_2):=C_R(\cals)$ f\"ur ein
 $R\in(0,\infty)$.
 Dann erf\"ullen die
 Koeffizienten von
 $E\in\hat{R}^{q,\Gamma_1}(Z)
 \cap \hat{D}^{q,\Gamma_2}(Z)$ aus
 $(\ref{koeffizienten0})$ und $(\ref{koeffizienten})$
 (sofern definiert)
 \beqa
  b_m(r)&=&-\i\ommeg{q-2} r a_m^D(r)\non\\
  (r^{N-q}b_m(r))'&=&r^{N-q}d_m^D(r)\non\\
  c_m(r)&=&\i\ommeg{q} rd_m^R(r)\non\\
  (r^qc_m(r))'&=&r^qa_m^R(r)\non\\
  (r^qd_m(r))'&=&-\i\omeg{q-1} r^{q-1}a_m(r)+r^qb_m^R(r)\non\\
  (r^{N-q}a_m(r))'&=&r^{N-q}c_m^D(r) +\left\{
   \begin{array}{ccl}
     \i\omeg{q-1} r^{N-q-1}d_m(r)&\mfur &m> N^{q-1}\\
     0&\mfur&m\le N^{q-1}
   \end{array}
   \right. \non
  \eeqa
\end{lemma}
\beweis
F\"ur $\psi\in \cnull{\infty }{\ }(I)$ gilt
nach $(\ref{eigenloesung})$ und $(\ref{walzez})$
\beqa
 \int_0^R r^{N-1}b_m(r){\psi}(r) dr
 &=& \spu{(b_m)}{(\ol{\psi}\delta_{nm})}{\call_2({N-1},R)} \non\\
 &=&\spu{\rho E}{\ol{\psi} H_m^{q-2}}{\calL^{q-1}}\non\\
 &=&-\i\ommeg{q-2} \spu{\rho E}{\rrot \ol{\psi} E_m^{q-2}}{\calL^{q-1}}
  \non\\
 &=&-\i\ommeg{q-2} \spu{\rho \div E}{\hat{M}\ol{\psi} E_m^{q-2}}{\calL^{q-2}}
  \non\\
  &=&-\i\ommeg{q-2} \int_0^R r^{N-1} r a_m^D(r){\psi}(r) dr \ .
  \non
\eeqa
Analog erhalten wir mit $(\ref{walzee})$ bis $(\ref{walzed})$
\beqa
 \int_0^R r^{N-1}b_m(r)r^{-(q-1)}&&\hspace{-1.0cm}(r^{q-1}{\psi}(r))'dr
    =-\int_0^R r^{N-1}d_m^D(r){\psi}(r) dr \non\\
 \int_0^R r^{N-1} c_m(r){\psi}(r)dr
       &=&\i\ommeg{q} \int_0^R r^{N-1}r d_m^R(r){\psi}(r) dr
    \non\\
\int_0^R r^{N-1} c_m(r) r^{-(N-q-1)}&&\hspace{-1.0cm}
(r^{N-q-1}{\psi}(r))' dr
    =-\int_0^R r^{N-1} a_m^R(r){\psi}(r)dr  \non\\
\int_0^R r^{N-1} d_m(r)r^{-(N-q-1)}&&\hspace{-1.0cm}
(r^{N-q-1}{\psi}(r))'dr\non\\
  &=&\i\omeg{q-1} \int_0^Rr^{N-1}r^{-1}a_m(r){\psi}(r) dr \non\\
  &&-\int_0^R r^{N-1}b_m^R(r){\psi}(r)dr\non\\
  \int_0^R r^{N-1}a_m(r)r^{-(q-1)}&&\hspace{-1.0cm}(r^{q-1}{\psi}(r))'dr
  \non\\
  &=&
-\int_0^Rr^{N-1}c_m^D(r){\psi}(r)dr
 \non\\
 &&+
\left\{ \begin{array}{cl}
0 &\mfur m\le N^{q-1} \\
 -\i\omeg{q-1} \int_0^R r^{N-1}r^{-1}d_m(r){\psi}(r)dr&
  \mfur m>N^{q-1}\ .
\end{array}\right. \non
%
\eeqa
\qed
\begin{lemma} \formel{lemma9}   \new
 F\"ur das Gebiet
   $\cals:=(\tilde{S},\gamma_1,\gamma_2)\in\calm_D(S_N)$
     und alle $q$ gelte
    die Aussage von Lemma  \ref{wwp15} im Falle $\eps=\id$, $\mu=\id$.
 Seien ferner
 $(Z,\Gamma_1,\Gamma_2):=C(\cals)$ und
 $\{E^\alpha\}$ eine
 beschr\"ankte Familie in $\hat{R}^{q,\Gamma_1}(Z)
 \cap \hat{D}^{q,\Gamma_2}(Z)$.
 Dann ist
 f\"ur alle $R<1$ die Familie  $\{E^\alpha\}$ in $L_2^q(Z_R)$
 relativ kompakt, wobei $Z_R:=C_R(\tilde S )$.
\end{lemma}
\beweis
Wegen der Isometrie der Abbildungen $A,B,C,D$ in $(\ref{abcd})$
und $\chtau$, $\chrho$ in $(\ref{wwsphiso})$
sind die Bildmengen von $\chrho A,\cdots,\chtau D$ abgeschlossen.
Wenn wir zeigen, da{\ss} f\"ur alle $R<1$ die Familien
$\{a_m^{\alpha}\},\cdots,\{d_m^{\alpha}\}$ der Koeffizienten von
  $E^{\alpha}$
  aus
  $(\ref{koeffizienten0})$ in ${\call}_2({N-1},R)$ relativ kompakt
  sind, so gilt dies auch f\"ur deren Bilder
  $\{\chrho A(a_m^{\alpha})\},\cdots,
 \{\chtau D(d_m^{\alpha})\} $ in $L_2^q(Z_R)$.
  Aus $(\ref{ezerlegung})$ folgt dann die Behauptung.
\abstandk
i) Wir betrachten zun\"achst die Koeffizienten
$b_m^{\alpha}$ und definieren
\beqa
 {\gamma}&:=&-N+2q-1\ ,\ \mu_m:=\omeg{q-2}\ ,\non\\
 u_m^{\alpha}(r)&:=&r^{N-q}b_m^{\alpha}(r)\non\\
f_m^{\alpha}(r)&:=&r^{N-q}d_m^{D,\alpha}(r)\ .\non
\eeqa
Dann gelten
\beqa
\normu{(u_m^{\alpha})}{\call_2({{\gamma}},1) }
   &=&\normu{(b_m^{\alpha})}{\call_2({N-1},1) }
   \le \normu{E^{\alpha}}{q,Z}\le c \non\\
  \normu{(f_m^{\alpha})}{\call_2({{\gamma}},1)}
   &=&\normu{(d_m^{D,\alpha})}{\call_2({N-1},1) }
   \le \normu{\div E^{\alpha}}{q-1,Z}\le c\ .\non\\
    \sum_{m>N^{q-2}}\int_0^1\mu_mr^{{\gamma}}|u_m^\alpha(r)|^2 dr
    &=& \sum_{m>N^{q-2}}\int_0^1\ommeg{q-2} r^{N+1}|a_m^D(r)|^2
    dr\non\\
    &\le&\sup_{m>N^{q-2}}\{\ommeg{q-2} \}\sum_{m>N^{q-2}}
    \int_0^1r^{N-1}|a_m^D(r)|^2dr\non\\
    &\le& c\ ,\non
  \eeqa
 wobei wir Lemma \ref{dgl}  und
\beqa
 \normu{(b_m^\alpha)}{\call _2({N-1},1)}
  \le \normu{\rho E^\alpha}{\calL ^{q-1}}
  \le \normu{E^\alpha}{q,Z} \non
\eeqa
(nach $(\ref{labelneu})$ und $(\ref{wwsphiso})$)
investiert
 haben. Die Absch\"atzung f\"ur $f_m^\alpha$ folgt
 analog. Nach Lemma \ref{dgl} und
 \cite[Lemma 8]{weck} besitzt $(u_m^{\alpha})_\alpha$
   eine in $\call_2({{\gamma}},1)$
 und damit auch $(b_m^\alpha)_\alpha$
eine in
 $\call_2({N-1},1)$ konvergente Teilfolge.
 \abstandk
ii) Die Koeffizienten $c_m^{\alpha}$ k\"onnen wir
 mit  ${\gamma}:=N-2q-1$, $u_m^{\alpha}(r):=r^qc_m^{\alpha}(r)$,
 $f_m^{\alpha}(r):=r^{q}a_m^{R,\alpha}(r)$ und
 \beqa
   \mu_m:=\Big\{
\begin{array}{ccc}
 0&\mfur &m\le N^q\\
 \omeg{q} &\mfur &m>N^q
\end{array} \non
\eeqa
analog behandeln und gewinnen somit eine in
$\call_2({N-1},1)$ konvergente
Teilfolge von $(c_m^{\alpha})_\alpha$.
\abstandk
iii) F\"ur $m>N^{q-1}$ definieren wir $\gamma:=N-2q-1$,
 $\beta_m:=\omeg{q-1} $ und
\beqa  \begin{array}{rclcrcl}
  u_m^{\alpha}(r)&:=&r^{q}d_m^{\alpha}(r)
  &,&
  v_m^{\alpha}(r)&:=&\i r^{q-1}a_m^{\alpha}(r)\\
 f_m^{\alpha}(r)&:=&r^{q}b_m^{R,\alpha}(r)
  &,&
 g_m^{\alpha}(r)&:=&\i r^{q+1}c_m^{D,\alpha}(r)\ .
  \end{array} \non
\eeqa
Diese erf\"ullen
\beqa  \begin{array}{lclclc}
   \normu{(u_m^{\alpha})}{\call_2({\gamma},1) }
   &=&\normu{(d_m^{\alpha})}{\call_2({N-1},1) }
   &\le& c&,\\
  \normu{(v_m^{\alpha})}{\call_2({\gamma+2},1) }
   &=&\normu{(a_m^{\alpha})}{\call_2({N-1},1) }
   &\le&  c&,\\
 \normu{(f_m^{\alpha})}{\call_2({\gamma},1) }
   &=&\normu{(b_m^{R,\alpha})}{\call_2({N-1},1) }
   &\le&  c&, \\
  \normu{(g_m^{\alpha})}{\call_2({\gamma-2},1) }
   &=&\normu{(c_m^{D,\alpha})}{\call_2({N-1},1) }
   &\le&  c& .
  \end{array} \non
\eeqa
Nach Lemma \ref{dgl} gelten die Differentialgleichungen
\beqa
 (u_m^{\alpha})'(r)&=&(r^{q}d_m^{\alpha}(r))'
  =-\i\omeg{q-1} r^{q-1}a_m^{\alpha}(r)+r^{q}b_m^{R,\alpha}(r)\non\\
  &=&-\beta_m v_m^{\alpha}(r)+f_m^{\alpha}(r)\non\\
 r^{-\gamma}(r^{{\gamma+2}}v_m^{\alpha})'(r)
  &=&r^{-(N-2q-1)}(\i r^{N-q}a_m^{\alpha}(r))'\non\\
  &=&-\omeg{q-1}r^qd_m^{\alpha}(r)+\i r^{q+1}c_m^{D,\alpha}(r)\non\\
  &=&-\beta_m u_m^{\alpha}(r)+g_m^{\alpha}(r)\ .\non
\eeqa
  Nach \cite[Lemma 9]{weck} enth\"alt
 $(u_m^\alpha)_{\alpha}$ bzw. $(v_m^\alpha)_\alpha$
 f\"ur alle $R<1$
 eine in $\call_2(\gamma,R)$ bzw. $\call_2({\gamma+2},R)$
 konvergente Teilfolge. Das gleiche gilt dann auch f\"ur
 $((\id-\Pi_{N^{q-1}})(a_m^\alpha))_\alpha$ bzw.
 $(d_m^\alpha)_\alpha$
 in $\call_2({N-1},R)$, wobei
 \beqa
  \Abb{\Pi_K}{\call_2({\gamma},R)}{\call_2({{\gamma}},R)}{
   (w_n)}{(\hat{w}_n)} \ ,
   \ \mmit \hat{w}_n:=\Big\{
\begin{array}{cl}
 w_n&\mfur  n\le K\\
 0& \mbox{ sonst }\ .
   \end{array}\non
\eeqa
 F\"ur $m\le N^{q-1} $ definieren wir ( vgl. i))
 ${\gamma}:=-N+2q-1$ und
 \beqa
  u_m^\alpha(r)&:=&\Big\{
                     \begin{array}{cl}
                        r^{N-q}a_m^\alpha(r) &\mfur n\le N^{q-1}\\
                        0 &\mbox{ sonst } \end{array} \non\\
 f_m^\alpha(r)&:=&\Big\{
                     \begin{array}{cl}
                        r^{N-q}c_m^{D,\alpha}(r)&\mfur n\le N^{q-1}\\
                        0 &\mbox{ sonst } \end{array} \non\\
 \mu_m&:=&\Big\{   \begin{array}{cl}
                        0&\mfur n\le N^{q-1}\\
                        m &\mbox{ sonst } \ .\end{array} \non
 \eeqa
 Nach \cite[Lemma 8]{weck} enth\"alt $(u_m^\alpha)_\alpha$,
 somit auch $(\Pi_{N^{q-1}}a_m^\alpha)_\alpha$ eine konvergente
 Teilfolge.
\qed
\begin{lemma} \formel{wwp12} \new
Die Kompaktheit der Einbettung in Satz \ref{einistkompakt} ist
unabh\"angig von $\eps$, $\mu$
und invariant unter Diffeomorphismen.
\end{lemma}
\beweis Die erste Aussage folgt wie
in  \cite{wwp} aus Lemma \ref{wwp11}, Lemma \ref{wwpc1} und einer
\"Uberlegung, analog zu \cite[Remark 2]{picard3}.
Die zweite Aussage ist dann eine Konsequenz aus der ersten
und Lemma \ref{taudivrot}.\qed
\begin{lemma}  \formel{indanfang} \new
 In eindimensionalen Z--Gebieten
 $(S,\Gamma_1,\Gamma_2)$ gelten
 $R^{q,\Gamma_1}(S)=\hat{R}^{q,\Gamma_1}(S)$ und
 die Aussage
 von Satz \ref{einistkompakt}.
\end{lemma}
\beweis  Nach Lemma \ref{wwp12} und der Bemerkung
 \ref{bemerk2}
 gen\"ugt es, statt eindimensionaler Z--Gebiete
 Intervalle zu betrachten und
 $\eps=\id$ anzunehmen.
 Die erste Behauptung folgt dann mit einer
 geeigneten Abschneidefunktion  aus den entsprechenden
 Aussagen f\"ur $\hnull{1}{\ }(S) $ und $H_1(S)$.
 F\"ur die zweite Behauptung
 definieren wir den Raum
 $Y^q:=R^{q,\Gamma_1}(S)\cap D^{q,\Gamma_2}(S)$.
 Aus $(\ref{risth})$ folgt $Y^0,*Y^1\subset H_1(S)$, und
 der Rellichsche Auswahlsatz liefert das Gew\"unschte.
 \qed
\abstand
{\bf Beweis} von Satz \ref{einistkompakt}.
Es gen\"ugt wieder $\eps=\id$ anzunehmen. Wir f\"uhren eine
Induktion \"uber die Raumdimension durch. Den Induktionsanfang
entnehmen wir Lemma \ref{indanfang}. Gelte
der Satz f\"ur $(N-1)$--dimensionale Z--Gebiete.
F\"ur eine Karte $(V,h)$
und eine in
$R^{q,h(\Gamma_1\cap V)}(h(S\cap V))
\cap D^{q,h(\Gamma_2\cap V)}(h(S\cap V))$
beschr\"ankte
Familie $\{E^\alpha \}$
mit kompaktem Tr\"ager in $U_N(1/3)\cap \ol{h(S\cap V)}$
gen\"ugt es nach Lemma \ref{wwp12} und $(\ref{suppined})$
zu zeigen,
da{\ss} f\"ur ein $R\in(1/3,1)$ die Familie
$\{E^\alpha\}$
in $L_2^q(U_N(R)\cap \ol{h(S\cap V)})$
relativ kompakt ist.\new
F\"ur innere Karten folgt dies direkt
aus \cite{weck}.
F\"ur eine Randkarte existiert nach Definition \ref{zgebiet}
ein Z--Gebiet
$\cals=(\tilde{S},\gamma_1,\gamma_2)$ mit
$(U_N^-,h(\Gamma_1\cap V),h(\Gamma_2\cap V))=C(\cals)$.
Nach Induktionsannahme hat dieses die Eigenschaft
aus Lemma \ref{wwp15}.
%
Wegen Bemerkung  \ref{bemerk2} und
$(S,\Gamma_1,\Gamma_2)\in\calm_D(M)$
finden wir
Folgen
$(\Phi^{i,\alpha}_k)\subset C_\infty^{q,C(\gamma_i)}(\ol{U_N^-})$
mit
\beqa
 \Phi^{1,\alpha}_k&\ra& E^{\alpha}\mbox{ in } R^q(U_N^-) \non\\
 \Phi^{2,\alpha}_k&\ra& E^{\alpha}\mbox{ in } D^q(U_N^-) \non
\eeqa
f\"ur $k\ra\infty$.
Damit sind die Voraussetzungen von Lemma \ref{lemma9} erf\"ullt.
Dies liefert schlie{\ss}lich die Behauptung. \qed
\abstandk
Als Folgerung notieren wir
\begin{korollar}
  \formel{rqpoincare}  \new
  Sei $(S,\Gamma_1,\Gamma_2)\in\calm_D(M)$.
  Ist die
  Einbettung
  $R^{q,\Gamma_1}(S)\cap D^{q,\Gamma_2}(S)\hookrightarrow L_2^q(S)$
  kompakt (dies ist insbesondere in Z--Gebieten erf\"ullt), so
  gelten:
 \begin{enumerate}
 \item[i)]
 $\gibt{c>0}\alle{\Phi\in(\rq_0(S)\cap\dq(S))\cap(\rq_0(S)\cap\dq_0(S))^{\perp}}
  \normu{\Phi}{q,S} \le c\normu{\div\Phi}{q-1,S}$
  \item[ii)]
 $\gibt{c>0}\alle{\Phi\in(\rq(S)\cap\dq_0(S))\cap(\rq_0(S)\cap\dq_0(S))^{\perp}}
 \normu{\Phi}{q,S} \le c \normu{\rot\Phi}{q+1,S}$
 \item[iii)]  Die R\"aume $\rot R^q(S)$ und $\div D^q(S)$ sind
  abgeschlossen in $L_2^q(S)$.
 \end{enumerate}
\end{korollar}
\beweis Lemma \ref{wwp9}.
\section{Ein Regularit\"atssatz}
\markboth{EIN REGULARIT\"ATSSATZ}{EIN REGULARIT\"ATSSATZ}
 \label{chregularitaet}
\begin{satz}   \formel{regularitaet}  \new
 Seien $S$ glatt, $m\in\nnn\cup\{0\}$,
  und $\eps$ eine zul\"assige $C_{m+1}(\ol{S})$--Transformation
  (d.h. $\eps$ ist zul\"assig und f\"ur alle Karten
  $(V,h)$ um $x$
  sind in der Matrixdarstellung
  $(\ref{matrixdar})$
  die Eintr\"age $\eps_{I,J}(x)$
  aus
  $C_{m+1}({\ol S}\cap V)$).
  Ferner erf\"ulle $H\in \rqn(S)\cap \eps^{-1}D^q(S)$
 \beqa
  \rot H&\in& H_m^{q+1}(S) \non\\
  \div\eps H&\in& H_m^{q-1}(S) \ .\non
 \eeqa
 Dann gelten $H\in H_{m+1}^{q}(S)$ und
 \beqa
 \gibt{c>0}\normu{H}{H_{m+1}^q(S)}\le
 c(\normu{H}{q,S} +\normu{\rot H}{H_m^{q+1}(S)}
   +\normu{\div\eps H}{H_m^{q-1}(S)} )\ .\non
\eeqa
\end{satz}
Im Falle $q=1$, $N=3$ wurde dies in \cite{weber} bewiesen.
Wir gehen einen \"ahnlichen Weg, k\"onnen jedoch nicht
$(\ref{risth})$ f\"ur $(q-1)$-- bzw. $(q+2)$--Formen benutzen.
Wir ben\"otigen einige Vorbereitungen:
\begin{lemma} \formel{spiegel} \new
 Seien $r>0$, $x':=(x_1,\cdots,x_{N-1})$
 und
 \beqa
  \Abb{\tau}{U_N^+(r)}{U_N^-(r)}{x}{(x',-x_N)\ .} \non
  \eeqa
  Der Spiegelungsoperator\label{pxspiegel}
 \beqa
 \Abb{S_{\mbox{{\tiny rot}}}}{R^q(U_N^-(r))}{R^q(U_N(r))}{H}{
 \left\{\begin{array}{ccc}
  H&\mbox{ in }&U_N^-(r) \\
  \tau^*H&\mbox{ in }& U_N^+(r)
  \end{array} \right.
  }                         \non
 \eeqa
 ist wohldefiniert, linear, stetig und hat die Eigenschaften
 \beqa
  \supp H\subset \ol{U_N^-(r')}
  &\Ra& \supp S_{\mbox{{\tiny rot}}}H \subset \ol{U_N(r')}\mfur r'<r\non\\
  \gibt{c>0}\normu{\rot S_{\mbox{{\tiny rot}}}H}{q+1,U_N(r)} &\le&
   c\normu{\rot H}{q+1,U_N^-(r)} \non\\
  \rot S_{\mbox{{\tiny rot}}}H&=&\left\{\begin{array}{ccc}
     \rot H&\mbox{ in }& U_N^-(r)\\
     \tau^*\rot H&\mbox{ in }& U_N^+(r)
  \end{array} \right.
   \ .\formel{srotsrot}
 \eeqa
\end{lemma}
\beweis
Wegen Lemma \ref{lemm0} gen\"ugt es,
$S_{\mbox{{\tiny rot}}} H\in R^q(U_N(r))$
sowie $(\ref{srotsrot})$ f\"ur Formen $H\in C_\infty^q(\ol{U_N^-})$
zu zeigen.
Die Aussagen \"uber den  Tr\"ager und die Stetigkeit
folgen dann direkt.
Seien $U_\pm:=U_N^{\pm}(r)$, $U_0:=U_N^0(r)$,
$U:=U_N(r)$ und $\iota_\pm,\iota_0$ die Einbettungen
$U_\pm,U_0\ra U$.
Wir beachten, da{\ss} $\tau$ die Orientierung ver\"andert
und berechnen
f\"ur $\Phi\in \cnull{\infty}{q+1}(U)$
\beqa
 (-1)^q\spu{S_{\mbox{{\tiny rot}}}H}{\div\Phi}{q,U}
  &=&\int_{U_-} H\wedge\iota_-^*\d *\ol{\Phi}
    +\int_{U_+} \tau^*H\wedge\iota_+^*\d *\ol{\Phi} \non\\
  &=&\int_{U_-} H\wedge\left(\iota_-^*\d *\ol{\Phi}
   -(\tau^{-1})^*\iota_+^*\d *\ol{\Phi}\right) \non\\
  &=&\int_{U_-}H\wedge\d\left(\iota_-^**\ol{\Phi}
    -(\tau^{-1})^*\iota_+^**\ol{\Phi}\right) \ .\non
\eeqa
Aus dem Satz von Stokes folgt
\beqa
 \spu{S_{\mbox{{\tiny rot}}}H}{\div\Phi}{q,U}
   &=& -\int_{U_-}\d H\wedge\left(\iota_-^**\ol{\Phi}
    -(\tau^{-1})^*\iota_+^**\ol{\Phi}\right)\non\\
   && +\int_{U_0}\iota_0^*H\wedge \iota_0^*
    (\iota_-^*-(\tau^{-1})^*\iota_+^*)*\ol{\Phi} \ .\non
\eeqa
Wegen $\iota_-\circ\iota_0-\iota_+\circ\tau^{-1}\circ\iota_0=0$
verschwindet das zweite Integral, und wir erhalten
\beqa
 \spu{S_{\mbox{{\tiny rot}}}H}{\div\Phi}{q,U}
  &=&-\left(\int_{U_-}\d H\wedge \iota_-^**\ol{\Phi}
  +\int_{U_+}\tau^*\d H\wedge \iota_+^**\ol{\Phi}\right)\non\\
  &=& -\spu{F}{\Phi}{q+1,U} \non\\
  &&\mmit F:= \left\{
\begin{array}{cc}
 \rot H&\mbox{in } U_-  \\
 \tau^*\rot H& \mbox{in } U_+ \ .
\end{array}          \right. \non
\eeqa
\qed
\abstand
Auf $D^q(U_N^-(r))$ k\"onnen wir den
Spiegelungsoperator mittels\label{pxspiegeld}
\beqa
 S_{\mbox{{\tiny div}}}H=\kappa_q*S_{\mbox{{\tiny rot}}}*H \non
\eeqa
erkl\"aren. Dieser hat dann die entsprechenden Eigenschaften.\new
\begin{lemma}  \formel{fourier} \new
  Seien $N\ge 3$ und $r>0$.
  Es gibt eine Konstante $c>0$ und zu jedem
   $H\in D^q_0(\rrr^N)$ mit $\supp H\subset U_N(r)$
   ein  $A\in H_1^{q+1}(\rrr^N)$ mit
   \beqa
  \div A=H \ ,
  \ \normu{A}{H_1^{q+1}(\rrr^N)}\le c\normu{H}{q,\rrr^N}\ .\non
 \eeqa
\end{lemma}
\beweis
Seien $x_i$ und
$y_i$, $i=1,\cdots,N$ kartesische Koordinaten.
Wir setzen f\"ur $\Phi=\sqqn \Phi_I\dx^I\in\cnull{\infty}{q} (\rrr^N)$
\label{pxfourier}
\beqa
 F\Phi(x)&:=&\sqqn \calf_0\Phi_I(x)\dx^I \non\\
 F^{-1}\Phi(y)&:=&\sqqn \calf_0^{-1}\Phi_I(y)\dy^{I}\ ,\non
\eeqa
wobei
\beqa
 \calf_0\Phi_I(x):=\frac{1}{(2\pi)^{N/2}}\int_{\rrr_N}
 e^{-\i\sp{y}{x} }\Phi_I(y)\dy \non
\eeqa
die Fouriertransformierte von
$\Phi_I$ ist.
Wir erinnern an die Operatoren $\hat{R}$, $\hat{T}$ und $m$ aus
$(\ref{wwoperator})$. Nach $(\ref{rotlokal})$ und
(\cite[Satz 10.5]{triebel}) gilt
\beqa
 F(\rot \Phi)&=&\i\hat{R}F\Phi \ .\non
\eeqa
Sei nun $H\in D^q_0(\rrr^N)$ mit $\supp H\subset U_N(r)$.
Wegen
\beqa
  |\calf_0 H_I(x)|\le \int_{U_N(r)}|H_I(y)|\dy\le
     c \normu{H}{q,\rrr^N} \non
\eeqa
(vgl. \cite[Satz 10.6]{triebel})
sind die Komponenten von $FH$ beschr\"ankt.
Wir definieren
$\hat{A}:=-\i m^{-2}\hat{R}FH$ (im Nullpunkt sei $\hat{A}=0$).
Aus
\beqa
 |\hat{A}_J(x)|&=&|(|x|^{-2}\sum_{n=1}^Nx_n\dx^n\wedge
  \sum_{I\in\cals(q,N)} \calf_0H_I(x)\dx^{I})_J| \non\\
   &\le& c\sum_{I\in\cals(q,N)} |x|^{-1}|\calf_0H_I(x)|
   \formel{didadu}
\eeqa
folgt, da{\ss} $\hat{A}$ und die Transformierte $F^{-1}\hat{A}$ in
$L_2^q(\rrr^N)$ liegen (in einer Nullumgebung  sch\"atzen
wir den Term $|\calf_0H_I(x)|$,
im Komplement den Term $|x|^{-1}$ durch eine
Konstante ab).
Mit $(\ref{wwrt})$ und \cite[10.25]{triebel}
erhalten wir  f\"ur $\Phi\in\cnull{\infty}{q-1} (\rrr^N)$
\beqa
 \spu{\hat{T}FH}{\Phi}{q-1,\rrr^N}
    &=&\int_{\rrr^N} \hat{T}FH\wedge *\ol{\Phi} \non\\
    &=&\int_{\rrr^N} FH\wedge *\hat{R}\ol{\Phi} \non\\
    &=&\i\spu{H}{\rot F^{-1}\Phi}{q+1,\rrr^N} \non\\
     &=&0\non
\eeqa
und f\"ur $\Phi\in\cnull{\infty}{q} (\rrr^N)$
\beqa
\spu{F^{-1}\hat{A}}{\rot\Phi}{q+1,\rrr^N}
  &=&-\int_{\rrr^N} m^{-2}\hat{R}FH\wedge *\hat{R}\ol{F\Phi} \non\\
 &=&-\int_{\rrr^N} m^{-2}\hat{T}\hat{R}FH\wedge *\ol{F\Phi}\non\\
 &=&-\spu{FH}{F\Phi}{q,\rrr^N} \non\\
 &=&-\spu{H}{\Phi}{q,\rrr^N}\ ,\non
\eeqa
wobei wir im vorletzten Schritt $(\ref{wwmz})$
benutzt haben. Es folgt $\div F^{-1}\hat{A}=H$.
Mit
\beqa
 \normu{{x_k} \hat{A}}{q+1,\rrr^N}  &\le&
  c \normu{H}{q,\rrr^N} \non
\eeqa
(vgl. $(\ref{didadu})$) erhalten wir $A:=F^{-1}\hat{A}\in H_1^{q+1}(\rrr^N)$ und
die Absch\"atzung.
\qed
Zur Vorbereitung des nachfolgenden Lemmas betrachten wir
eine Teilmenge $U$ des $\rrr^N$.
F\"ur eine Form
\beqa
        \Phi=\sum_{I\in\cals(q,N)}\Phi_I\dx^{I}\in L_2^q(U) \non
\eeqa
ist $\Phi=\Phi^\tau+\Phi^\rho$ mit
\beqa
 \Phi^\tau&:=&\sum_{I\in\cals(q,N-1)}\Phi_I\dx^{I}\non\\
 \Phi^\rho&:=&\sum_{I\in\cals(q,N),N\in I}\Phi_I\dx^{I}\non
\eeqa
eine orthogonale Zerlegung in $L_2^q(U)$.
\begin{lemma} \formel{h1schluss} \new
Seien $U\subset\rrr^N$,
$\eps$ eine zul\"assige $C_1(\ol{U})$--Transformation
und $H\in L_2^q(U)$.
Liegen  die Ausdr\"ucke $\partial_iH$ und
$\partial_i(\eps H)$ f\"ur $i=1,\cdots,N-1$
sowie $\partial_N H^\tau$
und
$\partial_N (\eps H)^{\rho}$
in $L_2^q(U)$,
so gilt $H\in H_1^q(U)$.\new
Hierbei wollen wir die Ableitungen komponentenweise verstehen.
\end{lemma}
\beweis
Aus den Vorausetzungen folgt $H^\tau\in H_1^q(U)$ und somit
auch
\beqa
 \partial_N(\eps H^{\rho})^{\rho}=
 \partial_N(\eps H)^{\rho} -\partial_N(\eps H^{\tau})^{\rho}
 \in L_2^q(U)\ .\non
\eeqa
Wir erhalten $(\eps H^{\rho})^{\rho}\in H_1^q(U)$.
Da die auf dem ''Normalenteil'' agierende
Einschr\"ankung $\eps^{\rho,\rho}$ von $\eps$
punktweise invertierbar
mit
$C_1(\ol{U})$ Eintr\"agen ist,
haben wir das Lemma bewiesen.
\qed
\abstand
Nun zum Beweis von Satz \ref{regularitaet}.
Wir wollen wieder lokalisieren und
beschr\"anken uns auf den schwierigeren
Fall der Randkarten.
Wir definieren $U:=U_N^-$ sowie $U(r):=U_N^-(r)$.
F\"ur Zahlen $r\in(0,\infty)$ sind nach \cite{weck}
und Lemma \ref{wwp9}
die
R\"aume $\rot R^q(U(r))$, $\rot\rnull{\ }{q}(U(r))$
usw. abgeschlossen. Dar\"uber hinaus gelten
die Poincar\'{e} Absch\"atzungen
in Lemma \ref{rqpoincare} mit
$(\Gamma_1,\Gamma_2)\in \{(\emptyset,\partial U),
(\partial U,\emptyset)\}$.
Aus der Segmenteigenschaft folgt dann
$R^{q,\partial U}(U)=\hat{R}^{q,\partial U}(U)
 =\rqn(U)$.
\abstandk
F\"ur eine zul\"assige $C_{m+1}(\ol{U})$--Transformation
$\eps$ gen\"ugt es zu zeigen
\beqa
 H&\in& H_{m+1}^q(U)
  \formel{hmqziel} \\
 \normu{H}{H_{m+1}^q(U)} &\le&
 c(\normu{H}{q,U} +\normu{\rot H}{H_m^{q+1}(U) }
 +\normu{\div\eps H}{H_m^{q-1}(U)} )  \non
\eeqa
f\"ur
\beqa
 H&\in&\rnull{\ }{q}(U)\cap\eps^{-1}D^q(U) \non\\
 \rot H&\in& H_m^{q+1}(U) \non\\
  \div \eps H&=& H_m^{q-1} (U) \non\\
  \supp H&\subset\subset& \ol{U(r)}\mbox{ f\"ur ein }
  r\in(0,1) \ .\non
\eeqa
Die Transformation
$\eps$ ist im allgemeinen nicht identisch mit
$\eps$ aus Satz \ref{regularitaet}, sondern das
Produkt $\eps_{h^{-1}} \cdot\eps\circ h^{-1}$ f\"ur die
Karte $(V,h)$.
Auf die Felder $(H_I)_I$ der Komponentenfunktionen
wirkt $\eps$ wie eine symmetrische
gleichm\"a{\ss}ig positiv definite
Matrix mit Eintr\"agen aus $C_{m+1}(\ol{U})$.
\abstandk
Sei zun\"achst $N\ge 3$.
Wir zeigen $(\ref{hmqziel})$
per Induktion \"uber $q$ und $m$. Da der Fall $q=0$ nach $(\ref{risth})$
schon bewiesen
ist, nehmen wir an, die Aussage gelte f\"ur $q-1$.
Sei also $m=0$.
Wir wollen zeigen
\beqa
 \partial_i(\eps H) &\in& L_2^q(U)
 \formel{tanabl} \\
 \normu{\partial_i\eps H}{q,U}
 &\le& c(\normu{H}{q,U} +\normu{\div\eps H}{q-1,U}
  +\normu{\rot H}{q+1,U} )
 \formel{tanabsch}
\eeqa
f\"ur $i=1,\cdots,N-1$,
wobei wir die Ableitung komponentenweise verstehen.
Aus Symmetriegr\"unden gen\"ugt es, den Fall
$i=1$ zu betrachten.
Wir w\"ahlen $\Delta\in(0,1)$ mit $r+4\Delta<1$ und setzen
$W_j:=U(r+j\Delta)$ sowie
\beqa
 \Abb{\tau_h}{\rrr^N_-}{\rrr^N_-}{x}{(x_1+h,x_2,\cdots,x_N)} \ ,
  0<|h|<\Delta \ .\non
\eeqa
Wegen
$\dy^{i}=\dx^{i}$ f\"ur
$y^{i}:=\tau_{h,i}(x)$ k\"onnen wir die Koordinaten
im Urbild und Zielbereich identifizieren. Somit ist
der Ausdruck
\beqa
 \delta_h H:=\frac{1}{h}(\tau_h^*-\id^*)H \non
\eeqa
wohldefiniert, wirkt auf die Komponenten wie
der Differenzenquotient und tauscht mit Rotation,
Sternoperator und Divergenz. Weiter gilt
f\"ur alle $F,G\in L_2^q(U)$ mit Tr\"ager in $\ol{W_3}$
\beqa
   |\spu{\delta_hF}{G}{q,U}|
  &=&|\spu{F}{\delta_{-h}G}{q,U}| \non\\
   \delta_h\eps F&=&\eps\delta_h F+(\delta_h\eps)\tau_h^*F \non\\
   \normu{(\delta_h\eps) \tau_h^*F}{q,U}
   &\le&c \normu{F}{q,U} \ .\formel{deltaepsab}
\eeqa
Hierbei sei
\beqa
 (\delta_h\eps)\Phi(x) := \sum_{I\in\cals(q,N)}
 \sum_{J\in\cals(q,N)} (\delta_h\tilde{\eps}_{I,J}(x))\Phi_J(x)\d x^{I} \non
\eeqa
f\"ur $\Phi(x)=\sum_{I\in\cals(q,N)}\Phi_I(x)\dx^{I}$ und
die oben erw\"ahnten Matrixdarstellung $\tilde{\eps}$ von
$\eps$.
Wie in \cite[Theorem 3.13]{agmon} zeigt man
f\"ur $m\in\nnn$
\beqa
 \alle{F\in H_m^q(U),\mbox{ {\footnotesize supp}} F\subset \ol{W_3}}
  \normu{\delta_h F}{H_{m-1}^q(U)}
&\le& \normu{F}{H_{m}^q(U)} \ .\formel{deltahnorm}
\eeqa
Die Form $H$ erf\"ullt
\beqa
 \supp \delta_h H\subset\subset \ol{W_{1}} \ ,\non
\eeqa
und nach \cite[Theorem 3.15]{agmon} und $(\ref{deltaepsab})$
gen\"ugt es,
die Absch\"atzung
\beqa
 |\spu{\delta_h\eps H}{\Phi}{q,W_{1}} |
 &\le& c(\normu{H}{q,U} +\normu{\div\eps H}{q-1,U}
  +\normu{\rot H}{q+1,U} )\normu{\Phi}{q,W_1} \non
\eeqa
f\"ur alle $\Phi\in \cnull{\infty}{q} (W_1)$ zu
zeigen. In $U$ zerlegen wir gem\"a{\ss} Lemma \ref{orthzerl}
\beqa
 \Phi&=&\hat{\Phi}_1+\eps^{-1}\hat{\Phi}_2 \non\\
   \hat{\Phi}_1&\in& \rnull{0}{q} (U)\non\\
  \hat{\Phi}_2&=&\div\Phi_2\ ,\
  \Phi_2 \in D^{q+1} (U)\cap \rot\rnull{\ }{q} (U)
  \formel{weberzerl1}\\
   &&\normu{\Phi_2}{D^{q+1}(U)} \le c\normu{\Phi}{q,W_1}
  \ .\non
\eeqa
Da alle Betti--Zahlen von $U$ verschwinden,
folgt aus \cite[Satz 1, Satz 2]{picard2} und Lemma \ref{orthzerl}
\beqa
 \rnull{0}{q} (U) = \rot(\rnull{\ }{q-1} (U)
  \cap \div D^{q}(U)) \ .\non
\eeqa
Nach Induktionsvoraussetzung
erhalten wir (vgl. Lemma \ref{randver})
\beqa
 \hat{\Phi}_1&=&\rot\Phi_1\ ,
 \ \Phi_1\in \rnull{\ }{q-1} (U)
  \cap \div D^q(U)  \non\\
   \mmit \chi\Phi_1&\in& H_1^{q-1}(U)
    \non\\
     \normu{\chi\Phi_1}{H_1^{q-1}(U)} &\le&c
     \normu{\Phi}{q,W_1} \formel{phieab}
\eeqa
f\"ur
$\chi\in\cnull{\infty}{\ } (U_N(r+2\Delta))$ mit
$\chi=1$ in $W_1$.\new
Die Form $\chi\Phi_2$ hat kompakten Tr\"ager
in $U\cup U_N^0$. Die Nullfortsetzung von $S_{\mbox{{\tiny div}}}\chi\Phi_2$
auf $\rrr^N$ liegt in $D^{q+1}(\rrr^N)$, und
$\tilde{\Phi}_2:=\div S_{\mbox{{\tiny div}}}\chi\Phi_2$
erf\"ullt
\beqa
 \tilde{\Phi}_2&=&\hat{\Phi}_2 \mbox{ in } W_1\non\\
 \div\tilde{\Phi}_2 &=&0\ ,
 \ \supp\tilde{\Phi}_2\subset\subset\rrr^N \ .\non
\eeqa
Nach Lemma $\ref{fourier}$ existiert $A\in H_1^{q+1}(\rrr^N)$
mit
\beqa
 \div A&=&\tilde{\Phi}_2 \non\\
 \normu{A}{H_1^{q+1}(\rrr^N)}
 &\le& c\normu{\tilde{\Phi}_2}{q,U_N} \non\\
 &=&c \normu{\div S_{\mbox{{\tiny div}}}\chi\Phi_2}{q,U_N} \non\\
  &\le& c\normu{\Phi_2}{D^{q+1}(U)} \non\\
  &\le&c\normu{\Phi}{q,W_1} \ .\formel{aabsch}
\eeqa
Es folgt
\beqa
 |\spu{\delta_h\eps H}{\Phi}{q,W_1}|
 &=& |\spu{\delta_h\eps H}{\rot\Phi_1+\eps^{-1}
  \div A}{q,U}| \non\\
 &=& |\spu{\delta_h\eps H}{\rot\chi\Phi_1+\eps^{-1}
 \div\chi A}{q,U} |\non\\
 &=& |\spu{\eps H}{\delta_{-h}(\rot\chi\Phi_1+\eps^{-1}
 \div\chi A)}{q,U} |\non\\
 &\le&|\spu{\eps H}{\rot\delta_{-h}\chi\Phi_1}{q,U} |
     +|\spu{H}{\div\delta_{-h}\chi A}{q,U} | \non\\
     &&+|\spu{\eps H}{
      (\delta_{-h}\eps^{-1})\tau_{-h}^*\div\chi A}{q,U} | \non\\
  &=:& I_1+I_2+I_3
      \ .\non
\eeqa
Den Term $I_3$ sch\"atzen wir mit $(\ref{deltaepsab})$
und $(\ref{aabsch})$ durch
$c\normu{H}{q,W_1} \normu{\Phi}{q,W_1} $ ab.
Mit $\Phi_1$ liegt auch
$\delta_{-h}\chi\Phi_1$ in $\rnull{\ }{q-1} (U)$.
Aus $(\ref{deltahnorm})$ und $(\ref{phieab})$ folgern wir
\beqa
 I_1\le c\normu{\div\eps H}{q-1,W_1}
  \normu{\Phi}{q,W_1} \ ,\formel{i1kleiner}
\eeqa
mit $\delta_{-h}\chi A\in D^{q+1}(U)$ und
$H\in \rnull{\ }{q} (U)$ aus $(\ref{deltahnorm})$ und
$(\ref{aabsch})$
\beqa
 I_2\le c\normu{\rot H}{q+1,W_1}
  \normu{\Phi}{q,W_1} \ .\non
\eeqa
Damit sind $(\ref{tanabl})$ und $(\ref{tanabsch})$
f\"ur $i=1,\cdots,N-1$
gezeigt.
Die Normalenableitungen erhalten wir wie folgt:
Wir wissen f\"ur $i=1,\cdots,N-1$
\beqa
 \partial_{i}\eps H&=&\eps\partial_iH +(\partial_i\eps)H
 \in L_2^q(U)
   \formel{regform1}\\
  \Ra \partial_i H&=&\eps^{-1} (\partial_i\eps H-(\partial_i\eps)H)
  \in L_2^q(U) \ ,i=1,\cdots,N-1 \ .\formel{regform2}
\eeqa
Nach $(\ref{rotlokal})$, $(\ref{divlokal})$ gelten dann
\beqa
 \partial_NH_I&=&\sigma(N,I)((\rot H)_{I+N}-\sum_{j\in I}
   \sigma(j,I+N-j)\partial_j H_{I+N-j} )\non\\
   && \in L_2(U) \mfur N\not \in I \formel{regform3}\\
 (\partial_N\eps H)_I&=&\sigma(N,I-N)((\div\eps H)_{I-N}
   -\sum_{j\not\in I}
  \sigma(j,I)(\partial_j\eps H)_{I-N+j})\non\\
  && \in L_2(U) \mfur N\in I \ .\formel{regform4}
\eeqa
Aus Lemma \ref{h1schluss} folgt
$H\in H_1^q(U)$.
Damit haben wir den Fall $m=0$ bewiesen. Gelte der Satz
f\"ur $m-1$. Aus den
Voraussetzungen im Falle $m$ folgt
\beqa
 H,\eps H&\in& H_m^q(U) \formel{indvor}\\
 \normu{H}{H_{m}^q(U)}&\le&
 c(\normu{H}{q,U} +\normu{\rot H}{H_{m}^{q+1}(U)}
   +\normu{\div\eps H}{H_{m}^{q-1}(U)} )\non\\
   &=:&C_{m}(H)\ .\formel{indvorzwei}
\eeqa
F\"ur einen Differentialoperator
$D_m:=\partial_1^{m_1}\cdots\partial_{N-1}^{m_{N-1}}$
mit $m=\sum_{j=1}^{N-1}m_j$
wollen wir
wie oben zeigen
\beqa
 |\spu{\delta_h \eps D_m H}{\Phi}{q,W_1} |
 \le c(\normu{H}{q,U} +\normu{\div\eps H}{H_m^{q-1}(U)}
  +\normu{\rot H}{H_m^{q+1}(U)} )\normu{\Phi}{q,W_1} \non
\eeqa
f\"ur alle $\Phi\in \cnull{\infty}{q} (W_1)$.
Mit der gleichen Zerlegung wie zuvor
erhalten wir
\beqa
  |\spu{\delta_h\eps D_mH}{\Phi}{q,W_1}|
 &\le&|\spu{\eps D_mH}{\rot\delta_{-h}\chi\Phi_1}{q,U} | \non\\
     &&+|\spu{D_mH}{\div\delta_{-h}\chi A}{q,U} | \non\\
     &&+|\spu{\eps D_m H}{
      (\delta_{-h}\eps^{-1})\tau_{-h}^*\div\chi A}{q,U} | \non\\
  &=:& I_1+I_2+I_3
      \ .\non
\eeqa
Wegen
\beqa
 D_m(\eps H)=\eps D_m H +r(\eps,H) \formel{restlein}
\eeqa
mit einem Ausdruck $r(\eps,H)\in H_1^q(U)$, dessen Tr\"ager
in $\ol{W_0}$ liegt, und der nur
Ableitungen von $\eps$ bis zur Ordnung
$m$ und von $H$ bis zur Ordnung $m-1$ enth\"alt, erhalten wir
\beqa
  I_1&\le&|\spu{D_m(\eps H)}{\rot\delta_{-h}\chi\Phi_1}{q,U} |
      +|\spu{\delta_hr(\eps,H)}{\rot\chi\Phi_1}{q,U} |\non\\
      &=:&I_{1,1}+I_{1,2} \ . \non
\eeqa
Den zweiten Term sch\"atzen wir nach
$(\ref{deltahnorm})$, $(\ref{phieab})$ und $(\ref{indvorzwei})$ ab durch
\beqa
 I_{1,2}&\le&c\normu{r(\eps,H)}{H_1^q(U)} \normu{\Phi}{q,W_1}
 \non\\
     &\le& c\normu{H}{H_m^q(U)} \normu{\Phi}{q,W_1} \non\\
     &\le& c\, C_{m}(H) \normu{\Phi}{q,W_1} \ .\non
\eeqa
Aus $\div\eps H\in H_m^q(U)$
folgt f\"ur $\Psi\in\cnull{\infty}{q-1} (U)$
\beqa
 \spu{D_m(\eps H)}{\rot\Psi}{q,U}
  &=& (-1)^{m}\spu{\eps H}{D_m\rot\Psi}{q,U} \non\\
   &=& (-1)^{m}\spu{\eps H}{\rot D_m\Psi}{q,U} \non\\
   &=& (-1)^{m+1}\spu{\div\eps H}{D_m\Psi}{q,U} \non\\
    &=& -\spu{D_m\div\eps H}{\Psi}{q,U} \ ,\non
\eeqa
also $\div D_m\eps H=D_m\div\eps H\in L_2^q(U)$.
Wir erhalten wie in $(\ref{i1kleiner})$
\beqa
 I_{1,1}&\le&|\spu{
  D_m\div(\eps H)}{\delta_{-h}\chi\Phi_1}{q,U} |
 \le c\normu{\div\eps H}{H_m^{q-1}(U)} \normu{\Phi}{q,W_1} \ .\non
\eeqa
Genauso verfahren wir mit $I_2$:
\beqa
 I_{2}\le c\normu{\rot H}{H_m^{q+1}(U)}\normu{\Phi}{q,W_1}
 \ .\non
\eeqa
Den letzten Term $I_3$ k\"onnen wir
wie im Fall $m=0$ behandeln und erhalten
nach entsprechenden \"Uberlegungen f\"ur
die anderen tangentialen Ableitungen mit $(\ref{restlein})$
\beqa
   \partial_iD_m\eps H&\in& L_2^q(U) \non\\
  \normu{\partial_iD_m\eps H}{q,U} &\le& c\,C_m(H) \non
\eeqa
f\"ur $i=1,\cdots,N-1$.
\abstandk
Um die Normalenableitungen zu untersuchen, setzen
wir f\"ur $1\le k\le m$
\beqa
 \tilde{H}^k:=\partial_N^k\partial^{\alpha_{m-k}} H\ ,
 \ |\alpha_j|=j\ ,\non
\eeqa
wobei $\alpha_j$ nur Tangentialkomponenten enth\"alt.
Unter Verwendung von $(\ref{indvor})$
und den Formeln $(\ref{regform1})$ bis
$(\ref{regform4})$
k\"onnen wir induktiv f\"ur $k=1,\cdots, m$ die
Voraussetzungen
von Lemma \ref{h1schluss} (mit $\tilde{H}^k$ statt $H$)
und damit
$\tilde H_k\in H_1^q(U)$
zeigen.
\abstandk
Damit ist der Fall $N\ge 3$ bewiesen.
In den F\"allen $N=1$, $q$ beliebig und
$N=2$, $q=0$ oder $q=2$ ist nichts zu zeigen.
Im Fall $N=2$, $q=1$ liegt die
Form ${\Phi}_2$ aus $(\ref{weberzerl1})$
wegen $(\ref{risth})$
in $H_1^2(U)$, so da{\ss} wir den Beweis
ohne Lemma \ref{fourier}
durchf\"uhren k\"onnen.
\qed
\abstand
Im Falle wechselnder Randbedingungen w\"urde der Beweis
scheitern, da die Anwendung von $\delta_h$ die Randbedingungen
nicht respektiert.
\section{Spurs\"atze}
\markboth{SPURS\"ATZE}{DIE SPURS\"ATZE IN $R^q$ UND $D^q$}
\label{chspursatz}
In diesem Kapitel betrachten wir Spuren und Teilspuren der
Formen aus $R^q$ und $D^q$, wobei wir die
folgenden S\"atze nur f\"ur den ersten Raum beweisen werden.
Die Behauptungen f\"ur den Raum $D^q$ sind
die dualen Resultate, die man mit Hilfe des Sternoperators
unter Verwendung von $(\ref{sternr})$, $(\ref{sterndv})$ und den
in diesem Zusammenhang bewiesenen Aussagen gewinnen kann.
\subsection{Die Spurs\"atze in $R^q$ und $D^q$}
\formel{chrqspur}\new
Sei in diesem Abschnitt $S$ glatt.
Aus Satz \ref{regularitaet} folgt, da{\ss}
Formen aus $\rnull{\ }{q} (S)\cap D^q(S)$ in $H_1^q(S)$
liegen, und wir k\"onnen auf solche die
Spuroperatoren $T$ aus Lemma \ref{hspur} und $N$
aus $(\ref{definspur})$ anwenden.
\abstandk
Nach $(\ref{phipsi3})$ und dem Satz von Stokes
gilt f\"ur
$\Phi\in\cq(\ol{S})$, $\Psi\in\cqp(\ol{S})$ und $0\le q<N$
\beqa
 \spu{\rot\Phi}{\Psi}{q+1,S} +\spu{\Phi}{\div \Psi}{q,S}
   &=&\int_S\d(\Phi\wedge *\ol{\Psi}) \non\\
  &=&\int_{\partial S} T(\Phi\wedge *\ol{\Psi})\non\\
  &=&\int_{\partial S} T\Phi\wedge \kappa'_q**T*\ol{\Psi}\non\\
  &=& \spu{T\Phi}{N\Psi}{q,\partial S}\ .\non
\eeqa
Mit $(\ref{hdicht})$ erhalten wir
\beqa
  \alle{\Phi\in H_1^q(S)} \alle{\Psi\in H_1^{q+1}(S)}
  \spu{\rot\Phi}{\Psi}{q+1,S} &+&\spu{\Phi}{\div \Psi}{q,S} \non\\
 &&=\spu{T\Phi}{N\Psi}{q,\partial S}\ .\formel{TN}
\eeqa
Dies motiviert, die Tangentialspur
$\gt E\in \hqehm(\partial S)$ einer Form $E\in R^q(S)$
wie folgt zu definieren\label{pxtanrotspur}
\beqa
  \alle{\vp\in\hqeh(\partial S)}
  \spu{\gt E}{\vp}{\hqehm(\partial S)}:=\spu{\rot E}{\chn\vp}{q+1,S}
  +\spu{E}{\div\chn\vp}{q,S} \ .\formel{tanspur}
\eeqa
Angewandt auf Formen
$E\in H_1^q(S)$ erf\"ullt diese
\beqa
   \alle{\vp\in\hqeh(\partial S)}
  \spu{\gt E}{\vp}{\hqehm(\partial S)} =
  \spu{T E}{\vp}{q,\partial S} \ .\non
\eeqa
Wir zeigen:
\begin{satz} \formel{rotspursatz} \new
F\"ur jedes $E\in R^q(S)$ liegt die Tangentialspur
$\gt E$ in ${R}^q_{-1/2}(\partial S)$ und es gelten:
\begin{enumerate}
  \item[i)]
   $\alle{E\in R^q(S)} \alle{\Psi\in H_1^{q+1}(S)}
   \spu{\gt E}{N\Psi}{H_{-1/2}^q(\partial S)} =
    \spu{\rot E}{\Psi}{q+1,S}
   +\spu{E}{\div \Psi}{q,S}
    $
     \item[ii)] $\alle{E\in R^q(S)} \rot\gt E =\gt \rot E$
  \item[iii)] Die Abbildung
  $\gt:R^q(S)\ra {R}^q_{-1/2}(\partial S)$ ist stetig.
    \end{enumerate}
\end{satz}
\beweis
i) Wegen $(\ref{TN})$ erf\"ullen $\Phi\in C_\infty^q(\ol{S})$
und $\Psi\in H_1^{q+1}(S)$
\beqa
 \spu{\rot \Phi}{\chn N\Psi}{q+1,S} &+&
 \spu{\Phi}{\div \chn N \Psi}{q,S}    \non\\
  &=&\spu{T\Phi}{N\Psi}{q,\partial S} \non\\
  &=&\spu{\rot \Phi}{\Psi}{q+1,S} +
 \spu{\Phi}{\div \Psi}{q,S} \ .\non
\eeqa
Da $C_\infty^q(\ol{S})$ dicht in $R^q(S)$ liegt, folgt
die Behauptung aus $(\ref{tanspur})$ .\new
\abstandk
ii),iii) F\"ur $E\in R^q(S)$ erhalten wir mit
$(\ref{tanspur})$, der Schwarzschen Ungleichung und der
Stetigkeit von $\chn$
\beqa
 |\spu{\gt E}{\vp}{\hqehm(\partial S)}|&\le& c \normu{E}{R^q(S)}
  \normu{\vp}{\hqeh(\partial S)}   \non\\
 \normu{\gt E}{\hqehm(\partial S)}&\le& c\normu{E}{R^q(S)} \ .
 \formel{gtenorm}
 \eeqa
Daraus folgt
f\"ur $\vp\in\hqpdh(\partial S)$ und eine Folge
$(E_k)\subset C_\infty^q(\ol{S})$ mit $E_k\ra E$ in $R^q(S)$
\beqa
   \spu{\gt E}{\div\vp}{\hqehm(\partial S)}&\leftarrow&
   \spu{\gt E_k}{\div\vp}{\hqehm(\partial S)} \non\\
   &=&\spu{T E_k}{\div\vp}{q,\partial S} \non\\
   &=&-\spu{\rot T E_k}{\vp}{q+1,\partial S} \non\\
   &=& -\spu{T\rot E_k}{\vp}{q+1,\partial S} \non\\
   &=& -\spu{\gt\rot E_k}{\vp}{\hqpehm(\partial S)} \non\\
   &\ra& -\spu{\gt\rot E}{\vp}{\hqpehm(\partial S)} \ ,\non
\eeqa
da $\rot :R^q(S)\ra R^{q+1}(S)$ stetig ist.
Nach (\ref{trot}) und $(\ref{gtenorm})$ gelten
\beqa
 \rot\gt E &=&\gt \rot E \in\hqpehm(\partial S)
 \non\\
 \normu{\rot\gt E}{\hqpehm(\partial S)}
  &=&\normu{\gt\rot E}{\hqpehm(\partial S)}
  \le c\normu{E}{R^q(S)} \ ,\formel{gterotnorm}
\eeqa
und wir erhalten Wohldefiniertheit und Stetigkeit der Abbildung
$\gt$.
\qed
\abstandk
Wir untersuchen noch die ''nat\"urlichen'' Eigenschaften des Spuroperators:
F\"ur $E\in \rqn(S)$ gilt
\beqa
 \alle{\psi\in H_{1/2}^q(\partial S)}
 \spu{\gt E}{\psi}{H_{-1/2}^q(\partial S)}
   =\spu{E}{\div \check{N}\psi}{q,S}
   +\spu{\rot E}{\check{N}\psi}{q+1,S} =0 \ .\non
\eeqa
Andererseits zieht $\gt E=0$ nach sich
\beqa
 \alle{\Phi\in C_\infty^{q+1}(\ol{S})}
  \spu{E}{\div \Phi}{q,S}
  +\spu{\rot E}{\Phi}{q+1,S} =
  \spu{\gt E}{N\Phi}{H_{-1/2}^q(\partial S)} =0 \ ,
  \non
\eeqa
also
\beqa
  \gt E=0\Leftrightarrow E\in \rnull{\ }{q }(S) \ .\formel{spurnull}
\eeqa
\abstandk
Mit Hilfe des Sternoperators k\"onnen wir
auf $D^q(S)$ einen Normalenspuroperator erkl\"aren:
\label{pxnordivspur}
\beqa
 \gn E:=\sigma_{q}*\gt *E \ .\non
\eeqa
Die Resultate aus Satz \ref{rotspursatz} werden  durch den
Sternoperator wie folgt \"ubertragen:
\begin{satz} \formel{divspursatz} \new
F\"ur jedes $E\in D^q(S)$ liegt die Normalenspur
$\gn E$ in ${D}^{q-1}_{-1/2}(\partial S)$ und es gelten:
\begin{enumerate}
  \item[i)]
   $\alle{E\in D^q(S)} \alle{\Phi\in H_1^{q-1}(S)}
   \spu{\div E}{\Phi}{q-1,S} +\spu{E}{\rot \Phi}{q,S}
    =\spu{\gn E}{T\Phi}{H_{-1/2}^{q-1}(\partial S)}$
  \item[ii)] $\alle{E\in D^q(S)} \div\gn E=-\gn\div E$
  \item[iii)] Die Abbildung
  $\gn:D^q(S)\ra {D}^{q-1}_{-1/2}(\partial S)$ ist stetig.
     \end{enumerate}
\end{satz}
\abstandk
Wir wollen einen Fortsetzungsoperator
${R}^q_{-1/2}(\partial S)\ra R^q(S)$
konstruieren und zeigen zun\"achst:
\begin{lemma} \formel{ninhnull2} \new
  F\"ur $\vp\in H_{1/2}^{q} (\partial S)$ gelten
  $\chn\vp\in \rnull{\ }{q+1} (S)$
  sowie $\cht\vp\in  \dnull{\ }{q} (S) $.
\end{lemma}
\beweis
F\"ur $\Psi\in C_\infty^{q+2}(\ol{S})$ erhalten wir
mit $(\ref{tnisnull})$ und $(\ref{TN})$
\beqa
 \spu{\chn \vp}{\div \Psi}{q+1,S}
  +\spu{\rot\chn \vp}{\Psi}{q+2,S}
  =\spu{T\chn \vp}{N\Psi}{q+1,\partial S}=0 \ .\non
\eeqa
und damit die erste Behauptung.
Die zweite Behauptung ist das duale Resultat.
\qed
\begin{satz} \formel{rotspursatz2} \new
\begin{enumerate}
  \item[i)]
   Es existiert ein linearer
  stetiger Fortsetzungsoperator \label{pxinvrotspur}
  \beqa
  \chgt :{R}^q_{-1/2}(\partial S)\ra R^q(S) \non
  \eeqa
   mit
  $\gt\chgt=\id$.
  \item[ii)]
   Es existiert ein linearer
  stetiger Fortsetzungsoperator\label{pxinvdivspur}
  \beqa
   \chgn :{D}^{q-1}_{-1/2}(\partial S)\ra D^q(S)\non
   \eeqa
    mit
  $\gn\chgn=\id$.
\end{enumerate}
\end{satz}
\beweis i):
Seien $\lambda\in {R}^q_{-1/2}(\partial S)$ und
(vgl.Korollar \ref{rqpoincare})\label{pageyq}
\beqa
 Y^q&:=&\rnull{ }{q} (S)\cap D^q(S) \non\\
 Y^q_1&:=&Y^q\cap\rot\rnull{ }{q-1} (S)
 =\rot\rnull{ }{q-1} (S)\cap D^q(S) \non\\
 Y^q_2&:=&Y^q\cap\div D^{q+1}(S)
  =\rnull{ }{q} (S)\cap \div D^{q+1}(S) \non\\
 Y^q_3&:=&Y^q\cap\rqnn(S)\cap D_0^q(S)
 =\rqnn(S)\cap D_0^q(S) \ , \non
\eeqa
versehen mit $\spu{\cdot}{\cdot}{R^q(S)\cap D^q(S)} $.
Nach Lemma \ref{orthzerl} mit $\Gamma_1=\partial S$,
$\Gamma_2=\emptyset$
und Korollar \ref{rqpoincare} gilt
\beqa
 Y^q=Y^q_1\oplus Y_2^q\oplus Y_3^q \non
\eeqa
orthogonal in $L_2^q(S)$,
und nach Satz \ref{regularitaet} liegen alle R\"aume in $H_1^q(S)$.
Wir betrachten das Problem (P1): Gesucht $W\in Y_1^{q+2}$ mit
\beqa
 \spu{\div W}{\div\Phi}{q+1,S} = \spu{\rot\lambda}{N\Phi}{\hqpehm(\partial S)}
 \mfur\malle \Phi\in Y_1^{q+2} \ .\formel{p1}
\eeqa
Nach Korollar
\ref{rqpoincare} ist die stetige Bilinearform auf der linken
Seite streng koerzitiv in $Y_1^{q+2}$. Die rechte Seite ist
ein
antilineares stetiges Funktional (Satz \ref{regularitaet}). Nach dem Satz von Lax--Milgram
k\"onnen wir (P1) l\"osen, und die L\"osung $W$ erf\"ullt
\beqa
    \normu{W}{D^{q+2}(S)}\le c\normu{\rot\lambda}{\hqpehm(\partial S)}
    \ . \formel{normw}
\eeqa
Ebenso l\"osen wir (P2): Gesucht $Q\in Y_1^{q+1}$ mit
\beqa
 \spu{\div Q}{\div\Phi}{q,S}
  &=&\spu{\div W}{\Phi}{q+1,S} +\spu{\eta(\lambda)}{\Phi}{q+1,S}
  \non\\
  &&-\spu{\lambda}{N\Phi}{\hqehm(\partial S)}
  \mfur\malle \Phi\in Y_1^{q+1} \formel{p2}\\
  \normu{Q}{D^{q+1}(S)}
  &\le& c (\normu{\div W}{q+1,S}
  +\normu{\lambda}{\hqehm(\partial S)} )
  \ ,\formel{normq}
\eeqa
wobei
\beqa
  \eta(\lambda):=\eta^{q+1}(\lambda):=\sum_{j=1}^{J^{q+1}}y_j^{q+1}
  \spu{\lambda}{Ny_j^{q+1}}{\hqehm(\partial S)} \non
\eeqa
f\"ur eine $L_2^{q+1}(S)$--Orthonormalbasis\label{pxyjq}
$\{y_j^{q+1}\mit j=1,\cdots,J^{q+1}\}$
von $Y_3^{q+1}$.
Testen der rechten Seite von $(\ref{p2})$ mit $y_j^{q+1}$ liefert
\beqa
 \spu{\div W}{y_j^{q+1}}{q+1,S} &+&\spu{\eta(\lambda)}{y_j^{q+1}}{q+1,S}
 -\spu{\lambda}{Ny_j^{q+1}}{\hqehm(\partial S)} \non\\
  &=&\spu{\lambda}{Ny_j^{q+1}}{\hqehm(\partial S)}
   -\spu{\lambda}{Ny_j^{q+1}}{\hqehm(\partial S)} \non\\
   &=&0\non
\eeqa
und damit $(\ref{p2})$ f\"ur alle
$\Phi\in Y_1^{q+1}\oplus Y_3^{q+1}$.
Um zu zeigen, da{\ss} die Formel f\"ur alle
Formen $\Phi\in \chn H_{1/2}^{q}(\partial S)
\subset Y^{q+1}$
(nach Lemma \ref{ninhnull2}) richtig ist,
verbleibt es diese mit
\beqa
        \Phi\in Y_2^{q+1}\ ,
        \ \Phi=\div\Psi \mbox{ mit o.B.d.A. }
        \Psi\in Y_1^{q+2} \non
\eeqa
zu testen.
Aus $\rot\Psi=0$, $\div\Psi=\Phi\in H_1^{q+1}(S)$ und
der Randbedingung folgt $\Psi\in H_2^{q+2(S)}$, und
die Formel (\ref{p1}) impliziert
\beqa
 \spu{\div W}{\Phi}{q+1,S} &+&\spu{\eta(\lambda)}{\Phi}{q+1,S}
  -\spu{\lambda}{N\Phi}{\hqehm(\partial S)} \non\\
  &=&\spu{\div W}{\div\Psi}{q+1,S} -\spu{\lambda}{N\div\Psi}{\hqehm(\partial S)} \non\\
 &=&\spu{\rot\lambda}{N\Psi}{\hqpehm(\partial S)}
 +\spu{\lambda}{\div N \Psi}{\hqehm(\partial S)} \non\\
 &=&0\ ,\formel{ndivrot}
\eeqa
also auch $(\ref{p2})$ f\"ur alle
$\Phi\in \chn H_{1/2}^{q}(\partial S)$.
Wir setzen $E:=-\div Q$ und erhalten
$\rot E=\div W+\eta(\lambda)$ und
$\gt E=\lambda$. Die Stetigkeit des
Fortsetzungsoperators
folgt aus $(\ref{normw})$ und $(\ref{normq})$. \abstandk
ii) Mit
\beqa
 \Abb{\chgn}{{D}_{-1/2}^{q-1}(\partial S)}{D^q(S)}{\lambda}{
   \kappa_{q}*\chgt*\lambda}  \non
\eeqa
gilt
\beqa
  \gn\chgn\lambda
    =\sigma_{q}*\gt\chgt*\lambda
   =\lambda \ .\non
\eeqa
Die Stetigkeit folgt aus i).
\qed
\subsection{Wechselnde Randbedingungen}
\markright{WECHSELNDE RANDBEDINGUNGEN}
\label{chrqgespur}
Sei
in diesem Abschnitt
${\cal S}:=(S,\Gamma_1,\Gamma_2)\in\calm_D(M)$ glatt.
Wir wollen die Tangentialspur
auf einem Randst\"uck
$\Gamma_1$ wie im homogenen Fall durch Volumenintegrale
beschreiben.
Dazu definieren wir \label{pxhsqg2}
\beqa
 H_s^{q,\Gamma_2}(\partial S)=\{\vp\in H_s^q(\partial S)
    \mit \vp=0 \mbox{ fast \"uberall in } \Gamma_2 \}
    \ ,\ s\in(0,\infty)\ .\non
\eeqa
Versehen mit der Norm $\normu{\cdot}{ H_s^q(\partial S)}$ ist
dieser ein abgeschlossener Unterraum von $H_s^{q}(\partial S)$.
Den Dualraum bezeichnen wir mit\label{pxhmsqg2}
$H_{-s}^{q,\Gamma_2}(\partial S)$ und erkl\"aren
die Rotation f\"ur $\lambda\in H_{-1/2}^{q,\Gamma_2}(\partial S)$
durch
\beqa
 \alle{\vp\in H_{3/2}^{q+1,\Gamma_2}(\partial S)}
 \spu{\rot\lambda}{\vp}{H_{-3/2}^{q+1,\Gamma_2}(\partial S)}
 :=-\spu{\lambda}{\div \vp}{H_{-1/2}^{q,\Gamma_2}(\partial S)} \ .\non
\eeqa
Ferner definieren wir\label{rmehq}
\beqa
 R_{-1/2}^{q,\Gamma_2} (\partial S)
 &:=&\{\lambda\in H_{-1/2}^{q,\Gamma_2} (\partial S) \mit
  \rot\lambda\in H_{-1/2}^{q+1,\Gamma_2} (\partial S)\} \non\\
 \normu{\lambda}{R_{-1/2}^{q,\Gamma_2} (\partial S)} &:=&
 \normu{\rot\lambda}{H_{-1/2}^{q+1,\Gamma_2} (\partial S)}
 +\normu{\lambda}{H_{-1/2}^{q,\Gamma_2} (\partial S)} \non
\eeqa
und zeigen zun\"achst:
 \begin{lemma}  \formel{ninhnull} \new
F\"ur $\vp\in H_{1/2}^{q,\Gamma_2} (\partial S)$ gelten
$\chn \vp\in D^{q+1,\Gamma_2}(S)\cap\rnull{\ }{q+1} (S)$
und
$\cht \vp\in R^{q,\Gamma_2} (S) \cap \dnull{\ }{q} (S)$.
\end{lemma}
\beweis
Sei
$\vp\in H_{1/2}^{q,\Gamma_2}(\partial S)$.
Aus $(\ref{TN})$ folgt f\"ur
$\Psi\in C_\infty^{q,\Gamma_1}(\ol{S})$
\beqa
 \spu{\chn \vp}{\rot \Psi}{q+1,S}
  +\spu{\div\chn \vp}{\Psi}{q,S}
  &=&\spu{N\chn \vp}{T\Psi}{q,\partial S} \non\\
  &=&\spu{\vp}{T\Psi}{q,\partial S} =0\ ,\non
\eeqa
also $\chn\vp\in D^{q+1,\Gamma_2}(S)$.
Lemma \ref{ninhnull2} und der Sternoperator liefern
die \"ubrigen Aussagen.\qed
\abstandk
F\"ur $E\in R^q(S)$ definieren wir die Tangentialspur\label{pxteilrotspur}
$\gt^{\Gamma_1}E$ in $\Gamma_1$ durch
\beqa
 \alle{\vp\in H_{1/2}^{q,\Gamma_2} (\partial S)}
 \spu{\gt^{\Gamma_1} E}{\vp}{H_{-1/2}^{q,\Gamma_2} (\partial S)}
  :=\spu{\gt E}{\vp}{H_{-1/2}^q(\partial S)} \ .\non
\eeqa
Diese hat folgende Eigenschaften:
\begin{satz}  \formel{rteilspur} \new
 F\"ur jedes $E\in R^q(S)$ liegt die Tangentialspur
$\gt^{\Gamma_1}E$ in ${R}^{q,\Gamma_2}_{-1/2}(\partial S)$ und es gelten:
 \begin{enumerate}
     \item[i)] $\alle{E\in R^q(S)}
     \rot\gt^{\Gamma_1} E =\gt^{\Gamma_1} \rot E$
  \item[ii)] Die Abbildung
  $ \gt^{\Gamma_1}:R^q(S)\ra R^{q,\Gamma_2}_{-1/2}(\partial S)$
  ist stetig.
  \item[iii)] $E\in R^{q,\Gamma_1} (S) \Lra \gt^{\Gamma_1} E =0 $.
 \end{enumerate}
\end{satz}
\beweis
Da $\gt^{\Gamma_1}$ eine Einschr\"ankung des Funktionals
$\gt$ ist folgen i) und ii) aus Satz \ref{rotspursatz}.
\abstandk
iii) F\"ur $E\in C_\infty^{q,\Gamma_1}(\ol{S})$ und
$\vp\in  H_{1/2}^{q,\Gamma_2} (\partial S)$
gilt nach Lemma \ref{ninhnull}
\beqa
 \spu{\gt^{\Gamma_1}E}{\vp}{ H_{-1/2}^{q+1,\Gamma_2} (\partial S)}
  =\spu{\rot E}{\chn \vp}{q+1,S}
   +\spu{E}{\div \chn \vp}{q,S}
    =0\ .\non
\eeqa
Andererseits impliziert $\gt^{\Gamma_1}E=0$ f\"ur
$\Phi\in C_\infty^{q+1,\Gamma_2}(\ol{S})$
($\Ra N\Phi\in H_{1/2}^{q,\Gamma_2} (\partial S)$)
\beqa
 \spu{\rot E}{\Phi}{q+1,S} +
 \spu{E}{\div\Phi}{q,S}
 &=&\spu{\gt E}{N\Phi}{ H_{-1/2}^{q} (\partial S)} \non\\
 &=&\spu{\gt^{\Gamma_1} E}{N\Phi}{ H_{-1/2}^{q,\Gamma_2} (\partial S)} \non\\
 &=&0 \ .\non
\eeqa
\qed
\abstandk
Entsprechende Aussagen erhalten wir f\"ur
die Teilspur von Formen $E$ aus $D^q(S)$.
Erkl\"aren wir den Sternoperator analog zu
$(\ref{sternrot})$ und\label{pxdteilspur}
\beqa
  \gn^{\Gamma_2} E&:=&\sigma_{q} *\gt^{\Gamma_2}* E \non\\
   D_{-1/2}^{q,\Gamma_1}(\partial S)&:=&
   *R_{-1/2}^{N-1-q,\Gamma_1}(\partial S) \non
\eeqa
\label{pxhatdgam}
so folgen die dualen Resultate:
\begin{satz}    \formel{kuhnsatzd} \new
F\"ur jedes $E\in D^q(S)$ liegt die Normalenspur
$\gn^{\Gamma_2} E$ in ${D}^{q-1,\Gamma_1}_{-1/2}(\partial S)$ und es gelten:
\begin{enumerate}
      \item[i)] $\alle{E\in D^q(S)}
      \div\gn^{\Gamma_2} E =-\gn^{\Gamma_2} \div E$
  \item[ii)] Die Abbildung
  $\gn^{\Gamma_2}:D^q(S)\ra {D}^{q-1,\Gamma_1}_{-1/2}(\partial S)$
  ist stetig.
      \item[iii)] $E\in D^{q,\Gamma_2} (S) \Lra \gn^{\Gamma_2} E =0 $.
\end{enumerate}
\end{satz}
\abstandk
\begin{bemerk} $\mbox{ }$ \new
Eine Konstruktion des Fortsetzungsoperators
wie in Abschnitt \ref{chrqspur} scheitert
im Falle wechselnder Randbedingungen daran,
da{\ss} Formen aus $R^{q,\Gamma_1}(S)\cap D^{q,\Gamma_2}(S)$
im allgemeinen nicht in $H_1^q(S)$ liegen. Diese
Eigenschaft kann man erwarten, wenn die Randst\"ucke
$\Gamma_1$ und $\Gamma_2$ senkrecht aufeinander treffen.
\end{bemerk}
\section{L\"osungstheorie} \label{chrwp}
\markboth{L\"OSUNGSTHEORIE}{L\"OSUNGSTHEORIE}
Seien die Bezeichnungen  wie in Abschnitt
\ref{chrqspur} und $S$ glatt.
Nach Satz
\ref{regularitaet} liegen Formen aus
$R^q(S)\cap D^q(S)$, welche die elektrische oder
magnetische homogene Randbedingung erf\"ullen, in $H_1^q(S)$.
Wir betrachten folgendes Problem:\new
Gesucht ist eine Form $E\in R^q(S)\cap D^q(S)$ mit
\beqa \left.\begin{array}{rcl}
 \rot E&=&F\\
 \div E&=&G\\
 \gt E&=&\lambda\\
 \spu{E}{y_j^q}{q,S} &=&\alpha_j\mfur j=1\cdots J^q \ .
 \end{array} \right\}
 \formel{problem1}
\eeqa
\begin{lemma}\formel{loesedivnull} \new
 Seien $F\in R_0^{q+1}(S)$ und $\lambda\in {R}_{-1/2}^q(\partial S)$
 mit
 \beqa
    \rot\lambda &=&\gt F  \formel{vorein}   \\
    \spu{F}{y}{q+1,S} &=& \spu{\lambda}{Ny}{H_{-1/2}^q(\partial S)}
    \mfur\malle y\in Y_3^{q+1}\ .\non
 \eeqa
 Dann existiert eine L\"osung $E\in R^q(S)\cap D^q(S)$ von
 \beqa \left.\begin{array}{rcl}
  \rot E&=& F\\
  \div E&=& 0\\
  \gt E&=&\lambda \ .
  \end{array}
  \right\} \formel{problem2}
 \eeqa
 Diese erf\"ullt
 \beqa
   \normu{E}{R^q(S)\cap D^q(S)}\le
       c(\normu{\lambda}{R_{-1/2}^q(\partial S)}
    +\normu{F}{q+1,S}) \ .\non
 \eeqa
\end{lemma}
\beweis  Wir suchen zun\"achst eine Form $Q\in Y_1^{q+1}$ mit
\beqa
 \alle{\Phi\in Y_1^{q+1}}
    \spu{\div Q}{\div\Phi}{q,S}
    = \spu{\lambda}{N\Phi}{H_{-1/2}^q(\partial S)}
    -\spu{F}{\Phi}{q+1,S} \ .\formel{problem4}
\eeqa
Nach Korollar \ref{rqpoincare} ist die linke Seite
eine stetige streng koerzitive Bilinearform \"uber $Y_1^{q+1}$.
Die rechte Seite ist ein antilineares stetiges Funktional.
Die nach dem Satz von Lax--Milgram existierende L\"osung
$Q$ erf\"ullt
\beqa
 \normu{Q}{D^{q+1}(S)}\le c(\normu{\lambda}{H_{-1/2}^q(\partial S)}
    +\normu{F}{q+1,S}) \ ,\non
\eeqa
und nach Voraussetzung gilt
$(\ref{problem4})$ sogar f\"ur alle
$\Phi\in Y_1^{q+1}\oplus Y_3^{q+1} $.
F\"ur $\Phi\in Y_2^{q+1}$ mit $\Phi=\div\Psi$ und
$\Psi\in Y_1^{q+2}\cap H_2^{q+2}(S)$
(vgl. Abschnitt \ref{chrqspur})
folgt aus Satz \ref{rotspursatz} und $(\ref{vorein})$
\beqa
  \spu{F}{\div\Psi}{q+1,S}
   &=&\spu{\gt F}{N\Psi}{H_{-1/2}^{q+1}(\partial S)} \non\\
   &=&\spu{\rot\lambda}{N\Psi}{H_{-1/2}^{q+1}(\partial S)} \non\\
   &=&\spu{\lambda}{N\div\Psi}{H_{-1/2}^q(\partial S)} \ ,\non
\eeqa
also auch $(\ref{problem4})$ f\"ur alle
$\Phi\in Y^{q+1}\supset \chn H_{1/2}^q(\partial S)$.
Somit ist $E=\div Q$
L\"osung zu $(\ref{problem2})$.
\qed
\begin{lemma} \formel{loeserotnull} \new
 Sei $G\in \div D^{q}(S)$.
 Dann existiert eine L\"osung $E\in R^q(S)\cap D^q(S)$ von
 \beqa \left.\begin{array}{rcl}
  \rot E&=& 0\\
  \div E&=& G\\
  \gt E&=&0 \ .
   \end{array}
  \right\} \formel{problem3}
 \eeqa
 Diese erf\"ullt
 \beqa
    \normu{E}{D^q(S)}&\le&c\normu{G}{q-1,S} \ .\non
 \eeqa
\end{lemma}
\beweis
Nach Korollar \ref{rqpoincare} k\"onnen wir ein
$Z\in Y_2^{q-1}$ finden mit
\beqa
 \alle{\Phi\in Y_2^{q-1}}
  \spu{\rot Z}{\rot\Phi}{q,S} &=& -\spu{G}{\Phi}{q-1,S}
  \formel{problem5} \\
   \normu{Z}{R^{q-1}(S)} &\le& c\normu{G}{q-1,S} \ .\formel{zstetig}
\eeqa
Wegen $G\in (Y_3^{q-1})^\perp$ ist
$(\ref{problem5})$ auch f\"ur alle
$\Phi\in Y_3^{q-1}$
erf\"ullt. F\"ur $\Phi\in Y_1^{q-1}$ gilt
die Gleichung $(\ref{problem5})$ wegen der Voraussetzung
an $G$. Wir folgern
\beqa
 \alle{\Phi\in Y^{q-1}}
  \spu{\rot Z}{\rot\Phi}{q,S} &=& -\spu{G}{\Phi}{q-1,S} \ .\non
\eeqa
Setzen wir
 \beqa
   E:=\rot Z\in \rot\rqmn(S) \subset \rnull{0 }{q}(S)\ ,\formel{randbe}
\eeqa
so folgt aus
$\cqn(S)\subset Y^{q-1}$,  $(\ref{zstetig})$ und
$(\ref{spurnull})$ die Behauptung.
\qed
\begin{satz}  \formel{loesereins} \new
 Notwendige und hinreichende Bedingungen f\"ur die L\"osbarkeit
 von $(\ref{problem1})$ sind
 \beqa
        F\in\rot R^q(S)\ ,
 \ G\in\div D^q(S)\ , \ \lambda\in {R}_{-1/2}^q(\partial S) \non
\eeqa
 und
 \beqa
   \rot\lambda&=&\gt F\non\\
   \alle{ y\in Y_3^{q+1}}
   \spu{F}{ y}{q+1,S} &=& \spu{\lambda}{N y}{
   H_{-1/2}^q(\partial S)} \ .\non
 \eeqa
 Das Problem ist dann
 eindeutig l\"osbar, und die L\"osung erf\"ullt
 \beqa
  \normu{E}{R^q(S)\cap D^q(S)} \le c(\normu{F}{q+1,S}
   +\normu{G}{q-1,S} +\normu{\lambda}{R^q_{-1/2}(\partial S)}
   +\sum_{j=1}^{J^q}|\alpha_j|) \ .\non
 \eeqa
\end{satz}
\beweis  Nach $(\ref{spurnull})$ liegt die L\"osung des
 homogenen Problems in $\rnull{0}{q} (S)\cap D_0^q(S)$.
Aus  $\spu{E}{ y_j^q}{q,S} =0$, $j=1,\cdots, J^q$ folgt $E=0$ und
  damit die Eindeutigkeit. \new
Seien $E_1$ und $E_2$ L\"osungen der Probleme
$(\ref{problem2})$ und $(\ref{problem3})$.
Setzen wir
\beqa
 E:=E_1+E_2+\sum_{k=1}^{J^q}(\alpha_k-\spu{E_1}{ y_k^q}{q,S})
   y_j^q \ ,\non
\eeqa
so gelten $\rot E=F$, $\div E=G$, und wegen $(\ref{randbe})$
und $(\ref{spurnull})$
ist die Randbedingung erf\"ullt.
Aus
\beqa
 \spu{E}{ y_j^q}{q,S} &=& \spu{E_1}{ y_j^q}{q,S}
  +(\alpha_j-\spu{E_1}{ y_j^q}{q,S}) =\alpha_j \non
\eeqa
folgt, da{\ss} $E$ das Problem l\"ost
(hier haben wir $E_2\in\rot\rqmn(S)$
nach $(\ref{randbe})$ benutzt). \new
Es verbleibt, die Notwendigkeit der Voraussetzungen zu zeigen.
Da{\ss} im Falle der L\"osbarkeit $F\in\rot R^q(S)$ und
$G\in \div D^q(S)$ gelten m\"ussen, ist klar.
Aus Satz \ref{rotspursatz} folgt
$\lambda\in {R}_{-1/2}^q(\partial S)$.
\abstandk
Nach Satz \ref{rotspursatz} ii) gilt
\beqa
 \rot\lambda =\rot\gt E =\gt\rot E
 =\gt F \non
\eeqa
und f\"ur $ y\in  \rnull{0}{q+1} (S)\cap D_0^{q+1} (S)$
\beqa
 \spu{F}{ y}{q+1,S}
  &=&\spu{\rot E}{ y}{q+1,S} +\spu{E}{\div  y}{q,S} \non\\
  &=&\spu{\gt E}{N y}{H_{-1/2}^q(\partial S)} \non\\
  &=&\spu{\lambda}{N y}{H_{-1/2}^q(\partial S)} \ .\non
\eeqa
\qed
\abstandk
Wir formulieren noch das duale Problem:
Gesucht
$E\in R^q(S)\cap D^q(S)$ mit
\beqa \left.\begin{array}{rcl}
 \rot E&=&F\\
 \div E&=&G\\
 \gn E&=&\lambda\\
 \spu{E}{* y_j^{N-q}}{q,S} &=&\alpha_j\mfur
  j=1\cdots J^{N-q} \ .
 \end{array} \right\}
 \formel{problemd1}
\eeqa
Das duale Resultat lautet dann:
\begin{satz}  \formel{loesereinsd} \new
 Notwendige und hinreichende Bedingungen f\"ur die L\"osbarkeit
 von $(\ref{problemd1})$ sind
 \beqa
        F\in\rot R^q(S)\ ,
 \ G\in\div D^q(S)\ ,\ \lambda\in {D}_{-1/2}^{q-1}(\partial S)
 \non
\eeqa
 und
 \beqa
   \div\lambda&=&-\gn G\non\\
   \alle{ y\in *Y_3^{N-(q-1)}}
   \spu{G}{ y}{q-1,S} &=& \spu{\lambda}{T y}{
   H_{-1/2}^{q-1}(\partial S)} \ .\non
 \eeqa
 Das Problem ist dann
 eindeutig l\"osbar, und die L\"osung erf\"ullt
 \beqa
  \normu{E}{R^q(S)\cap D^q(S)} \le c(\normu{F}{q+1,S}
   +\normu{G}{q-1,S} +\normu{\lambda}{D^{q-1}_{-1/2}(\partial S)}
   +\sum_{j=1}^{J^{N-q}}|\alpha_j|) \ .\non
 \eeqa
\end{satz}
\section{Die Dirichlet--Neumann Felder} \label{chfelder}
\markboth{DIE DIRICHLET -- NEUMANN FELDER}{DIE DIRICHLET -- NEUMANN FELDER}
\achtung{Wurden untersucht in \cite{kress} mit gemischten Randbedingungen
im $\rrr^3$ und in \cite{picard2} auf Mannigfaltigkeiten. }
\begin{satz}    \label{dirineumdim}
 Sei $(S,\Gamma_1,\Gamma_2)$ glatt und
 $\Gamma_1={\bigcup}_{k=1}^K\Gamma_{1,k}$,
 mit $\dist(\Gamma_{1,k},\Gamma_{1,l})>0$ f\"ur $l\not=k$ und
 $\Gamma_{1,k}$ offen, nicht leer und zusammenh\"angend.
 Ferner m\"oge die erste
 Betti--Zahl von $S$ verschwinden. Dann gilt
 \beqa
  \dim (R_0^{1,\Gamma_1}(S)\cap D_0^{1,\Gamma_2}(S)) =K-1 \ .\non
 \eeqa
\end{satz}
\abstandk
Um den Satz zu beweisen, ben\"otigen wir
zun\"achst einige Vorbereitungen:
\begin{lemma} \formel{nase} \new
  Aus den Voraussetzungen von Satz \ref{dirineumdim}
   folgt die Existenz
offener Mengen $S_k\subset  M \setminus \ol{S}$ mit
$\ol{S_k}\cap\ol{S}=\ol{\Gamma}_{1,k}$, $k=1,\cdots,K$ und
$\dist(S_k,S_l)>0$, $k\not= l$.
Wir k\"onnen die Mengen $S_k$ so w\"ahlen,
da{\ss} die Glattheitseigenschaft f\"ur
\beqa
 \hat{S}:=S\cup\Gamma_1\cup\bigcup_{k=1}^{K} S_k \non
\eeqa
erhalten bleibt
und
die erste Betti--Zahl von $\hat{S}$ verschwindet.\new
\end{lemma}
\beweis
Nach \cite[Theorem 1.1.7]{schwarz} oder auch
\cite{hirsch} existiert
ein Hom\"oomorphismus
$n$ von $\partial S\times (-1,1)$
auf eine in $M$ offene Umgebung
von $\partial S$ mit
\beqa
 \alle{x\in\partial S} n(x,0)=x \non
\eeqa
und der Eigenschaft, da{\ss} $n_{|_{\partial S\times (-1,0]}}$
bzw. $n_{|_{\partial S\times [0,1)}}$
Diffeomorphismen auf in $\ol{S}$ bzw. $M\setminus S$ offene
Umgebungen von $\partial S$
sind.
F\"ur eine Funktion $h\in C_\infty(\partial S)$
mit Werten in $[0,1/2]$ und $h_{|_{\Gamma_1}}>0$,
$h_{|_{\Gamma_2}}=0$ setzen wir
\beqa
 S_k&:=&\{n(x,th(x)) \mit 0<t<1,x\in\Gamma_{1,k}\} \non\\
  \tilde{S}&:=&\bigcup_{k=1}^{K} S_k\non\\
  \hat{S}&:=&S\cup\Gamma_1\cup \tilde{S} \non.
\eeqa
Gilt $\dist(S_k,S_l)=0$, so existieren Folgen
$(t_n),(s_n)\subset [0,1]$, $(x_n)\subset \Gamma_{1,k}$ und
$(y_n)\subset \Gamma_{1,l}$ mit
\beqa
 d_M(n(x_n,t_nh(x_n)),n(y_n,s_nh(y_n)))\ra 0 \ .\non
\eeqa
Dies impliziert $d_M(x_n,y_n)\ra 0$. Wegen
$\dist(\Gamma_{1,k},\Gamma_{1,l})>0$ f\"ur $k\not=l$
m\"ussen $k$ und $l$ \"ubereinstimmen. \"Ahnlich
zeigt man $\ol{S_k}\cap \ol{S}=\ol{\Gamma}_{1,k}$.
Da alle Punkte aus $\Gamma_1$ innere Punkte der Menge
$\hat{S}$ sind, ist $\hat{S}$ offen, und es gilt
wegen $S\cap \ol{\tilde{S}}=\ol{S}\cap\tilde{S}=\emptyset$
\beqa
 \partial\hat{S}&=&\ol{\hat{S}}\setminus \hat{S}\non\\
        &=&(\ol{S}\cup \ol{\tilde{S}})\setminus
        (S\cup\Gamma_1\cup\tilde{S})\non\\
       &=& (\partial {S}\setminus \Gamma_1)
        \cup  (\partial\tilde{S}\setminus \Gamma_1)\ .\non
\eeqa
Da $\tilde{S}$ das Bild der Menge
\beqa
 W:=\{(x,t)\in\Gamma_1 \times \rrr \mit 0<t<h(x)\} \non
\eeqa
unter dem Diffeomorphismus $n$ ist, erhalten wir mit
$\gamma:=\partial\Gamma_1=\partial\Gamma_2$
\beqa
 \partial\tilde{S} &=& n(\partial W)\non\\
 &=& n(\Gamma_1\times \{0\})\cup \{n(x,h(x)) \mit x\in\Gamma_1\}
  \cup n(\gamma\times \{0 \}) \non\\
 &=& \Gamma_1\cup  \{n(x,h(x)) \mit x\in\Gamma_1\} \cup\gamma\ .\non
\eeqa
Es folgt
\beqa
 \partial \hat{S}&=&\Gamma_2\cup \gamma\cup \{n(x,h(x)) \mit x\in\Gamma_1\}   \non\\
 &=& \{n(x,h(x)) \mit x\in\partial S\} \non
\eeqa
und damit die Glattheit des Randes $\partial\hat {S}$.
Setzen wir
\beqa
 S_k^-&:=&\{n(x,th(x)) \mit -1<t<0,x\in\Gamma_{1,k}\} \ ,\non
\eeqa
so stellt die Abbildung
\beqa
 x\mapsto\left\{
 \begin{array}{lcl}
  x&\mfur&x\not\in\bigcup_{k=1}^K S_k^-\\
  n(y,(1+2t)h(y))&\mfur& x=n(y,th(y)), t\in (-1,0)
 \end{array}  \right.      \non
\eeqa
einen Hom\"oomorphismus zwischen den Mengen $S$ und $\hat{S}$
dar. Somit verschwindet auch die
erste Betti--Zahl von $\hat{S}$
(\cite[Seite 160]{chillingworth} oder auch \cite[Seite 18]{hu}
in Verbindung mit \cite[Seite 42]{rosenberg}).
\qed
\begin{lemma}  \formel{nullinnase}  \new
 Die Nullfortsetzung $\hat{E}$ einer Form $E$ aus $R^{q,\Gamma_1}(S)$
 liegt in $R^q(\hat{S})$.
\end{lemma}
\beweis
Wegen $\Gamma_2\subset  M \setminus \hat{S}$ liegen Elemente
aus $\cnull{\infty}{q+1} (\hat{S})$ auch in
$C_{\infty}^{q+1,\Gamma_2} (S)$.
Somit gilt f\"ur $E\in R^{q,\Gamma_1}(S) $,
$\Phi\in\cnull{\infty}{q+1} (\hat{S})$ und
die Nullfortsetzung $F$ von $\rot E$
\beqa
 \spu{\hat{E}}{\div\Phi}{q,\hat{S}}
  &=&\spu{E}{\div \Phi}{q,S} \non\\
  &=&-\spu{\rot E}{\Phi}{q+1,S} \non\\
  &=&-\spu{F}{\Phi}{q+1,\hat{S}} \ ,\non
\eeqa
also auch $\hat{E}\in R^q(\hat{S})$.
\qed
\begin{lemma} \formel{nullinnase2}    \new
 Verschwindet $\hat{E}\in R^q(\hat{S})$ in $\hat{S}\setminus S$,
 so liegt
 $E:=\hat{E}_{|_S}$ in $R^{q,\Gamma_1}(S)$.
\end{lemma}
\beweis
Wegen $\partial S=\Gamma_1\cup \ol{\Gamma}_2$ und $S_k\subset  M \setminus\ol{S}$
gilt
\beqa
 \ol{S}\cap( M \setminus\hat{S})
  &=& \ol{S}\setminus(S\cup\Gamma_1\cup\bigcup_{k=1}^KS_k)\non\\
  &=&\partial S \cap \ol{\Gamma}_2\cap
   \ol{S}
  =\ol{\Gamma}_2 \ .\non
\eeqa
F\"ur eine Form $\Psi\in C_\infty^{q+1,\Gamma_2}(\ol{S})$
folgen
\beqa
 \supp\Psi\cap\ol{S}\cap ( M \setminus \hat{S}) &=&\emptyset \non\\
 \dist (\supp\Psi\cap\ol{S}, M \setminus\hat{S})&>&0\ \non
\eeqa
($\supp\Psi\cap\ol{S}$ ist kompakt, $ M \setminus\hat{S}$
abgeschlossen). Mit einer Funktion
$\chi\in\cnull{\infty}{\ } (\hat{S}) $,
$\chi=1$ in $\supp\Psi\cap\ol{S}$ gelten
$\chi\Psi\in \cnull{\infty}{q+1} (\hat{S})$
(nach Definition ist $\Psi$ auf ganz $M$ erkl\"art) und
\beqa
 \spu{E}{\div\Psi}{q,S}
  &=&\spu{\hat{E}}{\div\chi\Psi}{q,\hat{S}} \non\\
  &=&-\spu{\rot \hat{E}}{\chi\Psi}{q+1,\hat{S}} \non\\
  &=&-\spu{\rot E}{\Psi}{q+1,S} \ .\non
\eeqa
Wir erhalten $E\in R^{q,\Gamma_1}(S)$.
\qed
\abstandk
{\bf Beweis} von Satz \ref{dirineumdim}:\new
Wir verwenden die Bezeichnungen aus Lemma \ref{nase}.
F\"ur Funktionen $\Psi_k\in\cnull{\infty}{\ } ( M ) $ mit
$\Psi_k=\delta_{k,l}$ in $S_l$ und
\beqa
 \dist(\supp\rot\Psi_k, S_l)>0 \formel{distsupppsi}
\eeqa
f\"ur alle $k,l=1,\cdots, K$
existieren nach Korollar \ref{rqpoincare} und
dem Satz von Lax--Milgram Funktionen
$\eta_k$ mit $\eta_k-\Psi_k\in R^{0,\Gamma_1}(S)$
und
\beqa
 \alle{\Phi\in R^{0,\Gamma_1}(S)}
 \spu{\rot(\eta_k-\Psi_k)}{\rot\Phi}{1,S}
 =-\spu{\rot\Psi_k}{\rot\Phi}{1,S} \formel{etakphi}
\eeqa
(Beachte $R^{0,\Gamma_1}(S)=
R^{0,\Gamma_1}(S)\cap D_0^{0,\Gamma_2}(S)$).
Wir behaupten
\beqa
 \{\rot\eta_k \mit k=1,\cdots,K-1\} \formel{basis}
\eeqa
ist eine Basis von $R_0^{1,\Gamma_1}(S)\cap D_0^{1,\Gamma_2}(S)$.
Aus $(\ref{etakphi})$ folgt sofort
$\rot\eta_k\in D_0^{1,\Gamma_2}(S)$.
Wegen
\beqa
 \rot\eta_k=\rot(\eta_k-\Psi_k)+\rot\Psi_k \non
\eeqa
und $(\ref{distsupppsi})$ gilt
$\rot\eta_k\in R_0^{1,\Gamma_1}(S)$
(wie in Lemma \ref{supprandver} zeigt man
$E\in R^{q,\Gamma_1}(S)$ f\"ur Formen $E\in R^q(S)$ mit
$\dist(\supp E,\Gamma_1)>0$).
Um zu zeigen, da{\ss} $(\ref{basis})$ ein Erzeugendensystem
ist, w\"ahlen wir
$F\in R_0^{1,\Gamma_1}(S)\cap D_0^{1,\Gamma_2}(S) $. Die
Nullfortsetzung $\hat{F}$ liegt nach
Lemma \ref{nullinnase} in $R_0^1(\hat{S})$. Da die erste
Betti--Zahl von $\hat{S}$ verschwindet, existiert nach
\cite[Satz 1, Satz 2]{picard2} und Korollar \ref{rqpoincare}
eine Form $G\in R^0(\hat{S})$ mit
$\rot G=\hat{F}$. Aus $\rot G=0$ in $S_k$ und
$(\ref{risth})$ folgt
\beqa
 G(x)=c_k \mfur x\in S_k\ .\formel{hatgisnull}
\eeqa
Sei o.B.d.A. $c_K=0$. Wir definieren
\beqa
 \mu&:=&F-\sum_{k=1}^{K-1}c_k\rot\eta_k
  \in R_0^{1,\Gamma_1}(S)\cap D_0^{1,\Gamma_2}(S)\formel{muinrnull}\\
  H&:=& G-\sum_{k=1}^{K-1}c_k\Psi_k \ .\non
\eeqa
Wegen $(\ref{hatgisnull})$ verschwindet $H$ in $\hat{S}\setminus S$,
und wir erhalten nach Lemma \ref{nullinnase2}
\beqa
 H\in R^{0,\Gamma_1}(S)\ .\non
\eeqa
In $S$ gilt
\beqa
 \rot H= F-\sum_{k=1}^{K-1}c_k\rot\Psi_k
     =\mu
     -\sum_{k=1}^{K-1}c_k\rot(\Psi_k-\eta_k) \ ,\non
\eeqa
und es folgt
\beqa
 \mu=\rot H +\sum_{k=1}^{K-1}c_k\rot(\Psi_k-\eta_k)
  \in \rot R^{0,\Gamma_1}(S)\ ,\non
\eeqa
also $\mu=0$ (Lemma \ref{orthzerl} und  $(\ref{muinrnull})$).
Es verbleibt, die lineare Unabh\"angigkeit zu zeigen. Gelte
daher
\beqa
 \sum_{k=1}^{K-1}\alpha_k\rot\eta_k=0 \ .\non
\eeqa
Wir setzen
\beqa
 \hat{\eta_k}:=\left\{  \begin{array}{ccl}
 \delta_{k,l}&\mbox{ in }& S_l\cup\Gamma_{1,l} \ \mfur
 \  l=1,\cdots, K\\
 \eta_k&\mbox{ in }& S
 \end{array}\right.  \non
\eeqa
f\"ur $k=1,\cdots,K-1$. In $S_l\cup\Gamma_{1,l}$ gilt
$0=\Psi_k-\hat{\eta}_k$, und aus $\Psi_k-\eta_k\in R^{0,\Gamma_1}(S)$
folgt $\Psi_k-\hat{\eta}_k\in R^{0}(\hat{S})$, also auch
$\hat{\eta}_k\in R^{0}(\hat{S})$ und
\beqa
 \sum_{k=1}^{K-1}\alpha_k\rot\hat{\eta}_k=0 \ .\non
\eeqa
Dies impliziert
\beqa
 \sum_{k=1}^{K-1}\alpha_k\hat{\eta}_k=c \non
\eeqa
in $\hat{S}$. Wegen $\hat{\eta}_k=0$ in $S_K$
mu{\ss} die Konstante $c$ verschwinden, und
wir erhalten f\"ur
$l=1,\cdots, K-1$
\beqa
 0=\sum_{k=1}^{K-1}\alpha_k\hat{\eta}_k
  =\alpha_l \mbox{ in } S_l\ .\non
\eeqa
\qed
\section{Eigenformen } \label{cheigenformen}
\markboth{EIGENFORMEN}{DIE HALBE KREISLINIE}
Um Aussagen \"uber die G\"ute von L\"osungen der Maxwellgleichung
machen zu k\"onnen, kann man sich der Methode von Saranen
\cite{saranen} bedienen (in Kegelgebieten im $\rrr^3$).
Hierzu ist es notwendig, Eigenformen niederer
Dimensionen zu kennen.
Wir wollen f\"ur spezielle Kegelgebiete
$(S,\Gamma_1,\Gamma_2)\in\calm(M)$ die
Ortho\-normal\-systeme
aus Lemma \ref{wwp15}, die wir hier
unter Ber\"ucksichtigung der Raumdimension mit
$\{E_n^{N,q}\}$ bzw. $\{H_n^{N,q} \}$ bezeichnen, und
die Eigenwerte -- hier $\omega_{n}^{N,q}$ --
f\"ur $\eps=\id$ und $\mu=\id$ berechnen.
Dar\"uber hinaus untersuchen wir deren Regularit\"at.
Da{\ss} wir nach den Eigenformen
entwickeln k\"onnen, folgt aus Lemma \ref{wwp15}, wenn
wir die kompakte Einbettung
\beqa
 R^{q,\Gamma_1}(S)\cap D^{q,\Gamma_2}(S)\hookrightarrow
 L_2^q(S) \formel{hook}
\eeqa
und die Approximationseigenschaft
$(S,\Gamma_1,\Gamma_2)\in\calm_D(M)$ zeigen k\"onnen.
\subsection{Die halbe Kreislinie} \label{kreislinie}
Wir beginnen mit dem halben Kreisrand
\beqa
 K&:=&\{ \tau(\vp) \mit \vp\in(0,\pi) \}
 \ ,\ \gamma_1:=\tau(\pi)\ ,\ \gamma_2:=\tau(0) \non\\
 \mmit \tau(\vp)&:=& \zweivec{\cos\vp}{\sin\vp} \ .\non
\eeqa
Aus Lemma \ref{indanfang} folgen
$(K,\gamma_1,\gamma_2)\in\calm_D(S_2)$
und die Kompaktheit der  Einbettung in $(\ref{hook})$.
Sei zun\"achst $q=0$.
Aus $(\ref{risth})$ folgt, da{\ss} die Dirichlet--Neumann--Felder
verschwinden, d.h. $N^0=0$ in der
Terminologie von Lemma \ref{wwp15}.
Die Eigenformen $(E,H)$
erf\"ullen
\beqa
  \rot E+\i\omega H&=&0 \non\\
  \div H+\i\omega E&=&0 \ .\non
\eeqa
Mit der Koordinate $\vp:=\vp(x):=\arccos(x/|x|)$,
$E=e(\vp)$ folgt f\"ur $\Psi\in C_\infty^{0,\gamma_1}(\ol{K})$,
$\Psi=\psi(\vp)$ mit $(\ref{intphi})$ und $(\ref{rotlokal})$
\beqa
 0&=&\spu{\rot E}{\rot\Psi}{1,K}
  +\i\omega\spu{H}{\rot\Psi}{1,K} \non\\
  &=&\spu{\rot E}{\rot\Psi}{1,K}
  -\i\omega\spu{\div H}{\Psi}{0,K} \non\\
  &=&\spu{\rot E}{\rot\Psi}{1,K}
   -\omega^2\spu{E}{\Psi}{0,K} \non\\
  &=&\int_0^\pi e'(\vp)\psi'(\vp) d\vp
   -\omega^2\int_0^\pi e(\vp)\psi(\vp) \dvp \ .\non
\eeqa
Die Regularit\"atstheorie z.B. \cite[Theorem 6.2]{agmon}
liefert
\beqa
 e\in C_\infty(0,\pi)\ ,\ e''(\vp)+\omega^2 e(\vp)=0 \ .\non
\eeqa
Alle L\"osungen dieser Differentialgleichung haben die Form
\beqa
 e(\vp)=a\cos(\omega\vp)+b\sin(\omega\vp) \ .\formel{edgl5}
\eeqa
Testen wir mit $\Psi=\psi(\vp)\dvp\in C_\infty^{1,\gamma_2}(\ol{K})$,
$\psi(\pi)=1$, so impliziert
$E\in R^{0,\gamma_1}(K)$
\beqa
 0&=&\spu{E}{\div\Psi}{0,K}
   +\spu{\rot E}{\Psi}{1,K} \non\\
  &=&\int_0^\pi e(\vp)\psi'(\vp)\dvp
   +\int_0^\pi e'(\vp)\psi(\vp)\dvp =e(\pi) \ .\formel{randpi}
\eeqa
Wegen $\rot E=-\i\omega H\in D^{1,\gamma_2}(K)$ gilt f\"ur
$\Psi=\psi(\vp)\in C_\infty^{0,\gamma_1}(\ol{K})$, $\psi(0)=1$
\beqa
 0&=&\spu{\rot E}{\rot\Phi}{1,K}
  +\spu{\div\rot E}{\Phi}{0,K} \non\\
 &=&\int_0^\pi e'(\vp)\psi'(\vp)\dvp
  +\int_0^\pi e''(\vp)\psi(\vp)\dvp
 =-e'(0) \ .\formel{randpi2}
\eeqa
Mit $(\ref{edgl5})$, $(\ref{randpi})$, $(\ref{randpi2})$
und der Maxwellgleichung
erhalten wir Teil i) des folgenden Satzes
\begin{satz} \formel{dummy6} \new
\begin{enumerate}
 \item[i)] Bis auf normierende Konstanten gelten
\beqa
 E_n^{1,0}&=&\cos((n-\eh)\vp) \non\\
 H_n^{1,0}&=&-\i\sin((n-\eh)\vp) \dvp \non\\
 \omega_{n}^{1,0}&=&n-\eh \ \ \mfur n\in\nnn \ .\non
\eeqa
\item[ii)] F\"ur $q=1$ sind die Orthonormalsysteme leer.
\end{enumerate}
\end{satz}
\beweis ii) Aus $(\ref{risth})$ folgt, da{\ss} der Raum
$D_0^{1,\gamma_2}(K)$ nur aus der Nullabbildung  besteht. \qed
\subsection{Der Halbkreis}
\markright{DER HALBKREIS}
Mit den Bezeichnungen aus Abschnitt \ref{kreislinie}
sei
\beqa
 (S,\Gamma_1,\Gamma_2)=(C(K),C(\gamma_1)\cup K
       \cup\{(-1,0)\} ,C(\gamma_2))\ ,\non
\eeqa
der obere Halbkreis mit Radius 1 um den Ursprung.
Wir zeigen zun\"achst
\begin{lemma}    \formel{dichtekompakt}
  \begin{enumerate}
  \item[i)] Es gilt die Approximationsaussage
           $(S,\Gamma_1,\Gamma_2)\in\calm_D$.
  \item[ii)] Die Einbettung
  $R^{q,\Gamma_1}(S)\cap D^{q,\Gamma_2}(S)\hookrightarrow L_2^q(S)$
  ist kompakt.
  \end{enumerate}
\end{lemma}
\beweis i)
 F\"ur eine Karte $(V,h)$ um $x\in\ol{S}$ seien
\beqa \left.\begin{array}{rcl}
 G&:=&h(S\cap V)\\
 {\hat{\Gamma}}_i&:=&h(\Gamma_i\cap V)\\
  W&:=&U_2(1/3)\cap \ol{G}\ . \end{array} \right\}
 \formel{bezoben}
\eeqa
Es gen\"ugt wieder,
eine Form $E\in R^{q,{\hat{\Gamma}}_1} (G)$
mit kompaktem Tr\"ager in $W$
durch Formen aus $C_\infty^{q,\hat{\Gamma}_1} (\ol{G})$
mit kompaktem Tr\"ager in $W$ zu approximieren.
F\"ur alle $x\not=\gamma_2$
folgt dies aus $(\ref{lemm2})$ und
Satz \ref{dichtesatz}.
Zu $\gamma_2$ k\"onnen wir eine Karte $(V,h)$ finden
mit
\beqa \left. \begin{array}{rcl}
 G&=&\{x\in U_2 \mit x_1<0,x_2<0\} \\
 {\hat{\Gamma}}_1&=&\{x\in U_2 \mit x_1=0,x_2<0\}\\
\hat{\Gamma}_2&=&\{x\in U_2 \mit x_1<0,x_2=0\} \ .
\end{array} \right\} \formel{umgeb}
\eeqa
In
$\{x\in U_2 \mit x_1>0,x_2<0\}$
setzen wir wie in Lemma $\ref{nullinnase}$
$E$
zu Null fort, verschieben den Tr\"ager in
Richtung
$(-1,0)$ und approximieren dort
nach Lemma \ref{lemm3} durch Formen aus
$C_\infty^{q,{\hat{\Gamma}}_1}(\ol{G})$.
Nach Multiplikation mit einer geeigneten Abschneidefunktion
erhalten wir
$(S,\Gamma_1,\Gamma_2)\in\calm_D(S)$.
\abstandk
ii) Mit den
Bezeichnungen aus $(\ref{bezoben})$ zeigen wir
f\"ur eine Karte $(V,h)$ um $x$ und eine in
$R^{q,{\hat{\Gamma}}_1} (G)\cap
D^{q,{\hat{\Gamma}}_2} (G) $
beschr\"ankte Folge $(E^n)$
mit $\supp E^n \subset\subset W$,
da{\ss}  diese eine in $L_2^q(G)$ konvergente
Teilfolge besitzt.
F\"ur alle $x\not=\gamma_2$ folgt dies aus \cite{weck}
bzw. Satz \ref{einistkompakt}. Im Falle
$x=\gamma_2$ w\"ahlen wir die Umgebung $V$ mit $(\ref{umgeb})$.
Mit der Abbildung $\tau$ aus  Lemma \ref{spiegel}
gilt
$\mu_\tau=\eps_\tau=\id$ (siehe
$(\ref{transeins})$,$(\ref{taudx})$).
Wir setzen
\beqa
  \hat{G}&:=&G\cup\{x\in U_2 \mit x_1<0,x_2\ge 0\}\non
\eeqa
und zeigen f\"ur den (eingeschr\"ankten) Spiegelungsoperator
$S_{\mbox{{\tiny rot}}}$ aus Lemma \ref{spiegel}
\beqa
 S_{\mbox{{\tiny rot}}}E^n\in \rnull{\ }{q} (\hat{G})\cap D^q(\hat{G})\non\\
 \normu{S_{\mbox{{\tiny rot}}}E^n}{R^q(\hat{G})\cap D^q(\hat{G})} \le c \ .\non
\eeqa
Die Beschr\"anktheit folgt aus Lemma \ref{taudivrot}, die
Aussage $S_{\mbox{{\tiny rot}}}E^n\in D^q(\hat{G})$ durch
Approximation in $D^q(G)$ mit Formen aus
$C_\infty^{q,{\hat{\Gamma}}_2}
(\ol{G})$
(beachte $\eps=\id$).
Ebenfalls durch Approximation erhalten wir
$S_{\mbox{{\tiny rot}}}E^n\in R^{q,\hat{\Gamma}}(\hat{G})$ mit
\beqa
 \hat{\Gamma}:=
   \hat{\Gamma}_1\cup\{x\in U_2 \mit x_1=0,x_2\ge 0\}\ .\non
\eeqa
Die Lemmata \ref{taudivrot} und \ref{supprandver}
liefern $S_{\mbox{{\tiny rot}}}E^n\in \rnull{\ }{q}(\hat{G})$.
Nach $\cite{weck}$ besitzt $E^n$ eine
in $L_2^q(\hat{G})$,
damit auch in $L_2^q(G)$ konvergente Teilfolge. \qed
\abstandk
Wir definieren f\"ur $I:=(0,1)$
\beqa
 C_{1,q}&:=&\{\vp\in C_\infty((0,1]) \mit  \vp(1)=0,
   \hat{M}^{-1}\vp, \hat{M}^{-(q-1)}D\hat{M}^{q-1}\vp\in L_{2,2}(I) \} \non\\
 C_{2,q}&:=&\{\vp\in C_\infty((0,1)) \mit
  \hat{M}^{-(1-q)}D\hat{M}^{1-q}\vp,\hat{M}^{-1}\vp\in L_{2,2}(I) \} \non
\eeqa
Diese erf\"ullen
\begin{lemma}  \formel{dglkreis} \new
F\"ur die Koeffizienten aus
$(\ref{koeffizienten0})$, $(\ref{koeffizienten})$
zu einer Form
$F \in R^{q,\Gamma_1}(S)\cap D^{q,\Gamma_2}(S)$ gelten
\begin{enumerate}
 \item[i)] f\"ur alle $\vp\in C_{1,q}$
  \beqa
     \int_0^R r^{N-1}a_m(r)r^{-(q-1)}&&
         \hspace{-1.0cm}(r^{q-1}\vp(r))'dr\non\\
&=&-\i\omeg{q-1} \int_0^R r^{N-1}r^{-1}d_m(r)\vp(r)dr \non\\
 &&-\int_0^R r^{N-1}c_m^D(r)\vp(r)dr \mfur m>N^{q-1}\non
 \eeqa
 \item[ii)] f\"ur alle $\vp\in C_{2,q}$
 \beqa
  \int_0^R r^{N-1} d_m(r)r^{-(N-q-1)}&&\hspace{-1.0cm}
(r^{N-q-1}\vp(r))'dr\non\\
  &=&\i\omeg{q-1} \int_0^Rr^{N-1}r^{-1}a_m(r)\vp(r) dr \non\\
  &&-\int_0^R r^{N-1}b_m^R(r)\vp(r)dr
  \mfur m> N^{q-1}\ .\non
 \eeqa
\end{enumerate}
\end{lemma}
\beweis
Die beiden Aussagen k\"onnen wir ganz analog zu
Lemma \ref{dgl} beweisen, wenn wir zeigen, da{\ss}
die Formeln $(\ref{walzee})$ (f\"ur i)) und $(\ref{walzed})$
(f\"ur ii)) auch hier G\"ultigkeit haben.
Dazu vergleichen wir die Randbedingungen an $\vp$.
Beim Beweis von Lemma \ref{walze}
ist nicht benutzt worden, da{\ss} $\vp$ einen positiven Abstand zur
Null besitzt. Wir m\"ussen lediglich garantieren, da{\ss} die
von $r$ abh\"angigen Ausdr\"ucke auf den rechten Seiten
der Skalarprodukte so beschaffen sind, da{\ss} wir Lemma
\ref{approx} anwenden k\"onnen. Dies folgt aber aus den
Voraussetzungen an  $\vp$. \new
Ebenfalls entnehmen wir dem Beweis,
da{\ss} die Randbedingung $\vp(1)=0$ ausreichend ist. Dies liefert
i). Im Fall ii) \"ubernimmt $F\in {R}^{q,\Gamma_1}(S)$
die fehlende Randbedingung; der Deckel der Kegelspitze ist
hier in $\Gamma_1$ enthalten.
\qed
\abstandk
Wir untersuchen zun\"achst den Fall $q=0$.
Wie im Abschnitt \ref{kreislinie} liefert
$(\ref{risth})$, da{\ss} die Dirichlet--Neumann--Felder
verschwinden.
Seien $(E,H)$ L\"osungen von $(\ref{eigenloesung})$.
Mit den Resultaten
aus Abschnitt $\ref{kreislinie}$ zerlegen wir gem\"a{\ss}
$(\ref{ezerlegung})$
\beqa
 E&=&\chtau\sum_{n\ge 1} c_nE_n^{1,0} \formel{erstereihe}\\
 H&=&\chrho\sum_{n\ge 1} \tilde{a}_nE_n^{1,0}
   +\chtau\sum_{n\ge 1} \tilde{d}_nH_n^{1,0} \ ,
   \formel{zweitreihe}
\eeqa
wobei nach den Maxwellgleichungen die Beziehungen
\beqa
 a_n^R=-\i\omega\tilde{a}_n\ ,
 \ d_n^R=-\i\omega \tilde{d}_n\ ,
 \ \tilde{c}_n^{D}= -\i\omega c_n \formel{gleichu1}
\eeqa
gelten. Wir werden zeigen, da{\ss} die Koeffizienten
$c_n\in L_{2,2}(0,1)$
im Definitionsbereich des Abschlusses des Besseloperators
\label{pagebesselop1}
\beqa
 D(B^{\nu})&:=&\{v\in L_{2,2}(0,1) \mit v(r)=r^{\nu}u(r)\ ,
 \ u\in C_\infty([0,1])\ ,\ u'(0)=u(1)=0\} \non\\
 (B^{\nu}v)(r)&:=&-v''(r)-\frac{v'(r)}{r}+\frac{\nu^2}{r^2}
  v(r)\mmit \nu=\omega_n^{1,0} =n-\eh\ ,\ n\in\nnn\non
\eeqa
liegen.
Dieser ist nach \cite[Satz 27.2]{triebel} im Raum
$L_{2,2}(0,1)$ wesentlich
selbstadjungiert.
Funktionen $v\in D(B^\nu)$ erf\"ullen
$rv'\in C_{2,1}$ und $v\in C_{1,1}$.
Wir erhalten mit den Lemmata \ref{dgl} und \ref{dglkreis},
$(\ref{gleichu1})$
und $\tilde{b}_n^R=0$
\beqa
 \int_0^1 r(B^{\nu}v)(r)c_n(r)\dr
  &=&\i\nu^{-1}\int_0^1 r (B^{\nu}v)(r)rd_n^R(r) \dr \non\\
  &=&-\i\nu^{-1}\int_0^1r(rv'(r))'d_n^R(r) \dr
     +\i\nu^{-1}\int_0^1\frac{\nu^2}{r} v(r)r d_n^R(r)  \dr\non\\
&=&-\omega\nu^{-1}\int_0^1 r(rv'(r))'\tilde{d}_n(r) \dr
     +\i\nu\int_0^1 v(r) d_n^R(r)  \dr\non\\
      &=:&I_1+I_2\non\\
 I_1&=&-\i\omega\int_0^1\tilde{a}_n(r) rv'(r) \dr\non\\
  &=&-\nu\omega\int_0^1\tilde{d}_n(r)v(r)  \dr
   +\i\omega\int_0^1 r\tilde{c}_n^D(r) v(r) \dr \non\\
 &=& -I_2+\omega^2\int_0^1rv(r)c_n(r)\dr \ .\non
\eeqa
Bezeichnen wir mit
$J_\nu$ die Besselfunktion der Ordnung  $\nu$
\label{pagebesself} und mit
$\omega_{n,m}^{2,0}$ die positiven Nullstellen
der Besselfunktion $J_{n-1/2}$, so folgt aus
\cite[Satz 27.3]{triebel} bis auf normierende Konstanten
f\"ur ein $m\in\nnn$ entweder
\beqa
 \omega^2&=& (\omega_{n,m}^{2,0})^2 \non\\
 c_n(r)&=&c_{n,m}(r)=J_{n-\eh}(\omega_{n,m}^{2,0}r) \non\\
 \tilde{d_n}(r)&=&\tilde{d}_{n,m}(r)=\i(\omega_{n,m}^{2,0})^{-1}
  d_{n,m}^R(r) \non\\
   &&=\omega_n^{1,0}(\omega_{n,m}^{2,0})^{-1}r^{-1}c_{n,m}(r)
  \non\\
  \tilde{a}_n(r)&=&\tilde{a}_{n,m}(r)=\i(\omega_{n,m}^{2,0})^{-1}
   a_{n,m}^R(r) \non\\
   &&=\i(\omega_{n,m}^{2,0})^{-1} c_{n,m}'(r) \ ,\non
\eeqa
wobei die letzten Identit\"aten aus $(\ref{gleichu1})$ und
Lemma \ref{dgl} folgen,
oder $\omega>0$ beliebig und
\beqa
 c_n=\tilde{a}_n=\tilde{d_n}=0\ .\non
\eeqa
Da nach \cite[Seite 203]{porter} oder
\cite[Seite 485]{watson} die Besselfunktionen
$J_{n-1/2}$ f\"ur nat\"urliche Zahlen $n$ keine
gemeinsamen positiven Nullstellen haben, zerfallen die Reihen in
$(\ref{erstereihe})$ und $(\ref{zweitreihe})$ in eine
Komponente.
\abstandk
Im Falle $q=1$ k\"onnen wegen der unsymmetrischen
Randbedingungen nicht die dualen Resultate benutzt werden.
Da{\ss} die Dirichlet--Neumann--Felder verschwinden,
zeigt man analog zu Satz \ref{dirineumdim}.
Wie oben zerlegen
wir
\beqa
 E&=&\chrho\sum_{n\ge 1} \tilde{a}_nE_n^{1,0}
   +\chtau\sum_{n\ge 1} \tilde{d}_nH_n^{1,0}
   \formel{drittreihe}\\
    H&=&\chrho\sum_{n\ge 1} b_nH_n^{1,0} \ ,
    \formel{viertreihe}
\eeqa
erhalten die Gleichungen
\beqa
 a_n^D=-\i\omega\tilde{a}_n\ ,
 \ d_n^D=-\i\omega \tilde{d}_n\ ,
 \ \tilde{b}_n^{R}= -\i\omega b_n \formel{gleichu2}
\eeqa
und betrachten den nach \cite[Satz 27.4]{triebel}
wesentlich selbstadjungierten Besselschen Differentialoperator
\label{pagebesselop2}
\beqa
 D({\cal B}^{\nu})&:=&\{v\in L_{2,2}(0,1) \mit v(r)=r^{\nu}u(r)\ ,
 \ u\in C_\infty([0,1])\ ,\ u'(0)=v'(1)=0\} \non\\
 ({\cal B}^{\nu}v)(r)&:=&-v''(r)-\frac{v'(r)}{r}+\frac{\nu^2}{r^2}
  v(r)\mmit \nu=\omega_n^{1,0} =n-\eh\ ,\ n\in\nnn\ .\non
\eeqa
Elemente $v\in D({\cal B}^\nu)$ erf\"ullen
$rv'\in C_{1,1}$ und $v\in C_{2,1}$.
Daher erhalten wir wie oben (hier mit $(\ref{gleichu2})$ und
$\tilde{c}_n^D=0$)
\beqa
 \int_0^1 r({\cal B}^{\nu}v)(r)b_n(r)\dr
  &=&-\i\nu^{-1}\int_0^1 r ({\cal B}^{\nu}v)(r)ra_n^D(r) \dr \non\\
  &=&\i\nu^{-1}\int_0^1r(rv'(r))'a_n^D(r) \dr
     -\i\nu^{-1}\int_0^1\frac{\nu^2}{r} v(r)r a_n^D(r)  \dr\non\\
&=&\omega\nu^{-1}\int_0^1 r(rv'(r))'\tilde{a}_n(r) \dr
     -\i\nu\int_0^1 v(r) a_n^D(r)  \dr \non\\
     &=:&I_1+I_2\non\\
 I_1&=&-\i\omega\int_0^1\tilde{d}_n(r) rv'(r) \dr\non\\
  &=&\nu\omega\int_0^1\tilde{a}_n(r)v(r)  \dr
   +\i\omega\int_0^1 r\tilde{b}_n^R(r) v(r) \dr \non\\
 &=& -I_2+\omega^2\int_0^1rv(r)b_n(r)\dr \ ,\non
\eeqa
also $b_n\in D(\ol{{\cal B}^\nu})$.
Bezeichnen wir mit $\omega_{n,m}^{2,1}$ die positiven
und f\"ur $m\ra\infty$ wachsenden Nullstellen
von $J_{n-1/2}'$, so folgt aus \cite[Satz 27.4]{triebel},
$(\ref{gleichu2})$ und Lemma \ref{dgl} entweder
\beqa
  \omega^2&=& (\omega_{n,m}^{2,1})^2 \non\\
  b_n(r)&=&b_{n,m}(r)=J_{n-\eh}(\omega_{n,m}^{2,1}r) \non\\
  \tilde{a_n}(r)&=&\tilde{a}_{n,m}(r)=\i(\omega_{n,m}^{2,1})^{-1}
  a_{n,m}^D(r) \non\\
   &&=-\omega_n^{1,0}(\omega_{n,m}^{2,1})^{-1}r^{-1}b_{n,m}(r)
  \non\\
  \tilde{d}_n(r)&=&\tilde{d}_{n,m}(r)=\i(\omega_{n,m}^{2,1})^{-1}
   d_{n,m}^D(r) \non\\
   &&=\i(\omega_{n,m}^{2,1})^{-1} b_{n,m}'(r) \non
\eeqa
oder $\omega>0$ beliebig und
\beqa
  b_n=\tilde{a}_n=\tilde{d}_n=0 \ .\non
\eeqa
Herr Professor Dr. Hans Volkmer, University of Wisconsin,
Milwaukee/USA konnte w\"ahrend seines Gastaufenthaltes
in Essen zeigen, da{\ss} auch die Funktionen $J_{n-1/2}'$,
$n\in\nnn$ keine gemeinsamen Nullstellen haben:
Er ging den Weg von \cite[Seite 203]{porter},
wobei er investieren mu{\ss}te,
da{\ss} auch die Nullstellen der $J_{n-1/2}'$
trans\-zendente Zahlen sind (\cite[Seite 217]{shidlovskii}).
Damit zerfallen auch die Reihen
aus $(\ref{drittreihe})$ und $(\ref{viertreihe})$
in eine Komponente.\abstandk
Im Fall $q=2$ folgt aus $(\ref{risth})$,
da{\ss} der Raum $D_0^{2,\Gamma_2}(S)$ und
damit auch die Orthonormalsysteme verschwinden.
Wir erhalten:
\begin{satz} \formel{dummy4} \new
\begin{enumerate}
\item[i)] Bis auf normierende Konstanten gelten
\beqa
 E_{n,m}^{2,0}(r,\vp)\hspz&=&\hspz\chtau
    J_{n-\eh}(\omega_{n,m}^{2,0}r)\cos((n-\eh)\vp)
    \formel{eh1}\\
 H_{n,m}^{2,0}(r,\vp)\hspz&=&\hspz\chrho \i(\omega_{n,m}^{2,0})^{-1}
    (J_{n-\eh}(\omega_{n,m}^{2,0}r))'\cos((n-\eh)\vp)\non\\
    \hspz&-&\hspz\chtau \i(n-\eh)(\omega_{n,m}^{2,0})^{-1}
     r^{-1}J_{n-\eh}(\omega_{n,m}^{2,0}r)\sin((n-\eh)\vp)\dvp
     \formel{eh2}\\
 E_{n,m}^{2,1}(r,\vp)\hspz&=&\hspz-\chrho (n-\eh)
   (\omega_{n,m}^{2,1})^{-1}r^{-1}
    J_{n-\eh}(\omega_{n,m}^{2,1}r)\cos((n-\eh)\vp)\non\\
    \hspz&+&\hspz\chtau(\omega_{n,m}^{2,1})^{-1}
     (J_{n-\eh}(\omega_{n,m}^{2,1}r))'\sin((n-\eh)\vp)\dvp
     \formel{eh3}\\
  H_{n,m}^{2,1}(r,\vp)\hspz&=&\hspz-\chrho \i
    J_{n-\eh}(\omega_{n,m}^{2,1}r)\sin((n-\eh)\vp)\dvp
    \ .\formel{eh4}
 \eeqa
 \item[ii)] F\"ur $q=2$ sind die Orthonormalsysteme leer.
 \item[iii)] Die Eigenwerte $\omega_{n,m}^{2,0}$
 bzw. $\omega_{n,m}^{2,1}$ sind
  die positiven Nullstellen der Besselfunktion $J_{n-1/2}$
  bzw. $J_{n-1/2}'$.
\end{enumerate}
\end{satz}
Um die G\"ute der L\"osungen zu untersuchen, betrachten
wir zuerst mit der Polarkoordinatenabbildung
\beqa
 \Abb{\theta}{(0,1)\times (0,\pi)}{S}{(r,\vp)}{
 (r\cos(\vp),r\sin(\vp)}) \non
\eeqa
und $f\in C^1(S)$, $\theta^*f(r,\vp)=a(r)b(\vp)$ mit
$b\in C^1([0,\pi])$
\beqa
 \theta^*\partial_1f(r,\vp)
   &=&\cos(\vp)a'(r)b(\vp)-r^{-1}\sin(\vp)a(r)b'(\vp)
   \formel{thetad1f}\\
 \theta^*\partial_2f(r,\vp)
   &=&\sin(\vp)a'(r)b(\vp)+r^{-1}\cos(\vp)a(r)b'(\vp) \ .\non
\eeqa
Bekannt ist, da{\ss}
$f$ genau dann in $H_1(S)$ liegt, wenn
\beqa
 \int_0^{\pi}\int_0^1r|\theta^*f(r,\vp)|^2\dr\,\dvp\ ,\
 \int_0^{\pi}\int_0^1r|\theta^*\partial_if(r,\vp)|^2\dr\,\dvp
 <\infty\ ,\ i=1,2\formel{intkrit}
\eeqa
gilt.
Insbesondere folgt $(\ref{intkrit})$ aus
$a(r)\le cr^\nu $, $a'(r)\le cr^{\nu-1}$  im Falle
$\nu>0$.
Um die  Komponentenfunktionen bez\"uglich kartesischer
Koordinaten zu bestimmen, berechnen wir nach $(\ref{taudx})$
\beqa
 (r,\vp)&=&\theta^{-1}(x)=(|x|,\arccos(\frac{x_1}{|x|})) \non\\
   (\theta^{-1})^*\dr&=&
   \frac{x_1}{|x|}\dx^1+\frac{x_2}{|x|}\dx^2\non\\
  (\theta^{-1})^*\dvp
  &=& -\frac{x_2}{|x|^2}\dx^1+\frac{x_1}{|x|^2}\dx^2
     \non\\
  (\theta^{-1})^*\dr\wedge \dvp&=&\frac{1}{|x|}
  \dx^1\wedge \dx^2 \ .  \non
\eeqa
Im Falle $q=1$ erhalten wir f\"ur
$F=\chrho f_r+\chtau f_\vp\dvp$
\beqa
 (\theta^{-1})^*F (x)
  &=& f_1(x)\dx^1+f_2(x)\dx^2 \non\\
    \mmit f_1(x)&:=&(f_r\circ\theta^{-1})(x)\frac{x_1}{|x|}
     -(f_\vp\circ\theta^{-1})(x)\frac{x_2}{|x|} \non\\
     f_2(x)&:=&(f_r\circ\theta^{-1})(x)\frac{x_2}{|x|}
     +(f_\vp\circ\theta^{-1})(x)\frac{x_1}{|x|} \non\\
    (\theta^*f_1)(r,\vp)&=&
      f_r(r,\vp)\cos(\vp)-f_\vp(r,\vp)\sin(\vp) \formel{thetaf1}\\
    (\theta^*f_2)(r,\vp)&=&
      f_r(r,\vp)\sin(\vp)+f_\vp(r,\vp)\cos(\vp) \ .\non
  \eeqa
Im Fall $q=2$ gilt f\"ur $F=\chrho f_{r,\vp} \dvp$
\beqa
 (\theta^{-1})^*F&=&f_{1,2}(x)\dx^1\wedge\dx^2 \non\\
 \mmit f_{1,2}(x)&:=&(f_{r,\vp}\circ\theta^{-1})(x)\non\\
   (\theta^*f_{1,2})(r,\vp)&=& f_{r,\vp}(r,\vp) \ .\non
\eeqa
Da im Fall $q=0$ die Komponentenfunktion nicht ver\"andert wird,
k\"onnen wir uns nun den Ausdr\"ucken $(\ref{eh1})$ bis
$(\ref{eh4})$ zuwenden. Aus
\beqa
  (J_\nu (\omega r))'&=&\frac{\nu}{r}J_\nu(\omega r)
  -J_{\nu+1}(\omega r)\omega \ ,\ \omega>0\formel{jstrich}\\
 |J_\nu(r)|&\le&  r^\nu u(r) \mmit u\in C_\infty([0,1]), u(r)>0\non
\eeqa
(\cite[VII.2.(24)]{hilbert} bzw. \cite[Seite 346]{triebel})
und den Bemerkungen oben k\"onnen wir schlie{\ss}en,
da{\ss} die die Komponentenfunktionen der Transformationen
aller Ausdr\"ucke aus $(\ref{eh1})$
und
$(\ref{eh4})$ in $H_1(S)$ liegen. Das gleiche gilt im Falle
$n\ge 2$ f\"ur die Ausdr\"ucke in $(\ref{eh2})$
und  $(\ref{eh3})$.
Um eine Ausl\"oschung irregul\"arer Anteile zu
auszuschlie{\ss}en,
betrachten wir den Fall $n=1$ in $(\ref{eh2})$.
Bezeichnen wir mit $f_r$ den in der Zerlegung
$(\ref{jstrich})$ irregul\"aren Anteil des Normalenteils
und mit $f_\vp$ den
Tangentialteil, $f_1,f_2$ die Koeffizienten
bez\"uglich kartesischer Koordinaten (wie oben), so gilt
mit $\hat{c}:=1/2\cdot \i(\omega_{1,m}^{2,0})^{-1}$
\beqa
 f_r(r,\vp)&=& \hat{c}\frac{1}{r}J_{1/2}(\omega_{1,m}^{2,0}r)
  \cos(\eh\vp )\non\\
 f_\vp(r,\vp)&=& -\hat{c}\frac{1}{r}J_{1/2}(\omega_{1,m}^{2,0}r)
  \sin(\eh\vp )\ .\non
\eeqa
Aus $(\ref{thetaf1})$ und den
Additionstheoremen folgt
\beqa
 (\theta^*f_1)(r,\vp)&=&\hat{c}\frac{1}{r}J_{1/2}(\omega_{1,m}^{2,0}r)
  \cos(-\eh\vp ) \non
\eeqa
und aus $(\ref{thetad1f})$ (bis auf regul\"are Terme)
\beqa
  \theta^*\partial_1f_1(r,\vp) & \cong & \hat{c}
   ( -\frac{ J_{1/2}(\omega_{1,m}^{2,0}r) }{r^2}
   +\eh\frac{ J_{1/2}(\omega_{1,m}^{2,0}r)}{r^2} )\cos(\vp)
   \cos(-\eh\vp) \non\\
    && +\eh \hat{c}\frac{ J_{1/2}(\omega_{1,m}^{2,0}r)}{r^2}
    \sin(\vp)\sin(-\eh\vp) )\non\\
    &=& -\frac{1}{2} \hat{c}
    \frac{ J_{1/2}(\omega_{1,m}^{2,0}r)}{r^2}
     \cos(\frac{1}{2}\vp) \non\\
   &=& c
    r^{-3/2}\frac{ \sin(\omega_{1,m}^{2,0}r ) }{r}
     \cos(\frac{1}{2}\vp) \ ,\non
\eeqa
wobei wir
\beqa
 J_{1/2}(r)=c r^{1/2}\frac{\sin(r)}{r} \non
\eeqa
nach \cite[VII,(25)]{hilbert} benutzt haben.
Da hiermit die Bedingung  $(\ref{intkrit})$
verletzt ist, liegen die Transformationen der $H_{1,m}^{2,0}$
nicht in $H_1^{1}(S)$.
Analog verfahren wir mit $(\ref{eh3})$
und erhalten:
\begin{satz} \formel{dummy3} \new
Bis auf $H_{1,m}^{2,0}$ und
$E_{1,m}^{2,1}$, $m\in\nnn$ liegen die Komponentenfunktionen
der Transformationen aller
Ausdr\"ucke aus $(\ref{eh1})$ bis
$(\ref{eh4})$ in $H_1(S)$.
\end{satz}

\newpage

\subsection*{Danksagung}
\markboth{}{}
Ich bedanke mich bei Herrn Professor~Dr.~Norbert~Weck und
Herrn Professor~Dr.~Karl-Josef~Witsch
f{\"u}r zahlreiche Diskussionen und die Betreuung
meiner Dissertation.


\twocolumn
\section*{Symbole}
\addcontentsline{toc}{section}{\protect\numberline{}{Symbolverzeichnis}}
\pagestyle{myheadings}
\markboth{SYMBOLVERZEICHNIS}{SYMBOLVERZEICHNIS}
\beqa  \vspace{-1.0cm}
* &&\pageref{sterndef}\non\\
 \wedge &&\pageref{pxwedge}\non\\
 \oplus  &&\pageref{pxoplus}\non\\
 |\cdot|  &&\pageref{pxnormr}\non\\
 \subset\subset  &&\pageref{pxkompakt}\non\\
 \normu{\cdot}{H_m(S)}   &&\pageref{pxhm}\non\\
 \normu{\cdot}{H_m^q(S)}   &&\pageref{hnorm}\non\\
 \normu{\cdot}{\rvnull{-1/2}{q} (S)} ,  \normu{\cdot}{\dvnull{-1/2}{q} (S)} &&\pageref{pxrminusnorm}\non\\
 \normu{\cdot}{L_p(S)}   &&\pageref{pxlzwei}\non\\
 |\cdot|_\pm &&\pageref{pxhpmnorm}\non\\
 \spu{\cdot}{\cdot}{L_2(S)}  &&\pageref{pxlzwei}\non\\
 \spu{\cdot}{\cdot}{q,S}  &&\pageref{pxlzq}\non\\
 \spu{\cdot}{\cdot}{R^q(S)}  , \spu{\cdot}{\cdot}{D^q(S)}  &&\pageref{pxspurq}\non\\
 \spu{\cdot}{\cdot}{R^q(S)\cap D^q(S)}  &&\pageref{pxspurqdq}\non\\
 \spu{\cdot}{\cdot}{\hvnull{-m}{q} (S)}  &&\pageref{pxspuhminus}\non\\
\spu{\cdot}{\cdot}{\calL^q} &&\pageref{pxcalfq}\non\\
 \spu{\cdot}{\cdot}{_\pm}  &&\pageref{pxhpmspu}\non\\
 \spu{\cdot}{\cdot}{\mu,q,S}  &&\pageref{pxlzmuspu}\non\\
\spu{\cdot}{\cdot}{{\call}_2(\rho,R)} &&\pageref{pxspucall}\non\\
\gn &&\pageref{pxnordivspur}\non\\
\gn^{\Gamma_2} &&\pageref{pxdteilspur}\non\\
\chgn &&\pageref{pxinvdivspur}\non\\
\gt &&\pageref{pxtanrotspur}\non\\
\gt^{\Gamma_1} &&\pageref{pxteilrotspur}\non\\
\chgt &&\pageref{pxinvrotspur}\non\\
 \delta &&\pageref{pagedelta}\non\\
 \delta_{i,j} &&\pageref{pxkroneck}\non\\
  \eps,\eps_\tau,\eps^q_\tau &&\pageref{pxtrans}\non\\
\eta_s&&\pageref{lemm4}\non\\
\kappa_q,\kappa'_q &&\pageref{pxsigmaq}\non\\
\mu,\mu_\tau,\mu^q_\tau&&\pageref{pxtrans}\non\\
\rho&&\pageref{pxrhotau}\non\\
\chrho&&\pageref{pxchrhotau}\non\\
 \sigma(I) &&\pageref{pxsigma}\non\\
 \sigma_q,\sigma'_q &&\pageref{pxsigmaq}\non\\
\tau&&\pageref{pxrhotau}\non\\
\chtau&&\pageref{pxchrhotau}\non\\
 \tau_* &&\pageref{pxpushdef}\non\\
 \tau^* &&\pageref{pulldef}\non\\
 a_n&&\pageref{koeffizienten0}\non\\
 a_n^R,a_n^D &&\pageref{koeffizienten}\non\\
 A^* &&\pageref{pxadjungierter}\non\\
 A^q &&\pageref{pxaq}\non\\
 b_n&&\pageref{koeffizienten0}\non\\
 b_n^R,b_n^D &&\pageref{koeffizienten}\non\\
 c_n&&\pageref{koeffizienten0}\non\\
 c_n^R,c_n^D &&\pageref{koeffizienten}\non\\
 C_m(S) &&\pageref{pxcm}\non\\
 C_m^q(S) &&\pageref{pxcq}\non\\
 \cnull{m}{q} (S) &&\pageref{pxcqnull}\non\\
 C_m^q(\ol{S}) &&\pageref{pxcqol}\non\\
 C_\infty^{q,\Gamma}(\ol{S}) &&\pageref{pxcqgamma}\non\\
 C_R(S),C_R(S,\gamma_1,\gamma_2)&&\pageref{pxkegelspitze}\non\\
 \d &&\pageref{rotlokal}\non\\
 D &&\pageref{pxgd}\non\\
d_M&& \pageref{pxdm}\non\\
 d_n&&\pageref{koeffizienten0}\non\\
 d_n^R,d_n^D &&\pageref{koeffizienten}\non\\
 \ddiv &&\pageref{pxgdiv}\non\\
 \dq(S) &&\pageref{pxrqgamma}\non\\
 \hat{D}^{q,\Gamma_2}(S) &&\pageref{pxhrqgamma}\non\\
 D^q(S) ,  \dqn(S) &&\pageref{pxrq}\non\\
 D^q_0(S) , \dqn_0(S) &&\pageref{rqrest}\non\\
 \dq_0(S) , \hat{D}^{q,\Gamma_2}_0(S) &&\pageref{rqrest}\non\\
 \dvnull{-1/2}{q} (S) &&\pageref{pxrvnull}\non\\
D_{-1/2}^{q,\Gamma_1}(\partial S) &&\pageref{pxdteilspur}\non\\
 D(f)&&\pageref{pxdefinition}\non\\
\mbox{D--Gebiet}&&\pageref{dichtdefi}\non\\
 \div  &&\pageref{dddiv}ff,\pageref{pxdivzw},\pageref{pxdefhminus}\non\\
 \ddiv &&\pageref{pxgdiv}\non\\
 \dx^{i} &&\pageref{pxdxi}\non\\
 \dx^{I} &&\pageref{philokal}\non\\
 E_n^q &&\pageref{pxenq}\non\\
 \calf &&\pageref{calf}\non\\
\calL^q,\calL^q_R&&\pageref{pxcalfq}\non\\
 H_m(S) &&\pageref{pxhm}\non\\
 H_m^q(S) &&\pageref{pxhmq}\non\\
 \hnull m q (S) &&\pageref{pxhnullq}\non\\
H_s^{q,\Gamma_2}(\partial S) &&\pageref{pxhsqg2}\non\\
 \hvnull{-m}{q} (S) &&\pageref{pxhminus}\non\\
H_{-s}^{q,\Gamma_2}(\partial S) &&\pageref{pxhmsqg2}\non\\
 H_n^q &&\pageref{pxhnq}\non\\
 I' &&\pageref{pxistrich}\non\\
(I,J) &&\pageref{pxistrich}\non\\
 |I|,I-j,I+j &&\pageref{pxibetrag}\non\\
 \calj &&\pageref{pxcalj}\non\\
J^q&&\pageref{pxyjq}\non\\
{\call}_2(\rho,R)&&\pageref{pxcall}\non\\
 L_p(S) &&\pageref{pxlzwei}\non\\
L_{2,N}(I)&&\pageref{pagelzz}\non\\
 L_2^q(S) &&\pageref{pxlzq}\non\\
 L_{2,\mu}^q(S) &&\pageref{pxlzmu}\non\\
m&&\pageref{pxwwsph}\non\\
\hat{M}  &&\pageref{pxgm}\non\\
 \calm &&\pageref{pxcalm}\non\\
 \calm_D &&\pageref{dichtdefi}\non\\
 N &&\pageref{definspur}\non\\
 N^q&&\pageref{pxnq}\non\\
 \chn  &&\pageref{pxnspur}\non\\
\hat{R}&&\pageref{pxwwsph}\non\\
 R(f) &&\pageref{pxwerte}\non\\
 \rq(S) &&\pageref{pxrqgamma}\non\\
 \hat{R}^{q,\Gamma_1}(S) &&\pageref{pxhrqgamma}\non\\
 R^q(S) , \rqn(S) &&\pageref{pxrq}\non\\
 R^q_0(S) , \rqn_0(S)   &&\pageref{rqrest}\non\\
 \rq_0(S) , \hat{R}^{q,\Gamma_1}_0(S) &&\pageref{rqrest}\non\\
 \rvnull{-1/2}{q} (S) &&\pageref{pxrvnull}\non\\
R_{-1/2}^{q,\Gamma_2} (\partial S)&&\pageref{rmehq}\non\\
 \rot &&\pageref{rrrot}ff, \pageref{pxrotzw},\pageref{pxdefhminus}\non\\
 \rrot &&\pageref{pxgrot}\non\\
 \cals(q,N) &&\pageref{pxcals}\non\\
 S_{\mbox{{\tiny rot}}}  &&\pageref{pxspiegel}\non\\
 S_{\mbox{{\tiny div}}}  &&\pageref{pxspiegeld}\non\\
\mbox{S--Gebiet}&&\pageref{pxsgebiet}\non\\
 T &&\pageref{pxspur}\non\\
 \cht &&\pageref{pxinvspur}\non\\
\hat{T}&&\pageref{pxwwsph}\non\\
 {\cal T}M_x &&\pageref{tanm1}\non\\
 {\cal T}_x &&\pageref{pxtan1}\non\\
 U_N^{+,-,0}(R) &&\pageref{pxun}\non\\
 (V,h)  &&\pageref{pxkarte}\non\\
 \xi_k  &&\pageref{pxxik}\non\\
X&&\pageref{pxwwsph}\non\\
y_j^q&&\pageref{pxyjq}\non\\
Y^q,Y_i^q &&\pageref{pageyq}\non\\
\mbox{Z--Gebiet} &&\pageref{pxzgebiet}\non
\eeqa
\onecolumn


\newpage

\begin{thebibliography}{99}
\addcontentsline{toc}{section}{\protect\numberline{}{\bibname}}
\markboth{LITERATURVERZEICHNIS}{LITERATURVERZEICHNIS}
\bibitem{agmon} Agmon, S., {\it Lectures on
elliptic boundary value problems},
    Van Nostrand, New York, London, Toronto (1965).
\bibitem{alonso} Alonso, A. und Valli, A.,
 {\it Some Remarks on the Characterization of the Space
  of Tangential Traces of $H(rot;\Omega)$ and the Construction
  of an Extension Operator}, Manuscripta Math. {\bf 89},
   159--178 (1996).
\bibitem{bishop} Bishop, R.L. und Goldberg, S.I.,{\it Tensor
 Analysis on Manifolds}, Dover, New York (1968).
\bibitem{chillingworth} Chillingworth, D.R.J.,
{\it Differential topology with a view to applications},
  Pitman Publishing, London, San Francisco, Melbourne (1976).
\bibitem{hilbert} Courant, R., Hilbert, D.,
  {\it Methoden der Mathematischen Physik I}, Springer,
  Berlin, Heidelberg, New York (1924).
\bibitem{duff1} Duff, G. F. D., {\it Differential Forms in Manifolds
  with boundary}, Ann. of Math., {\bf 56}, 115--127 (1952).
\bibitem{duff2} Duff, G. F. D., Spencer, D. C., {\it Harmonic
 Tensors on Riemannian Manifolds with Boundary}, Ann. of. Math.,
 {\bf 56}, 128--156 (1952).
\bibitem{georgescu} Georgescu, V., {\it Some Boundary
    Value Problems for Differential Forms on
    Compact Riemannian Manifolds}, Ann. Mat. Pura Appl., {\bf 122},
    159--198 (1979).
\bibitem{hirsch} Hirsch, M.W., {\it Differential Topology},
Graduate Texts in Mathematics 33, Springer Verlag, New York (1976).
\bibitem{hu} Hu, S., {\it Homotopy Theorie}, Academic Press,
 New York, London (1959).
\bibitem{jaenich} J\"anich, K., {\it Vektoranalysis},
   Springer, Berlin, Heidelberg, New York (1992).
\bibitem{kress} Kress, R., {\it Ein kombiniertes
        Dirichlet--Neumannsches Randwertproblem bei harmonischen
        Vektorfeldern}, Arch. Rational Mech. Anal., {\bf 42},
        40--49 (1971).
\bibitem{kresspotential} Kress, R., {\it Potentialtheoretische
   Randwertprobleme bei Tensorfeldern beliebiger Dimensionen
   und beliebigen Ranges}, Arch. Rational Mech. Anal., {\bf 47},
   59--80 (1972).
\bibitem{leis} Leis, R., {\it Initial Boundary Value Problems
in Mathematical Physics}, Teubner, Stuttgart (1986).
\bibitem{martensen} Martensen, E., {\it Potentialtheorie},
 B. G. Teubner Stuttgart (1968).
\bibitem{paquet} Paquet, L., {\it Probl$\grave{{\it e}}$mes mixtes pour
 le syst$\grave{{\it e}}$me de Maxwell}, Ann. Fac. Sci. Toulouse Math.,
 {\bf 4}, 103-141 (1982).
\bibitem{picard2} Picard, R., {\it Zur Theorie der harmonischen
 Differentialformen}, Manuscr. Math. {\bf 27}, 31--45 (1979).
\bibitem{picard4} Picard, R., {\it On the boundary value
problems of electro and magnetostatics},
       Proc. R. Soc. Edinburgh, {\bf 92A}, 165--174 (1982).
\bibitem{picardhodge} Picard, R., {\it Ein Hodge--Satz f\"ur
  Mannigfaltigkeiten mit nicht--glattem Rand}, Math. Meth. in the
  Appl. Sci. {\bf 5}, 153--161 (1983).
\bibitem{picard3} Picard, R., {\it An Elementary Proof
for a Compact Imbedding Result in Generalized Electromagnetic
Theory}, Math. Z. {\bf 187}, 151--164 (1994).
\bibitem{wwp} Picard, R., Weck, N., Witsch, K.J.,
 {\it Time--Harmonic Maxwell Equations in the Exterior of
 Perfectly Conducting, Irregular Obstacles},
 eingereicht bei SIAM J. Math. Anal. (1999).
\bibitem{porter} Porter, M. B., {\it On the Roots of the
 Hypergeometric and Bessel's Function}, Amer. Journal of Math.,
 {\bf 20}, 193 -- 214 (1898).
\bibitem{rosenberg} Rosenberg, S., {\it The Laplacian on a Riemannian
 Manifold}, Cambridge University Press (1997).
\bibitem{saranen} Saranen, J.,
 {\"Uber das Verhalten der L\"osungen der Maxwellschen
 Randwertaufgabe in Gebieten mit Kegelspitzen}, Math. Meth. Appl.
 Sci., {\bf 2}, 235--250 (1980).
\bibitem{schwarz} Schwarz, G., {\it Hodge Decomposition - A Method
  For Solving Boundary Value Problems}, Lecture Notes Math., Vol. 1607,
  Springer Berlin (1995).
\bibitem{shidlovskii} Shidlovskii, A. B., {\it Transcendental Numbers},
 Walter de Gruyter, Berlin, New York (1989).
\bibitem{triebel} Triebel, H., {\it H\"ohere Analysis}, Harri Deutsch,
  Thun und Frankfurt am Main (1972).
\bibitem{watson} Watson, G. N., {\it A Treatise On The Theory
 of Bessel Functions}, Cambridge University Press (1922).
\bibitem{weber2} Weber, C., {\it A local compactness theorem
  for Maxwell's equations}, Math. Meth. Appl. Sci. {\bf 2},
  12--25 (1980).
\bibitem{weber} Weber, C., {Regularity Theorems for
 Maxwell's Equations}, Math. Meth. in the Appl. Sci. {\bf 3},
 523--536 (1981).
\bibitem{weck2} Weck, N., {\it Eine L\"osungstheorie f\"ur die
    Maxwellschen Gleichungen auf Riemannschen Mannigfaltigkeiten
    mit nicht glattem Rand}, Habilitationsschrift,
    Universit\"at Bonn (1972).
\bibitem{weck} Weck, N., {\it Maxwell's Boundary Value
 Problem on Riemannian Manifolds with Nonsmooth Boundaries},
  J. Math. Anal. Appl., {\bf 46}, 410--437 (1974).
\bibitem{wwspherical} Weck, N. und Witsch, K.J., {\it Generalized
    Spherical Harmonics and Exterior Differentiation in Weighted Sobolev
    Spaces}, Math. Meth. Appl. Sci., {\bf 17}, 1017--1043 (1994).
\bibitem{wentzig} Wentzig, S., {\it Eigenwerte des Maxwelloperators
   im wesentlichen Spektrum}, Dissertation, Essen (1995).
\bibitem{weyl} Weyl, H., {\it Die nat\"urlichen
 Randwertaufgaben im Au{\ss}enraum f\"ur Strahlungsfelder
 beliebiger Dimension und beliebigen Ranges},
 Math. Z., {\bf 56}, 105--119 (1952).
\bibitem{witsch} Witsch, K.J., {\it A
  Remark on a Compactness Result in Electromagnetic Theory},
  Math. Meth. Appl. Sci., {\bf 16}, 123--129 (1993).
\bibitem{wloka} Wloka, J., {\it Partielle Differentialgleichungen},
  Teubner, Stuttgart (1982).
\end{thebibliography}
\end{document}